\documentclass[11pt]{article}
\usepackage{amsmath}
\usepackage{amssymb}
\usepackage{amsfonts}
\usepackage{mathrsfs}
\usepackage{amsthm}
\usepackage{epsfig}
\usepackage{amsmath, amsthm, amssymb}
\usepackage[english]{babel}
\usepackage{graphics}
\usepackage{graphicx}
\usepackage{subfigure}
\usepackage{eucal}
\usepackage{amscd}
\usepackage{array}
\usepackage{verbatim}

\usepackage{multirow}

\newtheorem{thm}{Theorem}[subsection]

\newtheorem{lemma}[thm]{Lemma}

\newtheorem{prop}[thm]{Proposition}

\newtheorem{cor}[thm]{Corollary}
{\theoremstyle{definition}
\newtheorem{defn}[thm]{Definition}}

{\theoremstyle{remark}
\newtheorem{rem}[thm]{Remark}}

{\theoremstyle{remark}
\newtheorem{ex}[thm]{Example}}

\newenvironment{proofof}[2]{\begin{proof}[Proof of #1 \ref{#2}.]}{\end{proof}}
\newcommand{\be}{\begin{equation}}
\newcommand{\ee}{\end{equation}}
\newcommand{\bes}{\begin{equation*}}
\newcommand{\ees}{\end{equation*}}

\newcommand{\trajcdot}{}
\newcommand{\octcdot}{}

\newcommand{\Sing}{\mathcal{V}}

\newcommand{\LL}{{\bf  L}}
\newcommand{\A}{{\bf  A}}
\newcommand{\B}{{\bf{B}}}
\newcommand{\C}{{\bf{C}}}
\newcommand{\D}{{\bf  D}}

\newcommand{\g}[3]{\mathfrak{g}^{#2}_{#1} { #3  }}
\newcommand{\gop}[2]{\mathfrak{g}^{#2}_{#1}}
\newcommand{\PhiTraj}[2]{\Phi_{#1}  #2 }

\newcommand{\Die}{\mathrm{D}}
\newcommand{\Cdot}{}

\newcommand{\DD}{ {T_1 \mathbb{D}}}
\newcommand{\R}{\mathbb{R}}
\newcommand{\CC}{\mathbb{C}}
\newcommand{\Q}{\mathbb{Q}}
\newcommand{\Aff}{Af\!f}
\newcommand{\RP}{\mathbb{R}\mathbb{P}^1}
\newcommand{\CP}{\mathbb{C}\mathbb{P}^1}
\newcommand{\Disk}{\mathbb{D}}

\title{Symbolic coding for linear trajectories in the regular octagon}

\author{John Smillie\footnote{Work partially supported by NSF Grant DMS $0601299$.}
  \and  Corinna
Ulcigrai\footnote{Work partially supported by an RCUK Fellowship.} }
\date{}

\begin{document}

\maketitle

\begin{abstract}
We consider a symbolic coding of linear trajectories in the regular octagon with opposite sides identified (and more generally in regular 2n-gons).  Each infinite trajectory gives a cutting sequence corresponding to the sequence of sides hit. We give an explicit characterization of these cutting sequences. The cutting sequences for the square are the well studied Sturmian sequences which can be analyzed in terms of the continued fraction expansion of the slope. We introduce an analogous continued fraction algorithm which we use to connect the cutting sequence of a trajectory with its slope. Our continued fraction expansion of the slope gives an explicit sequence of substitution operations which generate the cutting sequences of trajectories with that slope. Our algorithm can be understood in terms of renormalization of the octagon translation surface by elements of the Veech group.
\end{abstract}

\section{Introduction}

In this paper we give  a complete characterization and an explicit description of the symbolic sequences which arise in coding bi-infinite linear trajectories on the regular octagon. The basic definitions and questions addressed are explained in \S\ref{questionssec}. The corresponding question in the case of the square gives rise to Sturmian sequences. These were considered by \cite{C:obs} and \cite{S:note} in the 1870's, by Morse and Hedlund \cite{MH:sym} in 1940 and by many authors since then (see \cite{Ar:stur} for a contemporary account). Sturmian sequences are interesting because of their geometric origin but are also of interest because they give the simplest non-periodic infinite sequences (see \cite{CH:seq}). The idea that the regular octagon and the square share some special property that might make their analysis easier first appears in the work of Veech \cite{Ve:tei}. Veech showed that the regular octagon and the square are both examples of lattice polygons (see section  \S \ref{defsec}).

The relation between the direction of a trajectory on the square and the symbol sequence is made by means of continued fractions (see e.g.~\cite{Se:mod}, \cite{Ra:low}).   We will make the connection between symbol sequences and directions of trajectories in the case of the octagon using an appropriate version of the continued fraction algorithm.
It is a classical fact  that continued fractions are related to the  Teichm\"uller  geodesic flow on the modular surface \cite{AF:mod, Se:mod}. 
 Our continued fraction algorithm is closely connected to the Teichm{\"u}ller geodesic flow on an appropriate Teichm\"uller curve that plays the role of the modular surface. 
 A type of continued fraction algorithm for the octagon was introduced by Arnoux and Hubert (\cite{AH:fra}). Our algorithm is different from that of Arnoux and Hubert in its details though similar in spirit. The dynamics of our continued fraction algorithm is also connected to the coding of this type of geodesic flow introduced by Caroline Series \cite{Se:sym}.

As Sturmian sequences also code orbits of rotations, the symbolic sequences that we consider give also a coding of a class  of interval exchange transformations. While multi-dimensional continued fractions  have been successfully used to analyze the space of all interval exchange transformations since the work of Rauzy and Veech \cite{Ra:ech, Ve:gau}, the challenge in the case of the octagon is that the corresponding interval exchange transformations constitute a family of measure zero. Our renormalization scheme allows us to obtain a characterization of the symbolic sequences coding  this family. 

We begin the paper with a discussion of the case of the square and Sturmian sequences, see \S\ref{torussec}.  This will serve to introduce key ideas and methods that we will use in the case of the octagon. Our treatment of this case is similar to that of Series \cite{Se:mod} but has the novel feature that it does not make use of the fact that the universal cover of the torus is $\R^2$. 
We formulate our main results for the octagon in \S\ref{mainresultssec}. Their proofs are presented in the following sections.
We end the paper with an analysis of the coding of trajectories for all regular polygons with an even number of sides, \S\ref{2ngonssec}. Veech showed that these are also lattice polygons.

\subsection{Cutting sequences: definition and questions}\label{questionssec}
\begin{figure}
\centering
{
\includegraphics[width=0.35\textwidth]{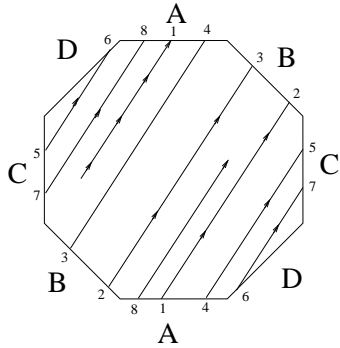}}
\caption{A linear trajectory in the octagon.\label{octagonflow}}
\end{figure}

\paragraph{Linear trajectories in the octagon.}
Let ${O}$ be a regular octagon. The boundary of $O$ consists of four pairs of parallel sides. We define a correspondence between points in opposite sides by using the   isometry between them which is the restriction of a translation. We define a \emph{linear trajectory} in direction $\theta$ to be a path which starts in the interior of the octagon and moves with constant velocity vector making an angle $\theta$ with the horizontal, until it hits the boundary at which time it re-enters the octagon at the corresponding point on the opposite side and continues traveling with the same velocity.  
For an example of a trajectory see Figure \ref{octagonflow}.
In this paper we will restrict ourselves to \emph{bi-infinite} trajectories, i.e.~trajectories which do not hit vertices of the octagon and thus are well defined as time tends to $ \pm\infty$. 
 In a forthcoming paper we will deal with the case of saddle connections. 

\paragraph{Cutting sequences.}
We now describe a \emph{symbolic coding} for trajectories.
 Let us label  each pair of opposite  sides of the octagon $O$ with a letter of an alphabet $\mathscr{A}=\{A,B,C,D\}$ as in Figure \ref{octagonflow}.
\begin{defn}\label{cuttingsequencedef}
The \emph{cutting sequence} $c(\tau)$ associated to the linear trajectory $\tau$ is the bi-infinite word  in the letters of the alphabet $\mathscr{A}$, which is obtained by reading off the labels of the pairs of identified sides crossed by the trajectory $\tau$ as time increases. 
\end{defn}
\begin{ex}
Consider the  piece of the trajectory $\tau_0$ shown in Figure \ref{octagonflow}. The sequence of sides hit during this orbit piece shows that    $c(\tau_0)$ contains $\ldots ABBACDCA \ldots$. 
\end{ex}
Let us remark that, since we are considering systems of entropy zero, the cutting sequences have very low complexity. More precisely, the number of subwords of length $n$ is at most  $3n+1$ (see Proposition \ref{growth} in \S\ref{wordcomplexitysec}).

\paragraph{Questions addressed.}
We are interested in describing explicitly the sequences that can arise as cutting sequences of linear trajectories in the octagon.
We list here some specific questions. 
Here and in the rest of the paper, $w$ denotes a word in the  alphabet $\mathscr{A}$, usually infinite unless differently specified.
\begin{enumerate}
\item[Q1.] Given a word $w$, 
 can we determine whether $w$ occurs as a cutting sequence? Can one give a \emph{characterization} of the infinite words  which arise as cutting sequences of linear trajectories in the octagon?

\item[Q2.] Does  a cutting sequence $c(\tau)$ uniquely determine the direction of the trajectory $\tau$? 
Or, given a finite subword of  $c(\tau)$, can one identify a sector of possible directions?

\item[Q3.] Given a finite word  $w$ which can appear in an infinite cutting sequence, can one give an explicit algorithm to construct a trajectory $\tau$ such that $w$ is a subword of $c(\tau)$?

\item[Q4.] Given a fixed direction, can one give an explicit algorithm to generate all possible words that occur as cutting sequences of trajectories with that direction? Or, more generally, can one algorithmically generate all cutting sequences corresponding to a sector  of directions?

\end{enumerate}
In this paper we give an explicit answer to these questions. It turns out that in the case of the octagon as in the case of the square the collection of cutting sequences is not closed in the space of all words with its usual topology. In  Theorem \ref{cuttingseqthm} we  give a \emph{characterization} of the \emph{closure} of the set of {cutting sequences} of trajectories. 
An algorithmic way to determine the direction(s) in Q$2$ is given by Theorem \ref{directionsthm}. The construction of a trajectories as in Q$3$ is given by Proposition \ref{finiterealization} and is used in the characterization of the closure.  An algorithm to generate cutting sequence as in Q$4$ is described in \S\ref{algorithmsec}. 

\subsection{Cutting sequences for the square}\label{torussec}
In order to motivate our approach to the octagon and give a warm-up example to the reader, let us recall here the characterization of cutting sequences for the torus. We follow the main ideas in Series \cite{Se:geo}, but sketch a proof in the spirit the proof of Proposition \ref{derivationiscuttingseqprop} for the octagon. 

\begin{figure}
\centering
\subfigure[$0\leq \theta \leq \pi/2$]{\label{square0}
\includegraphics[width=0.22\textwidth]{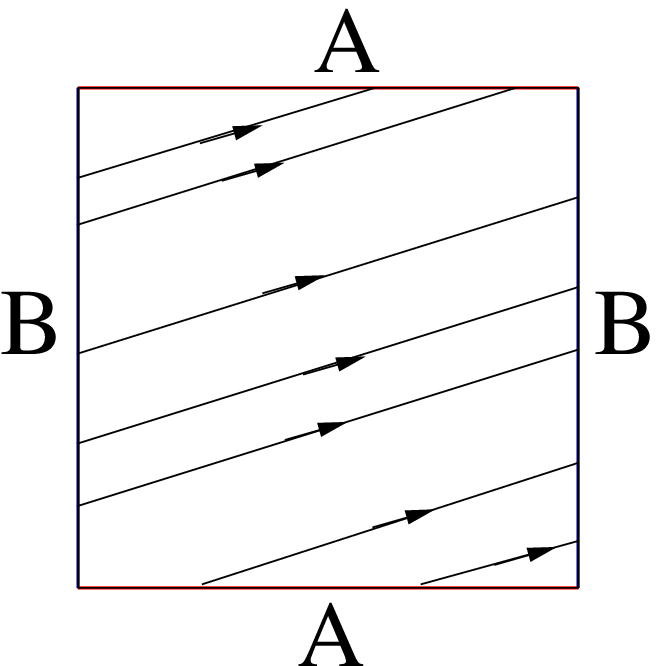}}
\hspace{2mm}
\subfigure[$\mathscr{D}'_0$ ]{\label{D0square}
\includegraphics[width=0.18\textwidth]{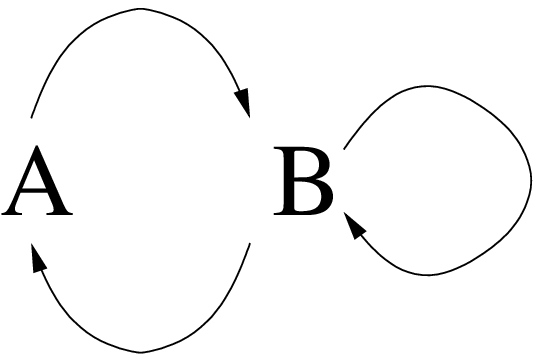}}
\hspace{2mm}
\subfigure[$\pi/2 \leq \theta \leq \pi$ ]{\label{square1}
\includegraphics[width=0.22\textwidth]{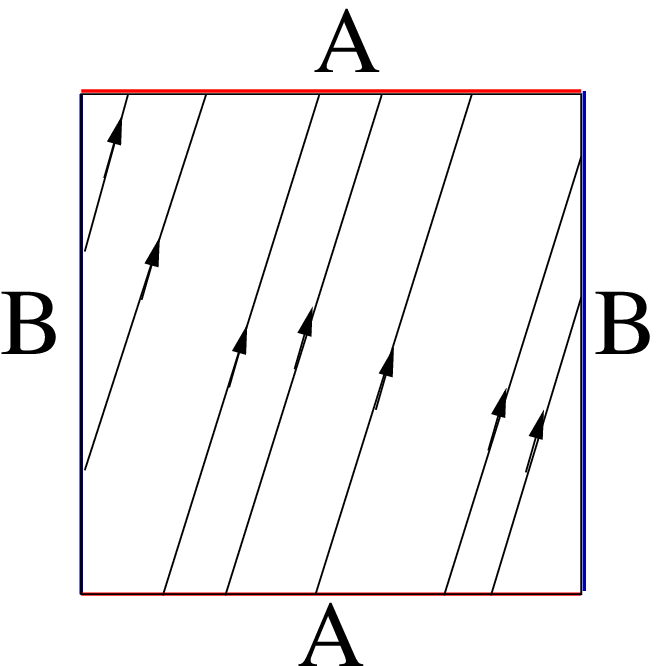}}
\hspace{2mm}
\subfigure[$\mathscr{D}'_1$ ]{\label{D1square}
\includegraphics[width=0.18\textwidth]{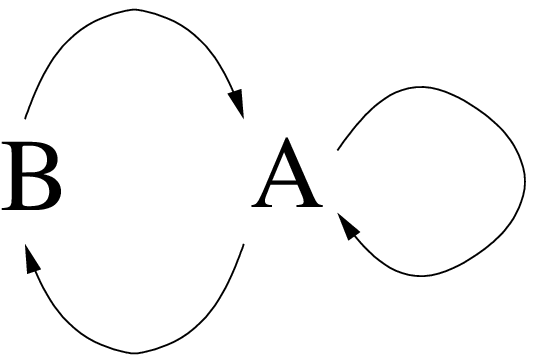}}
\caption{Possible transitions in the square.}
\end{figure}

Consider a square with opposite sides identified; as in Figure \ref{square0}, label by $A$ and $B$ respectively its horizontal and vertical sides\footnote{ Since squares (or, more generally, parallelograms) tile the plane by translation, the cutting sequence of a trajectory in a  square (parallelogram) is the same than the cutting sequence of a straight line in $\mathbb{R}^2$ with respect to a square (or affine) grid.}.  
We remark first that the group $\Die_4$ of isometries of a square acts on trajectories as well as acting on cutting sequences by permuting the letters;  the reflections in the vertical and horizontal  axes of the square preserve the labelling, while the reflection in the diagonals interchange $A$ and $B$.  The map from trajectories to cutting sequences is equivariant with respect to these actions.
Hence, given a cutting sequence $c(\tau)$, we we can assume without loss of generality that  $\tau$ is a trajectory with a direction in the first quadrant. Let us parametrize directions by measuring the angle that they make with the horizontal vector $e_1$. 
Thus, we can assume we have a trajectory in direction $\theta$ with $0\leq \theta < \pi/2$. Moreover the case $\pi/4 \leq \theta \leq \pi/2$ can be reduced to  $0\leq \theta < \pi/4$ by interchanging the role of $A$ and $B$.

Clearly if $0\leq \theta \leq \pi/4$, as in Figure \ref{square0}, the cutting sequence does not contain the subword $AA$  and if $\pi/4 \leq \theta \leq \pi/2$, as in Figure \ref{square1}, it does not contain the subword $BB$. Hence the only pairs of consecutive letters, or \emph{transitions}, that can occur are the ones shown in the diagram $\mathscr{D}'_0$ in Figure \ref{D0square} or respectively in the diagram $\mathscr{D}'_1$ in Figure \ref{D1square}. Let us say that a word $w \in \{A,B\}^{\mathbb{Z}}$ is \emph{admissible} if it corresponds to an infinite path on either the graph $\mathscr{D}'_0$ or the graph $\mathscr{D}'_1$.

Given an admissible word $w$, denote by $w'$ the \emph{derived sequence}\footnote{In this section, we are using the terminology from Series \cite{Se:geo}. } obtained by erasing one $B$ (respectively one $A$) from each block of consecutive $B$'s if  $w$ has no transitions $AA$ (respectively $BB$). 
\begin{ex}\label{sigmaex} For example, consider the following word $w$ and its derived sequence $w'$:
\begin{eqnarray}
w &=& \dots ABBBABBBBABBBABBBA BBBBA\dots  \nonumber \\ 
w'& = & \dots A BB\phantom{B}A BBB\phantom{B}ABB\phantom{B}ABB\phantom{B}ABBB\phantom{B}A \dots \nonumber
\end{eqnarray}
\end{ex}
\noindent A word is \emph{infinitely derivable} if it is admissible and each of its derived sequences is admissible. 

\begin{prop}[\cite{Se:geo}]\label{toruscutseq}
Cutting sequences of linear trajectories  on the square are infinitely derivable. 
\end{prop}

\noindent Moreover, the converse of Proposition \ref{toruscutseq} is \emph{almost} true; the exceptions, i.e. words in $\{A, B\} ^{\mathbb{Z}}$ which are infinitely derivable and are not cutting sequences such as $\dots BBBBBABBBBB\dots $, can be explicitly described. 
The space of words has a natural topology which makes it a compact space. The word given above is not a cutting sequence, but it has the property that any finite subword can be realized by a finite trajectory. This is equivalent to saying that it is in the closure of the space of cutting sequences. The closure of the space of cutting sequences is precisely the set of infinitely derivable sequences.

\begin{proofof}{Proposition}{toruscutseq} Let us prove that the derived sequence $w'$ of a cutting sequence $w=c(\tau)$ is again a cutting sequence of a linear trajectory on the square. Since a cutting sequence is admissible, this will suffice to show that every derived sequence is again admissible. 

By applying an element of $\Die_4$  we can assume that $0\leq \theta \leq \pi/4$. 
Let us add to the square the diagonal in Figure \ref{diagonalfig} labelled by the letter $c$.
We obtain an \emph{augmented cutting sequence} $\tilde{c}(\tau) \in \{A,B,c \}^{\mathbb{Z}}$ by recording the crossings of the sides of the square and the diagonal. We see  that   $c$ is crossed only during a $BB$ transition, not during a $BA$ or $AB$ transition, as in Figure \ref{diagonalfig}. 
Thus the auxiliary cutting sequence can be determined from cutting sequence of our trajectory without reference to the trajectory itself. We  construct it by adding a $c$ between any pair of $B$'s. 
\begin{ex}
The cutting sequence $w=c(\tau)$ in  Example \ref{sigmaex} give rise to the augmented sequence 
\bes
\dots ABcBcBABcBcBcBABcBcBABcBcBA BcBcBcBA \dots 
\ees
\end{ex}

\begin{figure}
\centering
\subfigure[Auxiliary diagonal]{\label{diagonalfig}
\includegraphics[width=0.215\textwidth]{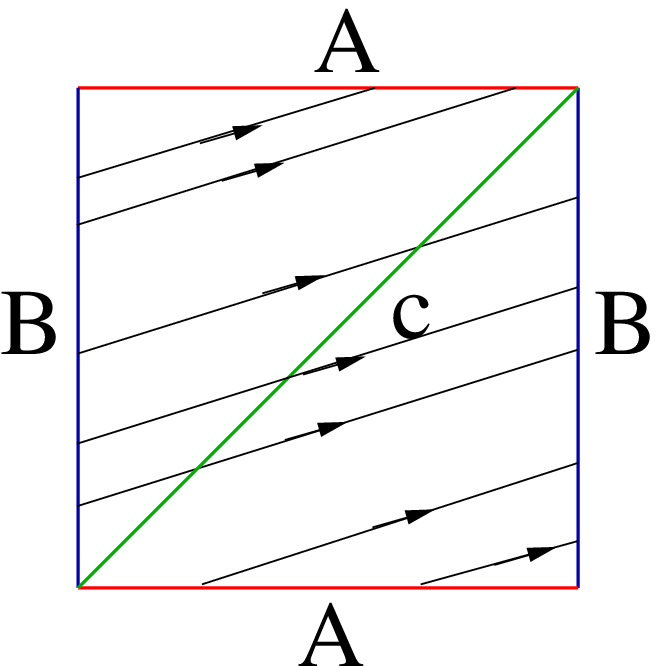}}
\hspace{2mm}
\subfigure[Parallelogram $\Pi$]{\label{pifig}
\includegraphics[width=0.33\textwidth]{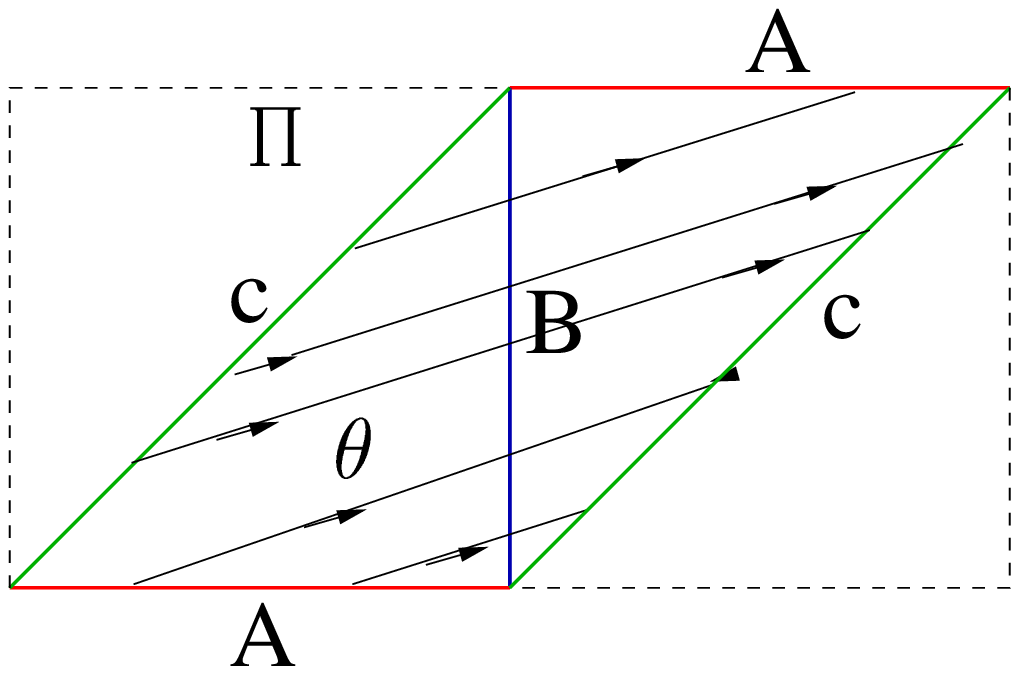}}
\hspace{2mm}
\subfigure[Renormalized flow]{\label{shearedfig}
\includegraphics[width=0.33\textwidth]{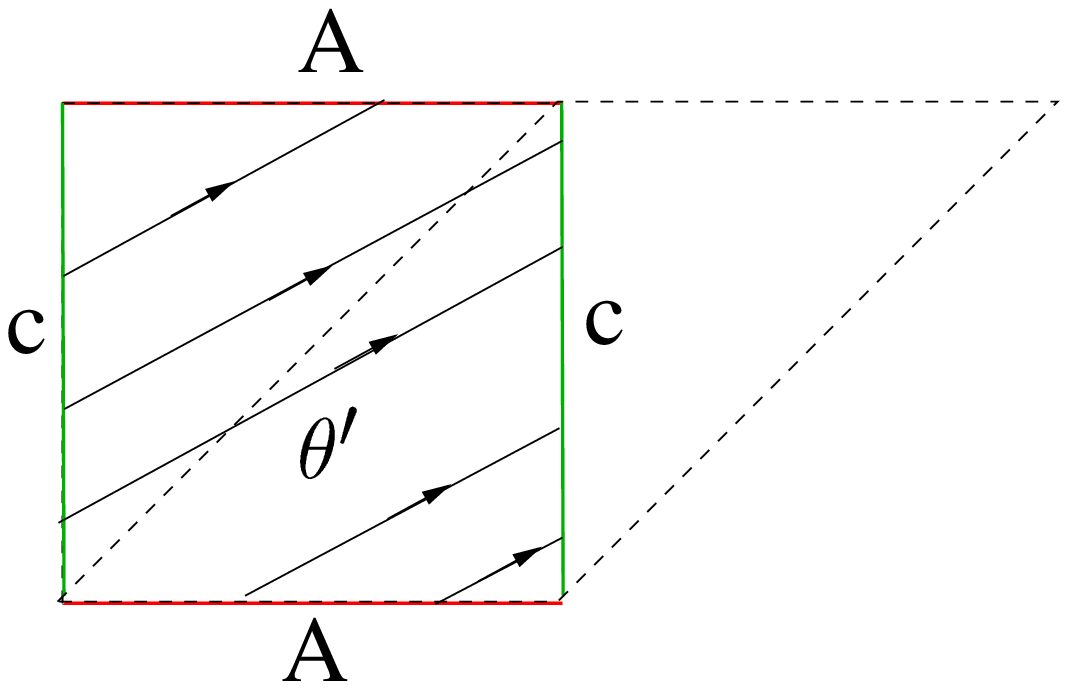}}
\hspace{2mm}
\caption{Geometric renormalization of cutting sequences of the square,  $0\leq \theta \leq \pi/2$.}
\label{cutandpaste}
\end{figure}

Let us cut the square along edge $c$ and reglue along edge $B$ to obtain  the parallelogram $\Pi$ in Figure \ref{pifig} whose sides are labelled by $A$ and $c$. The sequence obtained by erasing $B$ in the augmented sequence $\tilde{c}(\tau)$ is clearly the cutting sequence of $\Pi$. Let us \emph{renormalize} $\Pi$ by applying the linear map which sends it to a square. If we assume that the origin of the plane is at the lower left vertex of the square, this transformation is a shear, i.e.~the linear transformation given by  the matrix $\left( \begin{smallmatrix} 1& -1 \\ 0 & 1 \\ \end{smallmatrix} \right)$. Since the map is affine, linear trajectories in $\Pi$ are sent to linear trajectories (in a new direction) on a square (see Figure \ref{shearedfig}). Thus, if we replace the $c$'s with $B$'s, the sequence that we obtain is a cutting sequence of a trajectory  $\tau'$ in a new  direction $\theta'$. Moreover, the final effect  of the combinatorial operation, i.e.~replacing $BB$ by $BcB$ first and then dropping the $B$'s and replacing the $c$'s by $B$'s, is exactly the same as that of deriving the sequence. Thus, the new cutting $c(\tau')$ sequence coincides with the derived sequence $w'$.  
\end{proofof}
The crucial step in the previous proof is to show that the combinatorial operation of derivation of cutting sequences corresponds to a geometric operation of renormalization of trajectories. 

Let us consider the  map which associates to the direction $0\leq \theta\leq \pi/2$ of a trajectory $\tau$ the direction $\theta' $ of the renormalized trajectory $\tau'$ constructed the proof of Proposition \ref{toruscutseq}, which has the property that  $c(\tau')= c(\tau)'$. This map, that we call $F$, is defined on the chart of  $\RP$ which correspond to the first quadrant. If,  instead than the angle coordinate $\theta$, one can chooses for this chart the  coordinate $t$ obtained projecting radially  to the line $\{(1-t,t), \ 0\leq t \leq 1\}$, the map $F$ in this new coordinate is the   \emph{Farey map}, i.e.
\be\label{classicalFarey}
F(t) = \begin{cases}\frac{t}{1-t}  &\mathrm{if}\ 0\leq t \leq \frac{1 }{2} \\  \frac{1-t}{t}  &\mathrm{if}\ \frac{1}{2} \leq t \leq 1 \end{cases} , \qquad \mathrm{where} \quad t(\theta) = \frac{ \sin \theta }{\cos \theta + \sin \theta} . 
\ee
In this way, one can see that the number of $B$'s between two consecutive $A$'s is related to the number of times the iterates of $0< t(\theta)< 1/2$ under the Farey map lie in the sector $[0, 1/2]$. If we form a sequence $\{s_k\}_{k\in \mathbb{N}}$ where $s_k$ is $0$ when the $k$th iterate of $t(\theta)$ lies in $[0, 1/2]$ and $1$ otherwise then we have created the Farey expansion of $\theta$. We refer to this as the additive continued fraction expansion of the slope of $\theta$. The usual (mulitplicative) continued fraction expansion is obtained by counting the number of symbols in each block of $0$'s or $1$'s. Thus, one can recover the continued fraction expansion of the slope of $\theta$ from the cutting sequence of $\tau$. 

\section{Definitions and main results}\label{mainresultssec}
\subsection{Cutting sequences and derivation}
\begin{figure}[!h]
\centering
{\includegraphics[width=0.3\textwidth]{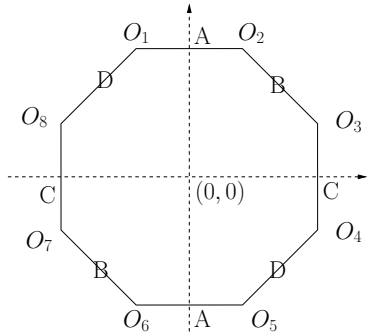}}
\caption{\label{octagonbasicfig}The labeled octagon  $O \subset \mathbb{R}^2$.}
\end{figure}

Let us begin, as in \S \ref{torussec}, with some preliminary observations about the action of the isometries of the octagon on cutting sequences.
It is useful to identify the octagon $O$  
with an explicit subset of $ \mathbb{R}^2$ as in Figure \ref{octagonbasicfig}, with center at the origin, sides of unit length and two horizontal sides. Let  $O_1, \dots, O_8$ be the vertices labelled clockwise,   $\overline{O_i O_{i+1}}$ be the side connecting $O_i$ with $O_{i+1}$ and  let $\overline{O_{1} O_{2}}$ be the horizontal top side. 
We denote by $S_O$ be the surface obtained from ${O}$ by identifying opposite sides by parallel translation, i.e. by identifying the side $\overline{O_i O_{i+1}}$ to $\overline{O_{i+4} O_{i+5}}$ for $i=1,2,3,4$.
We use the letters $A$, $B$, $C$ and $D$ to label the pair of sides $\overline{O_1 O_2}= \overline{O_6 O_5}$, $\overline{O_2 O_3}=\overline{O_7 O_6}$, $\overline{O_3 O_4}=\overline{O_8 O_7}$ and $\overline{O_4 O_5} = \overline{O_1 O_8 }$ in $S_O$, as shown in Figure \ref{octagonbasicfig}. 

\subsubsection{Action of the isometry group of the octagon on trajectories}\label{isometriesoctagonsec}
The group  $\Die_8$  of isometries of the octagon acts on the octagon preserving the pairing of sides so it acts on $S_O$ sending linear trajectories to linear trajectories and acts on cutting sequences by permutations of the letters as follows. Given $\nu \in \Die_8$, $\nu$ sends the trajectory $\tau:\R\to S_O$ to the trajectory given by the composition $\nu \tau:\R\to S_O$ which we denote by $\nu\trajcdot\tau$. 
\begin{defn}\label{inducedpermutation}
Let $\pi(\nu)$ be the permutation of the labels  $\{A,B,C,D\} $ defined by  $\pi(L)=L'$ iff $\nu$ maps a side labelled by $L$ to a side labelled by $L'$). We call $\pi(\nu)$ the \emph{permutation induced by} $\nu$.  
\end{defn}
\noindent If $\pi= \pi(\nu)$ then the cutting sequence $c(\nu\trajcdot \tau)$ is obtained from $w=c(\tau)$ by permuting its  letters according to $\pi$.   The action of $\pi$ on (finite or infinite) words $w$  will be denoted by $\pi \cdot w$. 

Since the rotation $\rho_{\pi}$ by $180^\circ$ is in $\Die_8$ and preserves the labels of sides of $O$ (i.e. the induced permutation $\pi(\rho_{\pi})$ is the identity), we can always assume without relabelling letters that a given cutting sequence $c(\tau)$ corresponds to a trajectory $\tau$ in direction $\theta\in [0,\pi]$.   
Moreover, up to relabelling of the letters, we can reduce to the case of a trajectory in direction $\theta$ with $\theta\in [0,\pi/8]$ as follows. 

We will denote by $\Sigma_i : = [i \pi /8 ,  (i+1)\pi/8)$ for $i=0,\dots, 6$ and $\Sigma_i : = [7 \pi /8 ,  \pi]$,  the eight sectors of length $\pi/8$ in $[0,\pi]$.  The set ${\Sigma}_0$ is a  fundamental domain for the action of $\Die_8$ on the circle of directions.  Let $\overline{\Sigma}_i : = [i \pi /8 ,  (i+1)\pi/8]$ denote the corresponding closed intervals. 
For each  $0\leq i \leq 7$  let $\nu_i\in\Die_8$  be the isometry which sends $\overline{\Sigma}_i$ to $\overline{\Sigma}_0$  and is given by:
\be \label{nujdef}
\begin{split} &
\nu_0= \begin{pmatrix} 1 & 0 \\ 0 & 1 \end{pmatrix} \phantom{-} \, \, 
\nu_1= \begin{pmatrix}  \frac{1}{\sqrt{2}} &  \frac{1}{\sqrt{2}} \\  \frac{1}{\sqrt{2}} & - \frac{1}{\sqrt{2}} \end{pmatrix}\, \, 
\nu_2= \begin{pmatrix} \frac{1}{\sqrt{2}} &  \frac{1}{\sqrt{2}} \\ - \frac{1}{\sqrt{2}} &  \frac{1}{\sqrt{2}} \end{pmatrix}  \, \, 
\phantom{-} 
\nu_3=\begin{pmatrix}0& 1 \\ 1& 0  \end{pmatrix} \\ & 
\nu_4=\begin{pmatrix} 0 & 1 \\ -1 & 0 \end{pmatrix} \, \, 
\nu_5= \begin{pmatrix} - \frac{1}{\sqrt{2}} &  \frac{1}{\sqrt{2}} \\  \frac{1}{\sqrt{2}} &  \frac{1}{\sqrt{2}} 
\end{pmatrix}\, \, 
\nu_6=  \begin{pmatrix} - \frac{1}{\sqrt{2}} &  \frac{1}{\sqrt{2}} \\ - \frac{1}{\sqrt{2}} & - \frac{1}{\sqrt{2}} \end{pmatrix}\, \, 
\nu_7= \begin{pmatrix} -1 & 0 \\ 0 & 1 \end{pmatrix}.
\end{split}
\ee
\noindent Since we assume $\theta \in  [0,\pi]$, $\theta$ can be put in $\overline{\Sigma}_0$ by applying one of the matrices $\nu_i$, $i=0, \dots, 7$. 
\begin{defn}\label{normalizedtrajdef}
Given a trajectory $\tau$ with direction in $\Sigma_k$, we call \emph{normal form} of the trajectory $\tau$ the trajectory $n(\tau):= \nu_k\trajcdot \tau$, obtained   by applying the isometry $\nu_k$ of $\Die_8$ which maps $\overline{\Sigma}_k$ to $\overline{\Sigma}_0$. 
\end{defn}

The corresponding eight permutations of $\{A,B,C,D\} $  will be denoted by $\pi_i = \pi (\nu_i)$ for $i=0,\dots, 7$ and are given explicitly by
\bes
\begin{array}{llll}
\pi_0= Id & \pi_1 = (AD)(BC) & \pi_2 = (ABCD) & \pi_3= (AC)\\
\pi_4=(AC)(BD) & \pi_5 = (AB)(CD) & \pi_6 = (ADCB) & \pi_7 = (BD).
\end{array}
\ees
Since  the action of permutations on cutting sequences is induced by the action of the isometries, we have the relation  
\be \label{relabellingeq}
c(\nu_k \trajcdot \tau) = \pi_k \cdot c(\tau), \qquad  \forall \quad 0\leq k \leq 7 .
\ee

\subsubsection{Transition diagrams}
As in the case of the torus, let us first determine the possible two letter words, or \emph{transitions}, which can occur in cutting sequences of  trajectories with directions in a specified sector.
\begin{figure}
\centering
\subfigure[ Starting with A]{\label{D0transitionsA}
\includegraphics[width=0.23\textwidth]{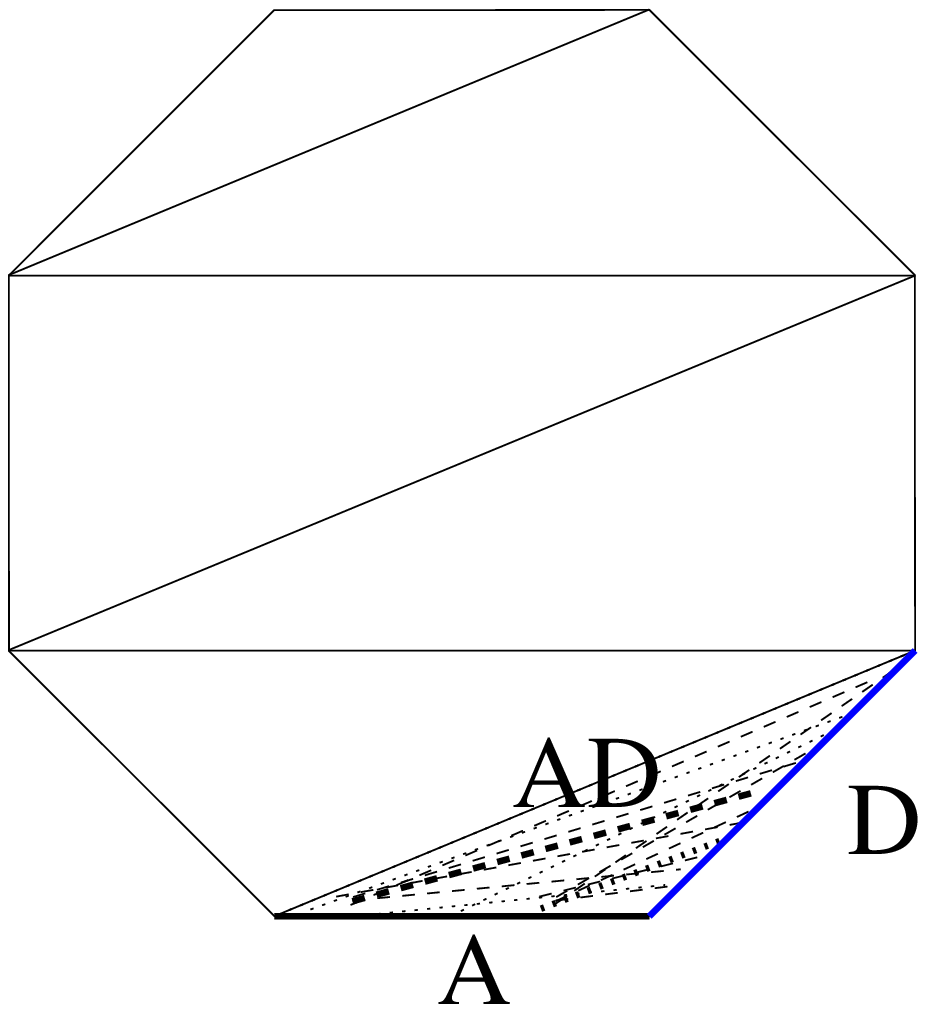}}
\subfigure[ Starting with B]{\label{D0transitionsB}
\includegraphics[width=0.23\textwidth]{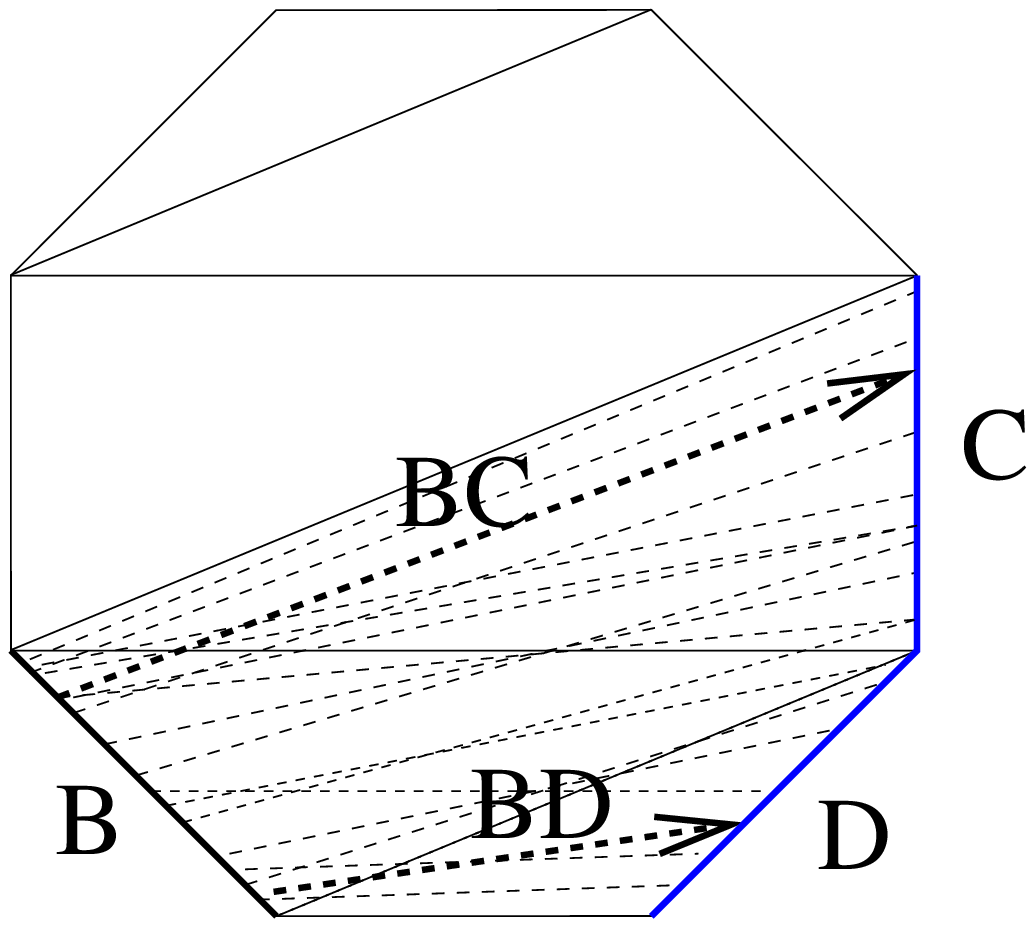}}
\subfigure[Starting with C]{\label{D0transitionsC}
\includegraphics[width=0.23\textwidth]{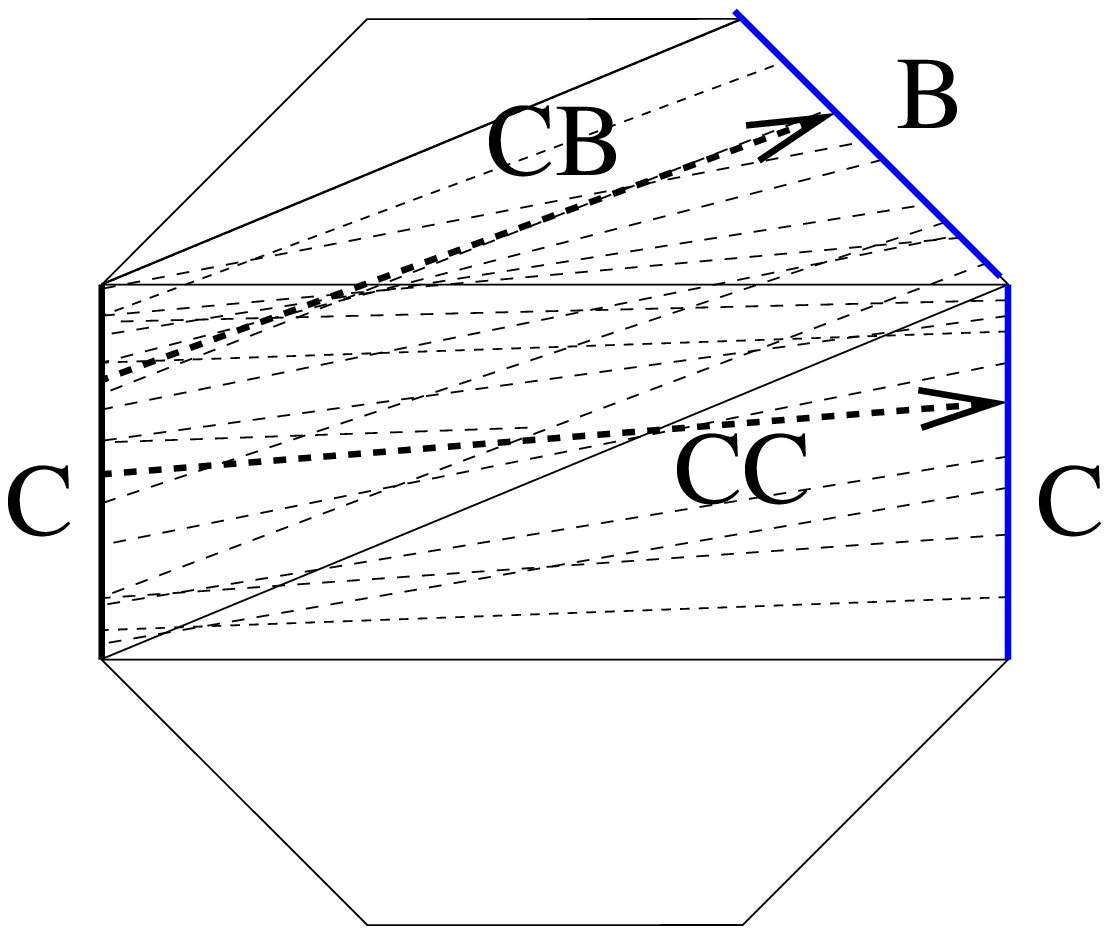}}
\subfigure[Starting with D]{\label{D0transitionsD}
\includegraphics[width=0.23\textwidth]{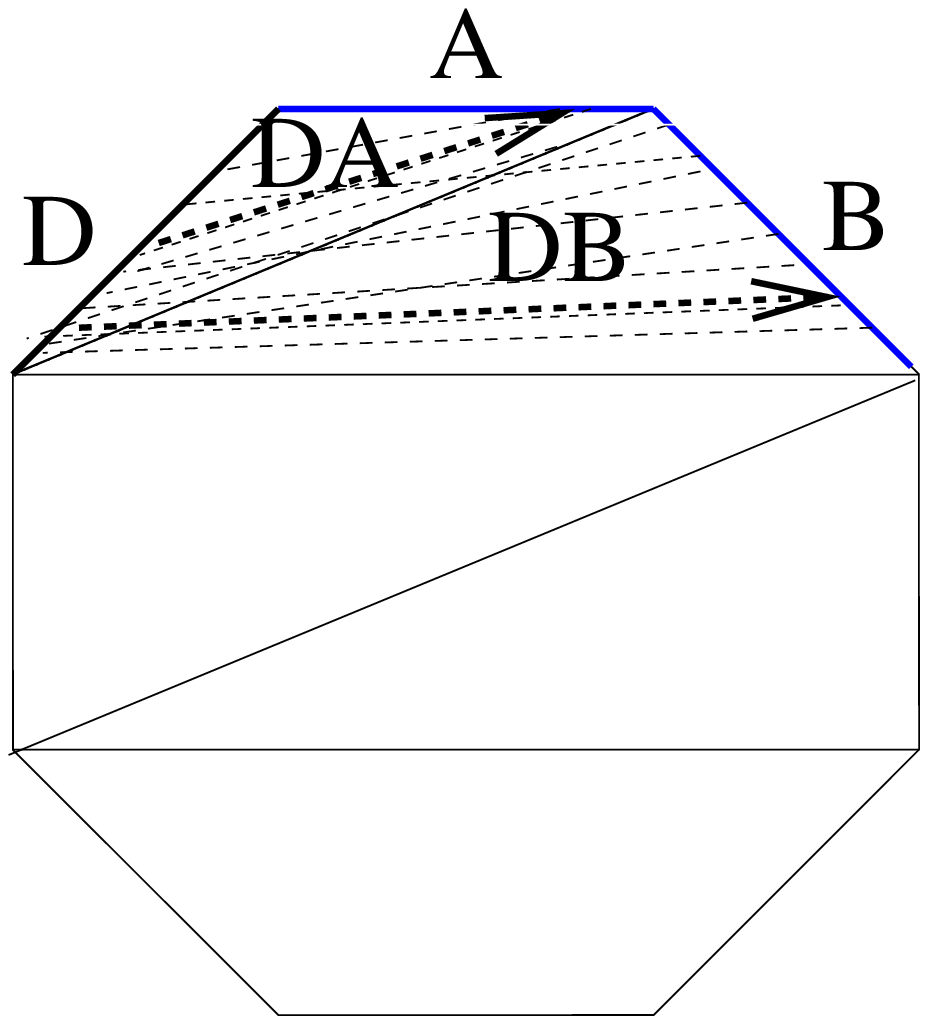}}
\caption{Possible transitions when   $\theta \in \Sigma_0$.}
\end{figure}
Let us first assume that $\theta \in \overline{\Sigma}_0$ and let us represent the transitions by a \emph{transition diagram} $\mathscr{D}_0$, whose vertices are labelled by $\mathscr{A}$ and whose oriented edges represent possible transitions.

\begin{figure}
\centering
\subfigure[$\mathscr{D}_0$\label{D0}]{\includegraphics[width=0.25\textwidth]{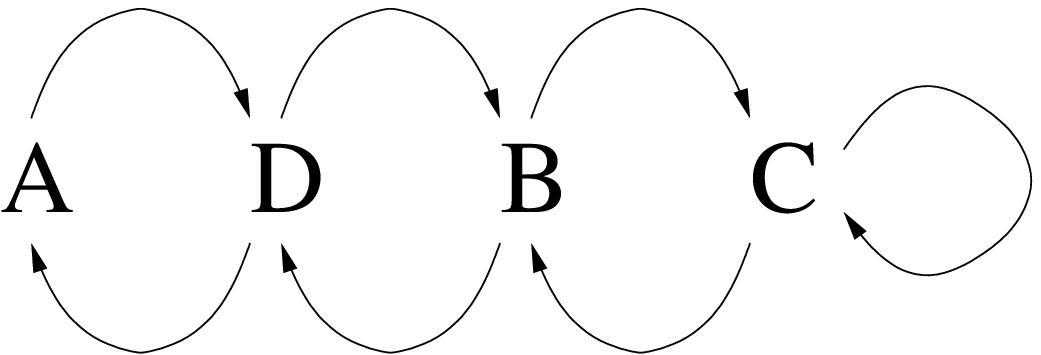}}
\hspace{30mm}
\subfigure[$\mathscr{D}_1$]{\includegraphics[width=0.25\textwidth]{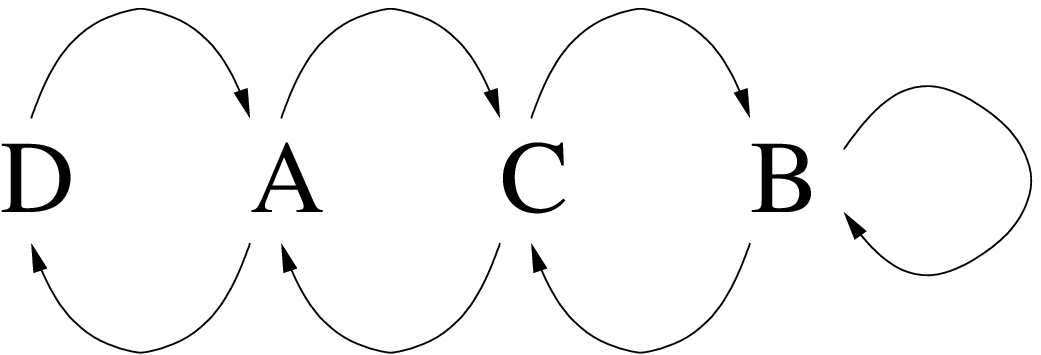}}

\subfigure[$\mathscr{D}_2$]{\includegraphics[width=0.25\textwidth]{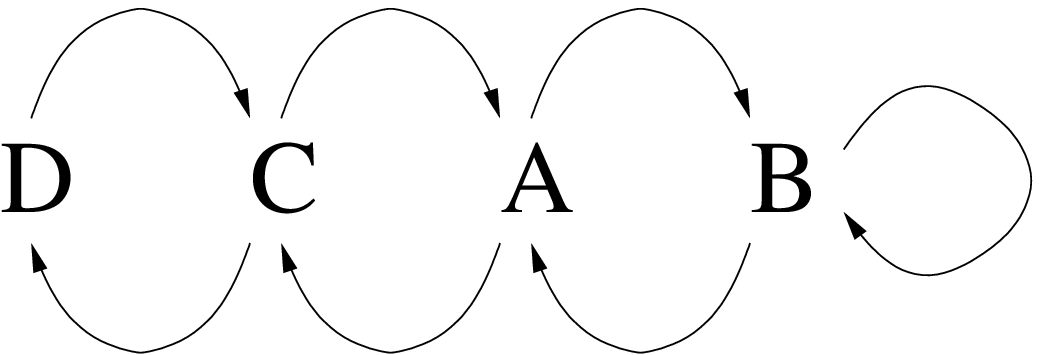}}
\hspace{30mm}
\subfigure[$\mathscr{D}_3$]{\includegraphics[width=0.25\textwidth]{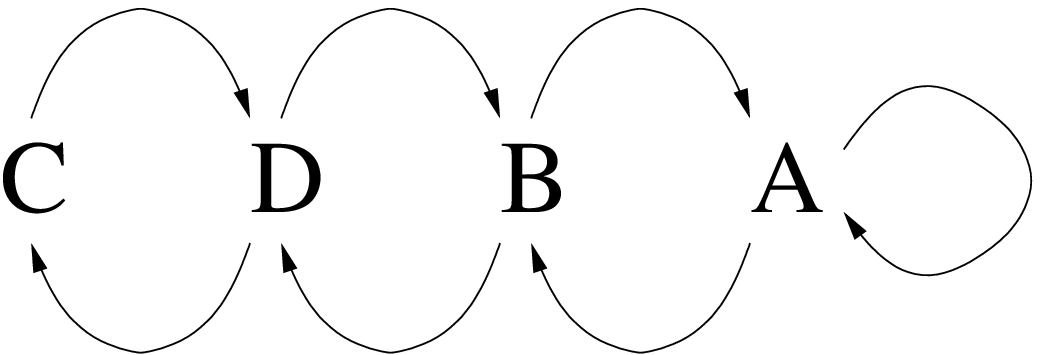}}

\subfigure[$\mathscr{D}_4$]{\includegraphics[width=0.25\textwidth]{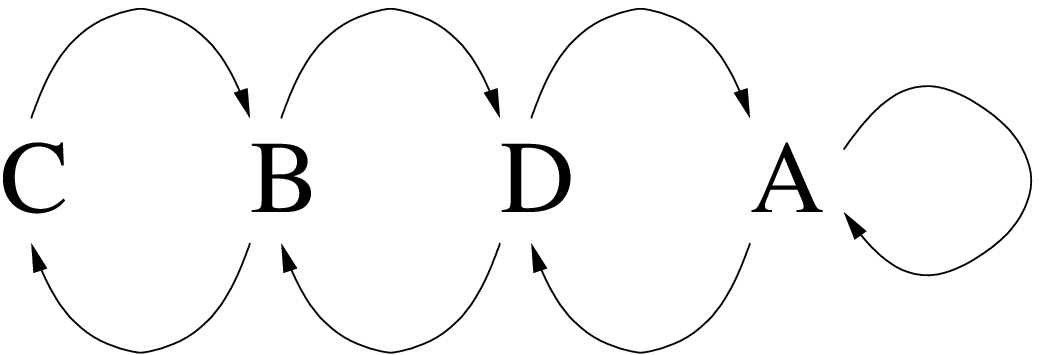}}
\hspace{30mm}
\subfigure[$\mathscr{D}_5$]{\includegraphics[width=0.25\textwidth]{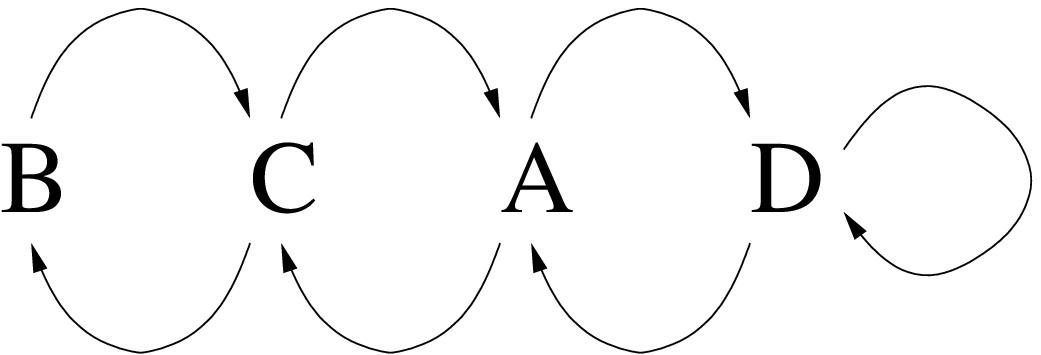}}

\subfigure[$\mathscr{D}_6$]{\includegraphics[width=0.25\textwidth]{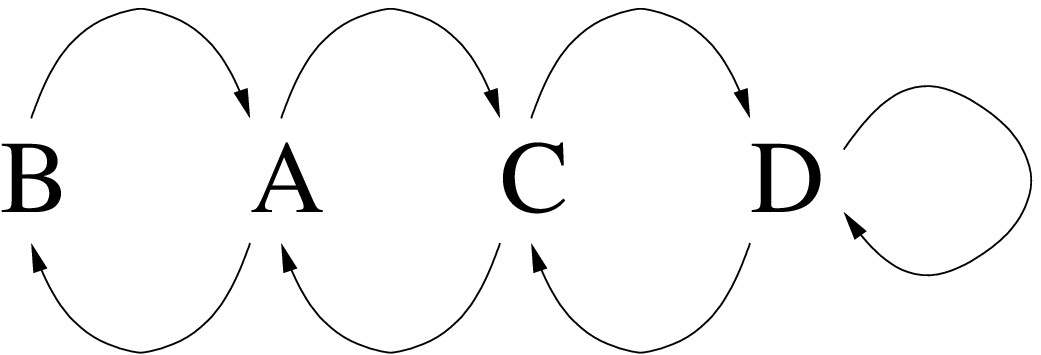}}
\hspace{30mm}
\subfigure[$\mathscr{D}_7$]{\includegraphics[width=0.25\textwidth]{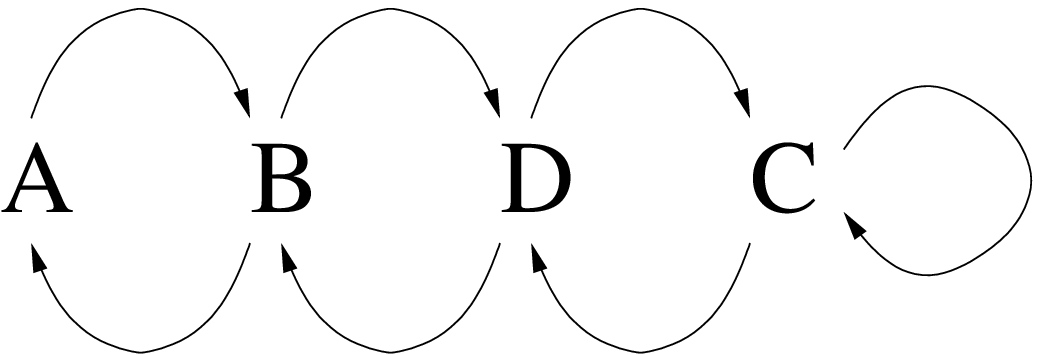}}
\caption{The transition diagrams $\mathscr{D}_0, \dots,\mathscr{D}_7$ corresponding to sectors  $\overline{\Sigma}_0,\dots, \overline{\Sigma}_7$. \label{diagrams}}
\end{figure}

\begin{lemma}\label{sector0transitions}
The diagram $\mathscr{D}_0$ in Figure \ref{D0} gives all  transitions which can occur when $\theta$ belongs to the sector $\overline{\Sigma}_0$. 
\end{lemma}
\begin{proof} A transition $LL'$ corresponds to segment of a linear trajectory contained  in $O$ and traveling from a side labelled $L$ to a side labelled $L'$.
Let us check transitions which begin with $A$.  Since  the segment has direction $\theta \geq 0$, it starts from the side $\overline{O_5O_6}$. Since moreover $\theta \in \overline{\Sigma}_0$, one can easily check that the segment has to hit the side $\overline{O_4O_5}$ labelled by $D$, hence the only possible transition starting with $A$ is $AD$ (see Figure \ref{D0transitionsA}). In the same way, if a cutting sequence contains a $B$, the corresponding segment starts from $\overline{O_6O_7}$. Under the direction assumptions, it can reach any point in the interior of the two sides  $\overline{O_3 O_4}$ or $\overline{O_4 O_5}$: for example, it can be arbitrarily close to the segment  through $O_7$ with $\theta= \pi/8$, which reaches $O_3$, or to the segment from $O_6$ with $\theta= 0$ which reaches  $O_5$, see Figure \ref{D0transitionsB}. 
Thus the two possible transitions starting with $B$ are $BC$ and $BD$ as shown in Figure \ref{D0transitionsB}. Reasoning in this way, one can obtain the complete diagram $\mathscr{D}_0$ in Figure \ref{D0}.   
\end{proof}
\noindent Let us consider now the other transition diagrams $\mathscr{D}_i$, $i=1,\dots, 7$, in Figure \ref{diagrams}. One can check the following.
\begin{lemma}\label{relabel}
The permutation $\pi_i$  associated to the isometry $\nu_i \in \Die_8$ which carries $\overline{\Sigma}_i$ to $\overline{\Sigma}_0$ acts on the labels of vertices of the diagram $\mathscr{D}_i$ by reducing it to $\mathscr{D}_0$.
\end{lemma}
\noindent  Combining Lemma \ref{sector0transitions} with Lemma \ref{relabel} we have the following.
\begin{cor}\label{sectorstransitions}
The transition diagrams $\mathscr{D}_i$, $i=0,\dots, 7$, in Figure \ref{diagrams} give the transitions which can occur in sectors  $\overline{\Sigma}_i$, $i=0,\dots, 7$ respectively. 
\end{cor}

\subsubsection{Admissibility of sequences}
To emphasize the connection between our treatment of the octagon and the treatment of the square we described in section \S~\ref{torussec} we borrow the terms ``admissible sequence" and ``derivation"  that we used in the square case (see \S \ref{torussec}), giving these words appropriate new meanings in the following sections. 
 Let $\mathscr{A}=\{A,B,C,D\}$ and let $\mathscr{A}^{\mathbb{Z}}$ be the space of bi-infinite words $w$ in the letters of $\mathscr{A}$. 
We will use the convention that, given letters $L_i \in \mathscr{A}$, $1\leq i\leq n$, the notation $\overline{L_1 \dots L_n}$  denotes the periodic word of period $n$ obtained by repeating the sequence $L_1 \dots L_n$.



\begin{defn}[Admissibility]\label{admissibledef} Let us say that the word  $w$ is \emph{admissible} if there exists a diagram $\mathscr{D}_i$ for $i\in \{0, \dots, 7\}$ such that all transitions in $w$ correspond to labels of edges of $\mathscr{D}_i$. In this case, we will say that it is \emph{admissible in diagram $i$}. Equivalently, the sequence is admissible in diagram $i$ if it describes an infinite path on $\mathscr{D}_i$. Similarly, a finite word $u$ is admissible (admissible in diagram $i$) if it describes a finite path on a diagram (on $\mathscr{D}_i$). 
\end{defn}
From the previous paragraphs considerations, we have that:
\begin{lemma}\label{admissiblelemma}
Cutting sequences are admissible.
\end{lemma}
\begin{proof}
Since any cutting sequence $w$ is  of the form $c(\tau)$ for a trajectory $\tau$ in direction $\theta$ with  $\theta \in \Sigma_i$ for some $0\leq i \leq 7$, $w$ is admissible by Corollary \ref{sectorstransitions}.
\end{proof}

\begin{figure}
\centering
\subfigure[$\mathscr{D}_{0,1}$\label{D01}]{\includegraphics[width=0.1\textwidth]{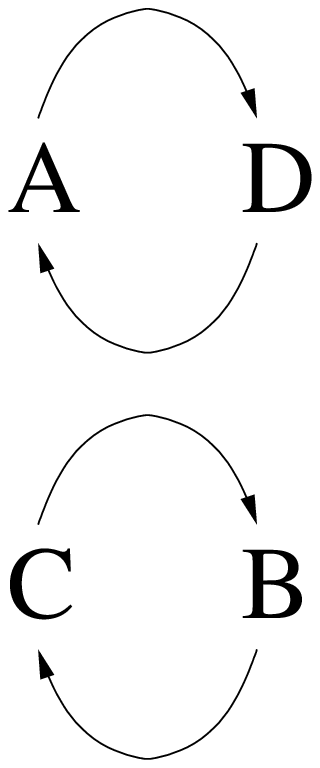}}
\hspace{2mm}
\subfigure[$\mathscr{D}_{1,2}$]{\includegraphics[width=0.1\textwidth]{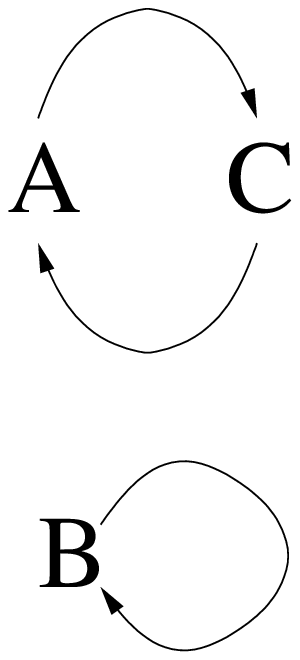}}
\hspace{2mm}
\subfigure[$\mathscr{D}_{2,3}$]{\includegraphics[width=0.1\textwidth]{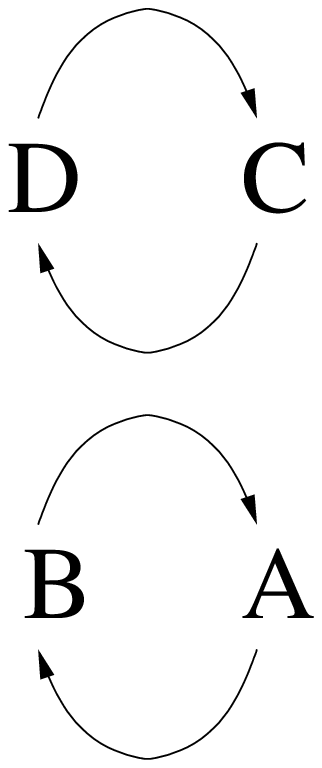}}
\hspace{2mm}
\subfigure[$\mathscr{D}_{3,4}$]{\includegraphics[width=0.1\textwidth]{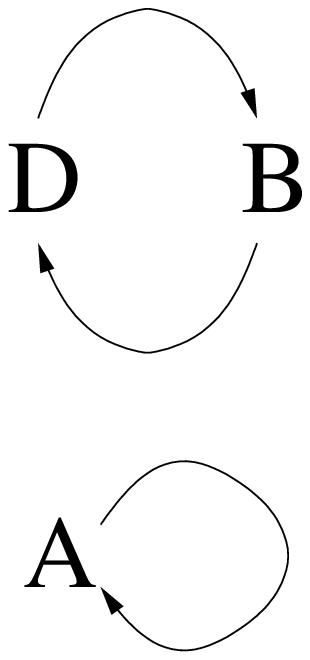}}
\hspace{2mm}
\subfigure[$\mathscr{D}_{4,5}$]{\includegraphics[width=0.1\textwidth]{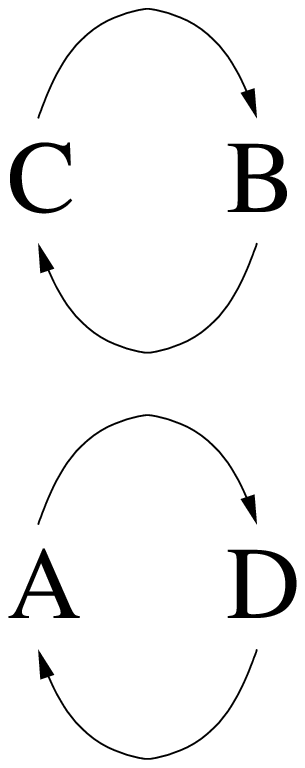}}
\hspace{2mm}
\subfigure[$\mathscr{D}_{5,6}$\label{D56fig}]{\includegraphics[width=0.1\textwidth]{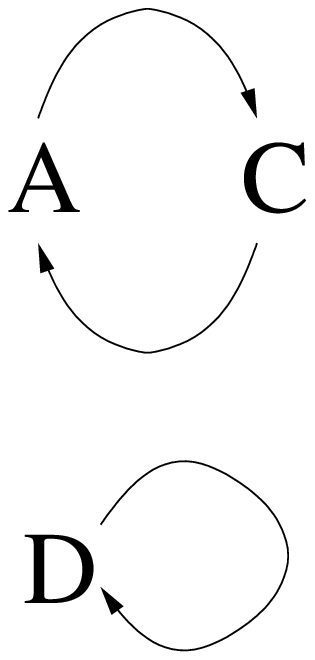}}
\hspace{2mm}
\subfigure[$\mathscr{D}_{6,7}$\label{D67fig}]{\includegraphics[width=0.1\textwidth]{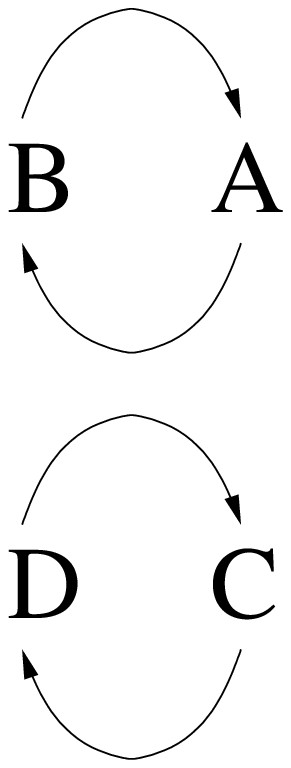}}
\hspace{2mm}
\subfigure[$\mathscr{D}_{7,0}$]{\includegraphics[width=0.1\textwidth]{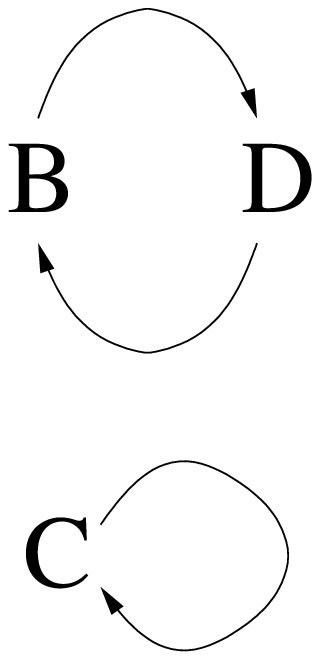}}
\caption{The transitions diagrams $\mathscr{D}_{k-1,k}$   for directions $\theta = k\pi/8$, $k=1,\dots, 7$. \label{boundarydiagrams}}
\end{figure}

 We remark that some words are admissible in more than one diagram.   
 \begin{ex}\label{boundaryambiguityex} The periodic sequence $w$ of period one given by $\overline{ B} $ is admissible both in diagram $1$ and in diagram $2$. 
\end{ex} 
\noindent If $\theta = k\pi/8$ for some $k=1,\dots, 7$, then $\theta= \overline{\Sigma}_{k-1}\cap \overline{\Sigma}_{k}$ is a boundary point of two adjacent sectors and $c(\tau)$ is admissible in both diagram $k$ and $k-1$. 
 The non-uniqueness of the diagram  can also occur in non-neighboring sectors.
\begin{ex}\label{truambiguityex}
The periodic sequence $w = \overline{ ADBD} $ of period 4 is admissible both in diagram $0$ and in diagram $4$. The periodic word $\overline{AD}$ is admissible in diagrams $0$, $1$, $4$ and $5$.
\end{ex}
\noindent If a sequence is admissible at the same time  both in diagram $k-1$ and in diagram ${k}$ 
only the transitions which are allowed in both $\mathscr{D}_{k-1}$ and $\mathscr{D}_{k}$ can occur. Figure \ref{boundarydiagrams} shows the eight corresponding diagrams, $\mathscr{D}_{k-1,k}$,  $k=1,\dots, 7$ and $\mathscr{D}_{7,0}$. 
The diagrams show that all such sequences are periodic sequences of period one or two. Furthermore, simply by looking at the octagon, one can immediately observe the following stronger fact.
\begin{rem}\label{periodicmultipleofpi8}
All infinite trajectories in directions of the form $\theta = k\pi/8$, $k\in \mathbb{N}$ are periodic trajectories.   
\end{rem}

\begin{lemma}\label{all_transitions} If a cutting sequence $w$ contains all the transitions in some diagram $\mathscr{D}_i$,  then $w$ is admissible \emph{only} in diagram $i$.
\end{lemma}
\begin{proof} The proof follows by inspecting the diagrams in Figure \ref{diagrams}.\end{proof}

Let us define a \emph{normal form} for  words which is analogous to the normal form for trajectories.
\begin{defn} \label{normalword} Say that a word $w$ is admissible in a \emph{unique} diagram $\mathscr{D}_i$, for some $0\leq i \leq 7$. The \emph{normal form} of $w$ is the word $n(w):=\pi_i \cdot w_0$.\end{defn}

\subsubsection{Derivation and infinitely derivable sequences}
Let us define a combinatorial operation, that we call \emph{derivation}, on cutting sequences of linear trajectories in the octagon.  
\begin{defn}[Sandwiched property] Let us say that a letter 
 in a  sequence  $w\in \mathscr{A}^{\mathbb{Z}}$ is \emph{sandwiched} if it is preceded and followed by the same letter. Moreover, if $L\in \mathscr{A}$ is sandwiched and $L' \in \mathscr{A}$ is the letter which precedes and follows $L$, we will say that $L$ is  \emph{$L'$-sandwiched}.
\end{defn}
\begin{ex}
In the finite subword $CACCCDBD$, the letter $A$ is $C$-sandwiched since it is followed and preceded by $C$; the letter $B$ is followed and preceded by $D$, hence it is $D$-sandwiched; similarly the third $C$ is $C$-sandwiched.
\end{ex}
\noindent Let us remark that in a finite word it is not possible to determine whether the first and last letter are sandwiched. This is one of the reasons to define derivation only for bi-infinite words. We consider finite cutting sequences in a forthcoming paper.   

\begin{defn}[Derivation]
Given a word $w \in \mathscr{A}^{\mathbb{Z}}$ the \emph{derived sequence}, which we denote by  $w'$, is the word obtained by keeping only the letters of $w$ which are sandwiched. 
\end{defn}
\begin{ex}
If $w$ contains the finite word $CACCCDBDCDC$ the derived sequence $w'$ contains the word $ACBCD$, since these are sandwiched letters in the string. 
\end{ex}
\begin{rem}\label{commute}
The property of being sandwiched is invariant under permutations of the letters. Hence, the operation of derivation on a word commutes with the action of permutations of the letters, i.e.~$(\pi \cdot w)' = \pi \cdot w'$.  
\end{rem} 
\begin{defn}[Derivable and infinitely derivable]
A word $w \in \mathscr{A}^{\mathbb{Z}}$ is \emph{derivable} if it is admissible 
and its derived sequence $w'$ is admissible. 

A word  $w \in \mathscr{A}^{\mathbb{Z}}$ is  \emph{infinitely derivable} if it is derivable and for each $n$ the result of deriving the sequence $n$ times is again derivable. 
\end{defn}

\noindent Using this notion of derivability, 
we have the following key Proposition.

\begin{prop}\label{derivationiscuttingseqprop}
Given a cutting sequence, the derived  sequence is again a cutting sequence. 
\end{prop}
\noindent 
Section \S \ref{derivationsec} is devoted to the proof of Proposition \ref{derivationiscuttingseqprop}. 
The following necessary condition on cutting sequences follows  immediately, as in the case of the square, 
from Proposition \ref{derivationiscuttingseqprop}. 
\begin{thm}\label{cutseqinfderiv}
A  cutting sequence $c(\tau)$ is infinitely derivable.
\end{thm}
\begin{proof} A cutting sequence is infinitely derivable since by Proposition \ref{derivationiscuttingseqprop} each derived sequence is again a cutting sequence and by Lemma \ref{admissiblelemma} cutting sequences are admissible. 
\end{proof}
\noindent 


\subsection{Renormalization schemes}\label{renormalizationsec}  
In this section we define three \emph{renormalization operators}, one combinatorial or at the level  of symbolic sequences (\S\ref{combrenormsec}), one acting on the space on trajectories (\S\ref{trajrenormsec}) and one on directions, given by iterations of an expanding map (the octagon Farey map, see \S\ref{Fareymapsec}).
These renormalization schemes are useful in formulating some of our results (such as the direction recognition algorithm in Theorem \ref{directionsthm} and the algorithm to generate cutting sequences in \ref{algorithmsec}) and they turn out to be crucial ingredients in the proofs. In the following sections we also state how they are related to one another. 

\subsubsection{Combinatorial renormalization scheme}\label{combrenormsec}
Given a word $w$ which is infinitely derivable such as a cutting sequence $w=c(\tau)$ (which is infinitely derivable by Theorem \ref{cutseqinfderiv}) it will be useful to have information about the sequence of diagrams in which the derived sequences are admissible. It is convenient for us to put the sequences into a {normal form}  after each derivation, as explained below. The use of this convention when recording admissible diagrams will allow us to use a relatively simple continued fraction map (in \S \ref{Fareymapsec}).
 
Set $w_0:= w$. Let us assume for now that $w_0$ is admissible in an \emph{unique} diagram $d_0$, with $0\leq d_0 \leq 7$. The first step of the renormalization scheme  consists of taking the derived sequence of the normal form of $w_0$ (see Definition \ref{normalword}) and setting  $w_1 := n(w_0)'$.
If after $k$ iterations the sequence obtained is again admissible in an \emph{unique} diagram, we  define
\be\label{renormalizedseq}
 w_{k+1}:= n(w_k)'.
\ee

The above assumptions hold for a large class of sequences. By \emph{non-periodic} sequence, we mean as usual that the letters in the sequence do not repeat periodically. 
\begin{prop}\label{admsectorsfornonper}
If $w$ is a \emph{non-periodic cutting sequence}, then the renormalization scheme described above is well defined for all $k\in \mathbb{N}$.   In this case, $w$ determines a unique infinite sequence $\{d_k\}_{k\in \mathbb{N}}$ such that $w_k$ is admissible in diagram  $d_k$. 
\end{prop}
\begin{defn}\label{admissiblesectorsdef} We refer to the sequence  $\{d_k\}_{k\in\mathbb{N}}$ as the \emph{sequence of admissible diagrams}. 
\end{defn}
The proof of Proposition \ref{admsectorsfornonper} is given in \S\ref{directionsthmproof}.  
\noindent The sequence of admissible diagrams of a cutting sequence $c(\tau)$ can be used to recover the direction of the trajectory $\tau$, see Theorem \ref{directionsthm}.  
We give here an example of computation of the first entries of the sequence of admissible diagrams.

\begin{ex}\label{typesex} Consider a sequence $w$ which contain the following finite subsequence 
\vspace{-1.6mm} \be\begin{split}
w =  \dots & AADBDA   A A A B D B C B D  B D A A A A   D B D A A A A \\ & D B D A A A ADBDBCBDBDAAADBDB  DAAADB \dots  
\end{split}
\ee
\noindent This sequence is admissible in $\mathscr{D}_4$. Since it contains $AA$, present only in $\mathscr{D}_3$ and $\mathscr{D}_4$, and $CB$, which is not a transition in $\mathscr{D}_4$, it is  admissible only in diagram $4$.  Thus, $d_0=4$ and $n(w)$ is obtained applying $\pi_4 = (AC)(BD)$, which gives \\
\vspace{-1.6mm} 
\be\begin{split}
n(w)=\dots & CCBDBCCCCDBDADBDBCCCCBDBCCCC\\ & BDBCCCC BDB DADBDBCCCBDBDBCCCBD\dots . 
\end{split}
\ee
\vspace{-3mm}

\noindent
Thus the derived sequence $w_1= n(w)'$ contains 
\vspace{-1.6mm}
\be
 \dots DCCDBABDCCDCCDCCDBABDCDBDC\dots
\ee
\noindent
 and, since it contains $CC$ (which is an admissible transition only in diagrams $7$ and $0$) and $AB$ (which excludes sector $0$), it is of type $7$, so $d_1=7$. 
To find $d_2$, we put this word in normal form  by applying $\pi_7=(BD)$. We obtain  
\vspace{-1.5mm}  \be 
 n(w_1)= \dots BCCBDADBCCBCCBCCBDADBCBDBC \dots.  
\ee 
\noindent  Deriving once more, we get $w_2= n(w_1)'= \dots ABBACD \dots$ and since this sequence is admissible only in diagram $2$, we get $d_2=2$ and $n(w_2)= \cdots BCCBDA \cdots$ and so on. 
\end{ex}
We have already observerd (see Examples \ref{boundaryambiguityex}, \ref{truambiguityex}) that a periodic trajectory can be admissible in more than one diagram, so the assumption of the renormalization scheme does not hold for periodic words. 
In Example \ref{boundaryambiguityex}, $w$ is admissible in two possible diagrams, but one can check that if we act by $\pi_1$ or $\pi_2$ we obtain the same sequence, so one could still uniquely define a normal form. On the other hand, if we consider $w = \overline{ ADBD} $ in  Example \ref{truambiguityex}, $\pi_0 \cdot w$ and $\pi_4\cdot w$ are different and hence to define the normalization one would need to choose one of the two as normal form.
We show in the next section that one can still define a non-ambiguous renormalization scheme for periodic cutting sequences if the direction of the periodic cutting sequence is given (see Remark \ref{renormalizeperiodic}).  

\subsubsection{Renormalization on trajectories}\label{trajrenormsec}
Let us define a renormalization scheme on trajectories which is the natural counterpart of the combinatorial operation on their cutting sequences. Given a trajectory $\tau$ in direction $\theta$, in order to refer to the sector to which $\theta$ belong, let us introduce the following notation. 
\begin{defn}[Sector of a trajectory]\label{sectordef}
Given a trajectory  $\tau$ in direction $\theta$, let $s(\tau)=k \in \{ 0, \dots, 7 \}$ be such that $\theta \in \Sigma_{k}$. We say that $s(\tau)$ is the \emph{sector} of $\tau$. 
\end{defn}
\noindent Clearly, if the \emph{sector} of $\tau$ is $k$, the cutting sequence  $c(\tau)$ is \emph{admissible} in diagram $k$. 
The converse is not true (see Examples \ref{boundaryambiguityex} and \ref{truambiguityex}), but the following simple remark holds.
\begin{rem}\label{uniquenesstyperk} If there exists a unique $k\in \{0, \dots, 7\} $ such that $w=c(\tau)$ is admissible in diagram $k$, then $s(\tau)=k$. In this case, we can identify the sector $s(\tau)$ only from the knowledge of the combinatorial sequence $w$, without a priori knowing its direction. 
\end{rem}

Let us recall that given a trajectory $\tau$, we defined in Definition \ref{normalizedtrajdef} its normal form $n(\tau)$.   
We recall also that given a cutting sequence $c(\tau)$ of a trajectory $\tau$, $c(\tau)'$ is also a cutting sequence (see Proposition \ref{derivationiscuttingseqprop}). More precisely,  given a trajectory $\tau$ with $s(\tau)=0$, in Proposition \ref{derivationiscuttingseqpropsec0} we \emph{define} a trajectory $\tau'$ such that
\be\label{derivedseqeq}
 c(\tau)'=c(\tau').
\ee 
\begin{defn}\label{derivedtrajdef}
We refer to the trajectory $\tau'$  which satisfies (\ref{derivedseqeq}) and is explicitly defined in Proposition \ref{derivationiscuttingseqpropsec0} as the \emph{derived trajectory} of $\tau$.
\end{defn}

The renormalization scheme on the space of trajectories is obtained by alternatively putting in normal form and deriving trajectories as follows. 
\begin{defn}[Renormalization on trajectories]\label{renormalizationtrajectoriesdef}
Given a trajectory $\tau$, let us recursively define a sequence of renormalized trajectories $\{\tau_k\}_{k\in \mathbb{N}}$ by:
\bes
\tau_0:= \tau; \qquad \tau_{k+1}:= n( \tau_k)', \quad k\in \mathbb{N}. 
\ees
\end{defn}
\begin{defn}[Sequence of sectors]\label{sectorsdef}
Given a trajectory $\tau$, we denote by $\{s_k(\tau)\}_{k\in \mathbb{N}}$ the \emph{sequence of sectors} given by
\be\label{sectorsseqdef}
 s_k(\tau):= s(\tau_k),
\ee
where $s(\tau_k)$ is the sector of the $k^{th}$ renormalized trajectory $\tau_k$ in Definition \ref{renormalizationtrajectoriesdef}.
\end{defn}
\noindent In words, the sequence $\{s_k(\tau)\}_{k \in \mathbb{N}}$ is obtained  by 
recording the sectors of the renormalized trajectories. 
Corollary \ref{typesequenceslemma} 
 gives a characterization of sequences in $\{0,\dots, 7\}^{\mathbb{N}}$ which can occur as sequences of types.

\paragraph{Trajectories versus sequences renormalizations.} 
The following Lemma shows that renormalization on sequences, when it is well defined, is the natural combinatorial counterpart of renormalization on trajectories. 
\begin{lemma}\label{cutseqrelationlemma}
Let $w$ be the cutting sequence $c(\tau)$ of a trajectory $\tau$. 
Assume that the $k^{th}$ renormalized sequence $w_k$ given by (\ref{renormalizedseq}) is well defined. Then
\be\label{cuttseqrelation}
c(\tau_k) = w_k.
\ee
\end{lemma}
\begin{proof}
The proof is by induction on $ k$. If $k=0$, then $c(\tau_0)=c(\tau)=w=w_0$ by definition. Assume (\ref{cuttseqrelation}) holds for $k$.  
 Since by Definition \ref{sectorsdef} the sector of $\tau_k$ is $s_k=s_k(\tau)$, we have  by Definition \ref{normalizedtrajdef} that $n(\tau_k)=\nu_{s_k} \trajcdot \tau_k$. Moreover, since  by induction we know that $w_k = c(\tau_k)$, $w_k$ is admissible in diagram $s_k$. Since by assumption this is the unique diagram in which it is admissible, we have $d_k(w)= s_k$ and, by (\ref{renormalizedseq}), $n(w_k)=\pi_{s_k} \cdot w_k  $.  Since  by (\ref{relabellingeq}) we get $c( \nu_{s_k} \trajcdot \tau_k )= \pi_{s_k} \cdot w_k $, this shows that $c(n(\tau_k))   =  n(w_k)$. Moreover, by 
 (\ref{derivedseqeq}), 
  $c(  n(\tau_k)  '   ) =  c(n(\tau_k))'$. Combining these two observations with the definitions $\tau_{k+1} = n(\tau_k)'$ and $w_{k+1} = n(w_k)'$ (Definition \ref{renormalizationtrajectoriesdef} and (\ref{renormalizedseq}) respectively), we get 
$c( \tau_{k+1}) = c( n(\tau_k)' )= c(n(\tau_k))' = n(w_{k})'= w_{k+1}$, as desired.
\end{proof} 

\begin{rem}\label{renormalizeperiodic}
Since $\tau_k$ is well defined for all $\tau$,  equation (\ref{cuttseqrelation}) can be used to define the renormalized sequences $w_k$ when $w =c(\tau)$ is a periodic cutting sequence. 
\end{rem}

\subsubsection{The octagon Farey map and renormalization on directions.}\label{Fareymapsec} 

Let us define an additive continued fraction algorithm. We will eventually show that this algorithm plays the role for the octagon that the standard additive continued fraction plays for the square. In particular we will be able to relate the cutting sequence of a trajectory in a given direction to the appropriate continued fraction expansion of the direction. The origin of this map will be clear in \S\ref{derivationsec}. 

Let $\RP$ be the space of directions. Let $\widetilde{\Sigma}_i \subset \RP$ be the set of directions corresponding to the sector $\overline{\Sigma}_i$,  $i=0,\dots, 7$. The linear maps $\nu_i$, $i=0,\dots, 7$,
 given in  (\ref{nujdef}) induce an action on $\RP$ so that  $\nu_i:  \widetilde{\Sigma}_i \to \widetilde{\Sigma}_0$. 
\bes
\gamma : = \begin{pmatrix} -1 &  2(1+\sqrt{2})  \\ 0  & 1  \end{pmatrix}. 
\ees
\noindent The map $\gamma$ takes $\widetilde{\Sigma}_0$ to the union of $\widetilde{\Sigma}_j$ for $j\ge 1$. Let $F_i: \widetilde{\Sigma}_i \to\RP$ be the map induced by the linear map $\gamma\nu_i$,  $i=0,\dots, 7$. Define the \emph{octagon Farey map} $F : \RP\to\RP$ to be the piecewise-projective map, whose action on $\widetilde{\Sigma}_i$ is given by $F_i$. Observe that these branches fit together so that $F$ is continuous.

There are two coordinate systems on $\RP$ which will prove to be useful for describing directions of trajectories.  The first is the the \emph{inverse-slope coordinate} $\mu$, where $\mu((x,y))=x/y$. In these coordinates  $F$ is a piecewise linear fractional transformation. 

 If $ \mu  \in  \left[ \cot ({\pi i }/{8}), \cot  (i+1)\pi/8 \right) $ we have  
\bes \label{Fareymapdefinvslopes}
 F(\mu)  =  \frac{a_i \mu + b_i}{c_i \mu + d_i}, \qquad \mathrm{where} \, \, \begin{pmatrix} a_i & b_i \\ c_i & d_i \end{pmatrix}:= \gamma \nu_i  . 
\ees
\noindent It is also useful to describe $F$ in terms of  \emph{angle coordinate} $\theta$, where $\theta = \cot(\mu) \in [0,\pi]$. The action in this coordinates is obtained by conjugating by $\cot : [0,\pi] \rightarrow \mathbb{R}$ as follows.  
Given, $\nu \in GL(2, \mathbb{R})$ on $\mathbb{R}$, we denote by $\nu [ y]$ the conjugated  action by fractional transformation given by
\be\label{actionslopes}
\nu [\theta] =   \cot^{-1}\left( \frac{a \cot(\theta) + b}{c \cot(\theta) + d}\right),  \qquad  \mathrm{if} \,    \qquad \nu :=  \begin{pmatrix} a & b \\ c & d \end{pmatrix} , \qquad \theta \in  [0,\pi].
\ee 
Thus, the octagon Farey map in the coordinate $\theta$ is given by  
\begin{equation}\label{Fareymapdef}
 F(\theta) = F_i (\theta) :=  (\gamma \nu_i)[\theta],   \qquad  \mathrm{if} \quad   \theta  \in  \overline{\Sigma}_i,  \quad i=0, \dots, 7.  
\ee
The map $F$ is well defined and continuous since $F_i(\theta) = F_{i+1}(\theta)$ if $\theta \in \overline{\Sigma}_i \cap \overline{\Sigma}_{i+1}$. 
The graph of $F$ in angle coordinates is shown in Figure \ref{Fareygraph}. This map is pointwise expanding but not uniformly expanding.
Clearly, since for each $ i=0, \dots, 7$ the branch $F_i$ of $F$ is monotonic,  the inverse maps $F_i^{-1}: [\frac{\pi}{8}, \pi] \to \overline{\Sigma}_i$ are well defined. 
\begin{figure}
\centering
{\includegraphics[width=0.6\textwidth, height=0.6\textwidth
]{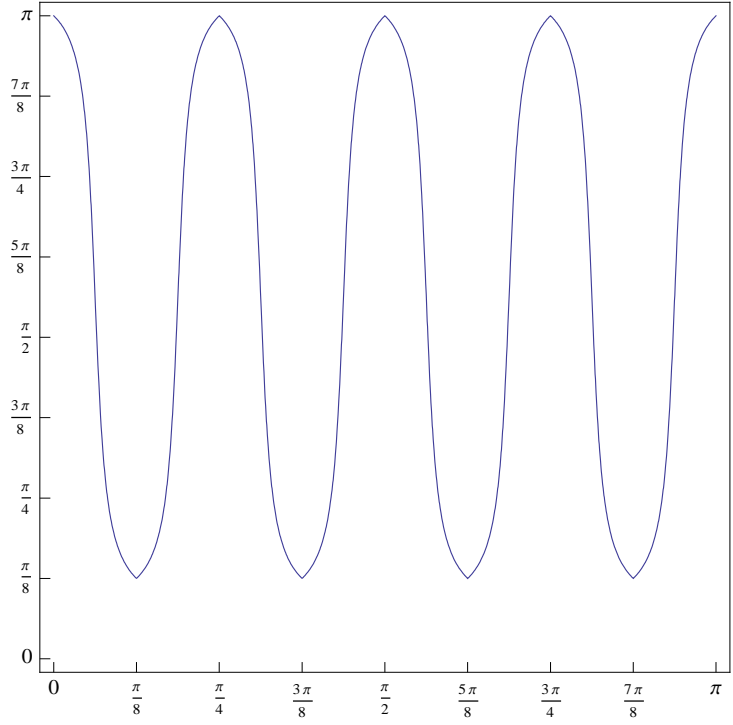}}
\caption{The graph of the octagon Farey map $F$.\label{Fareygraph}}
\end{figure}
\begin{defn}[Itinerary]\label{itinerarydef}
To any $\theta\in[0,\pi]$ we can assign the sequence $s_0,s_1,s_2 \ldots $ defined by  $F^j(\theta)\in\Sigma_{s_j}$ for each $j \in \mathbb{N}$. 
We call this sequence the \emph{itinerary} of $\theta$ under $F$. 
\end{defn}
\noindent The itinerary gives  the \emph{symbolic coding} of the map $F$ with respect to the natural partition  $\{ \Sigma_0, \cdots, \Sigma_7\}$.     
Let us remark that $F([0,\pi])\subset [\pi/8,\pi]$. This immediately gives the following.
\begin{lemma} 
\label{itineraryproperties} 
If $s_0,s_1,s_2,\dots$ is an itinerary, then it satisfies the following condition that
 $s_0\in \{0,\dots, 7\} $ and $s_k \in \{ 1,\dots, 7\} $ for all $k\ge1$. 
 \end{lemma}
 
\paragraph{The octagon additive continued fraction expansion.}
The map $F$ gives an additive continued fraction algorithm for numbers in the interval $[0,\pi]$. 
\begin{defn} \label{Sstardef} Let $S^*$ to be the set of all sequences $s=\{s_k\}_{k \in \mathbb{N}}$ that satisfy the condition $s_k\in \{0,\dots, 7\} $ and $s_k=0$ implies $k=0$.
\end{defn}
\begin{defn}Given  $\{s_k\}_{k \in \mathbb{N}} \in S^* $, let 
\be\label{sectorCFdef}
\overline{\Sigma}[s_0; s_1, \dots, s_k] :=  \bigcap_{k \in \mathbb{N}} F_{s_0 }^{-1} F_{s_1 }^{-1} \cdots  F_{s_k}^{-1} [0, \pi]  . 
\ee
\end{defn}
\noindent One can check that $ \overline{\Sigma}[s_0; s_1, \dots, s_k] $ is a non-empty closed interval. The  itinerary  of the points in the interior of $ {\overline{\Sigma}}[s_0; s_1, \dots, s_k] $  begins by $s_0, \dots, s_k$ by construction.
\begin{lemma}
\label{contfracprop1} 
Given $s= \{s_k\}_{k\in\mathbb{N}}  \in S^*$,  the intersection
\be\label{directionasintersection}
 \bigcap_{k \in \mathbb{N}} \overline{\Sigma}[s_0; s_1, \dots, s_k]  
\ee
is non empty and consists of an unique number $\theta \in [0,\pi]$.
\end{lemma}
\begin{defn}[Octagon additive continued fraction expansion]\label{CFexpansion} If 
\bes
\{ \theta \} = \bigcap_{k \in \mathbb{N}} \overline{\Sigma}[s_0; s_1, \dots, s_k]   , \qquad \mathrm{with} \quad \{s_k\}_{k\in\mathbb{N}} \in S^*,
\ees
 we write\footnote{We chose this notation to emphasize that the sequence determines an expansion of the angle $\theta$. The reader should be aware though that the entries of the usual continued fraction give the itinerary of the Gauss map, while the entries $s_0$, $s_1, \ldots$ in Definition \ref{CFexpansion} should be compared to the itinerary of the Farey map in (\ref{classicalFarey}), which corresponds to the \emph{additive} continued fraction.  The (multiplicative) continued fraction entries can be obtained by grouping together consecutive $0$s in the Farey itinerary Farey itinerary entries, which are either $0$ or $1$. Similarly, one could obtain a multiplicative version of the octagon continued fraction by grouping together consecutive $0$s and consecutive $7$s in the octagon Farey itinerary but we do not make use of this version in this paper.} 
$\theta := [s_0; s_1, s_2, \ldots ]$ and say that $[s_0; s_1, \ldots ]$ is an \emph{octagon additive continued fraction expansion} of $\theta$ and that $s_0,s_1,\ldots$ are its \emph{entries}. 
\end{defn}

Let us remark that the itinerary of a point with respect to the semi-open intervals $\Sigma_i$ is uniquely defined since these intervals are disjoint. On the other hand, in the additive octagon continued fraction Definition \ref{CFexpansion},  we consider itineraries with respect to the closed intervals $\overline{\Sigma}_i$, since we use (\ref{sectorCFdef}) and branches $F_i$ include endpoints. Thus ambiguity is possible, as  explicitly described in the following Lemma. 

\begin{lemma}\label{contfracprop2}
Every $\theta \in [0,\pi]$ has at least one and at most two expansions  as in Definition \ref{CFexpansion}. More precisely:
\begin{itemize}
\item[(i)] 
If $s_0,s_1,s_2,\dots$ is the itinerary of $\theta$ under $F$, then an octagon additive continued fraction expansion of $\theta$  is given by  $[s_0; s_1, s_2, \ldots ]$;
\item[(ii)]
If $\theta$ has two additive octagon continued fraction expansions, the entries are either eventually all 
$1$'s or eventually all $7$'s and the two expansions  have one of the two following forms:
\begin{eqnarray}
\theta &= &[s_0; s_1, \dots, s_{k}, 1, 1, \dots ] = [s_0;s_1, \dots, s_{k}-1, 1, 1, \dots ]  \quad  s_k \in \{ 3,5,7\} 
\label{2CFeven}\\
 \theta &= &[s_0;s_1 \dots, s_{k}, 7, 7, \dots ] = [s_0; s_1, \dots, s_{k}-1, 7, 7, \dots ] \quad  s_k \in \{ 2,4,6\} 
 \label{2CFodd}.
 \end{eqnarray}
In this case, only the first expansion in (\ref{2CFeven}, \ref{2CFodd})  is given by the itinerary of $\theta$. 
\end{itemize}
\end{lemma}
\noindent
Lemma \ref{contfracprop1} and Lemma \ref{contfracprop2} follow from standard techniques in the theory of symbolic coding for piecewise expanding map.  For completeness, we included their proofs in \S\ref{contfracpropproof}.

\begin{defn}[Terminating directions]\label{terminatingdef}
Let us call a direction $\theta$  \emph{terminating} if the  entries of the octagon additive continued fraction expansion  of $\theta$ are eventually $1$ or eventually $7$.  
\end{defn}
 With the exception of $\pi/8$ and $\pi$,   terminating  directions have two octagon additive continued fraction expansions.  
The ambiguity of coding of terminating directions is intrinsic in our 
algorithm. It is a analogous to the intrinsic ambiguity of the continued fraction expansion of a rational number. In the case of the classical continued fraction algorithm, the convention of terminating the expansions of rationals is usually adopted.   This is our reason for the choice of the name \emph{terminating}. 

\paragraph{Direction versus trajectories renormalization.} 
The following Proposition shows that the renormalization scheme on directions given by iterations of the octagon Farey map is the natural counterpart of the renormalization scheme on trajectories defined in \S\ref{trajrenormsec}.

\begin{prop}\label{typesequalCF} 
Let $\tau$ be a trajectory in direction $\theta$. For each $k\in\mathbb{N}$, the direction $\theta_k$ of the $k^{th}$ renormalized trajectory $\tau_k$ in Definition \ref{renormalizationtrajectoriesdef} is given by 
\be
\theta_k = F^k (\theta).
\ee
In particular, the sequence $\{s_k(\tau)\}_{k\in\mathbb{N}}$ of sectors of a trajectory $\tau$  in direction $\theta$ coincides with the itinerary of $\theta$ under the octagon Farey map $F$,
 i.e.
\bes
F^{k} (\theta) \in \Sigma_{s_k(\tau)}, \qquad  \forall   \, \, k\in \mathbb{N}. 
\ees
\end{prop}
\noindent Combining Lemma \ref{itineraryproperties}, Lemma \ref{contfracprop1} and Lemma \ref{contfracprop2} with Proposition \ref{typesequalCF} we get the characterization of sequences of sectors. 
\begin{cor}[Characterization of sequences of sectors]\label{typesequenceslemma}
A sequence $s=\{ s_k \}_{k \in \mathbb{N}}$ is the sequence of sectors of a cutting sequence $c(\tau)$ if and only if the following hold:
\begin{itemize}
\item[(i)] $s \in S^*$ (see  Definition \ref{Sstardef} for the definition of $S^*$);
\item[(i)] if for some $k_0\in \mathbb{N}$, $s_{k}=1$ for all $k >  k_0$ (eventually $1$), then  $s_{k_0} \in \{ 1, 3,5,7\}$;
\item[(ii)]  if for some $k_0\in \mathbb{N}$, $s_{k}=7$ for all $k >  k_0$ (eventually $7$), then  $s_{k_0} \in \{ 2,4,6 \}$;
\end{itemize}
\end{cor}



\subsection{From sequences to directions}
In this section we formulate two results (Theorem \ref{directionsthm} and Proposition \ref{3periodic}) which show that from a cutting sequence $w$, using the renormalization schemes defined in \S\ref{renormalizationsec}, 
  one can recover information about the direction of the trajectories $\tau$ such that $w=c(\tau)$.
\subsubsection{Direction recognition}\label{directionrecognitionsec}
 We assume here that $w$ is a cutting sequence of trajectory $\tau$ in some unknown direction $\theta$ and we would like to determine $\theta$ from the knowledge of the symbolic sequence (this is question Q$2$ in \S\ref{questionssec}). 
In this section we assume that  $w=c(\tau)$ is \emph{non-periodic}.  Thus, by Proposition \ref{admsectorsfornonper}, we can  associate to it an unique sequence $\{d_k(w)\}_{k\in\mathbb{N}}$ of  admissible diagrams   which is determined   through the combinatorial renormalization scheme in \S\ref{combrenormsec}. 

\begin{thm}[Direction recognition]\label{directionsthm}
If $w$ is a \emph{non-periodic} cutting sequence, the direction of trajectories $\tau$ such that $w=c(\tau)$
is uniquely determined and given by 
\bes
\theta=[d_0(w);d_1(w),\ldots, d_k(w), \ldots  ], 
\ees 
where  $\{d_k(w)\}_{k\in \mathbb{N}}$ 
is the sequence of admissible diagrams. 
\end{thm}

\noindent Theorem \ref{directionsthm} is proved in \S\ref{directionsthmproof}. 
From Theorem \ref{directionsthm}, one can get an effective algorithm to determine the direction up to finite precision in finitely many steps.  
If after deriving  a cutting sequence $c(\tau)$ finitely many times, we determine that the sequence of admissible diagrams begins with $d_0, d_1, \dots, d_{k_0} $, then the direction $\theta$ of $\tau$ belongs to the \emph{sector} of possible directions $\overline{\Sigma}(d_0,d_1,\dots,d_{k_0})$.

\subsubsection{Characterization of terminating directions}\label{charperiodic}
Terminating directions can be characterized and recognized through periodicity of the corresponding cutting sequences as follows. 
\begin{prop}\label{3periodic}
Given a trajectory $\tau$ in direction $\theta$, the following properties are equivalent:
\begin{itemize}
\item[(i)] The trajectory $\tau$ is {periodic}; 
\item[(ii)] The cutting sequence $c(\tau)$ is a periodic sequence; 
\item[(iii)] The direction $\theta$ is {terminating}.
\end{itemize}
\end{prop}
\noindent Proposition \ref{3periodic} is proved in \S\ref{3periodicproof}.

In the case of the octagon, terminating directions can also be characterized algebraically as follows.
\begin{thm}\label{Calta}
The direction $\theta$ is terminating if and only if its slope is in the field $ \mathbb{Q}(\sqrt{2})$.
\end{thm}
\noindent This result (mentioned also in \cite{AH:fra}) follows from combining Proposition \ref{3periodic} with the characterization of periodic directions for translation surfaces which have genus two and one singular point. Calta \cite{C:vee} and McMullen \cite{Mc:tei} showed that in appropriate coordinates these are just the directions with slopes in the field  $\mathbb{Q}(\sqrt{2})$  (see also  \cite{HS:inf}).  While many of the results of this paper concerning  the octagon generalize to $2n$-gons, this algebraic characterization of terminating sequences does not hold for general 2n-gons (see Remark \ref{holonomy} in \S\ref{2ngonssec}).

\subsection{Characterization of the closure of cutting sequences}
Unlike the case of the square, the property of  infinite derivability does not suffice to characterize cutting sequences (see Example \ref{incoherentex}, which is given in this section  after Proposition \ref{infcoherentprop}).  In order to characterize cutting sequences we need to require an additional condition. 

\subsubsection{Coherence}
In this section we introduce the notion of \emph{coherence}\footnote{This condition guarantees that an admissible word whose derivatives are all admissible is not obtained by putting together finite words which are all cutting sequences of finite trajectories, but relative to different sectors (as in Example \ref{incoherentex}).}, which  is used to give a characterization of the closure of the set of cutting sequences. 
We stress that in this definition $w$ is any word (not necessarily a cutting sequence). In the following definitions we are taking some care so that the notion of coherence applies to words which do not uniquely define their sector. The definition would be simpler if we were only considering the case in which words do uniquely define their sector.

\begin{defn}[Coherence]\label{coherencedef}
Let us say that a \emph{derivable} word $w$ is \emph{coherent} if there exists a pair $(i,j)$ where $i, j \in \{0,\dots, 7\}$ such that the following conditions $($C$0),($C$1),($C$2),($C$3) $ hold. 
In this case we also say that $w$ is \emph{coherent with respect to} the pair $(i,j)$. 
\begin{itemize}
\item[(C0)] The word $w$ is admissible in diagram $i$;
\end{itemize}
If we normalize $w$ by setting $n(w):=\pi_i \cdot w$, then:
\begin{itemize}
\item[(C1)] The sandwiched letters which occur in $n(w)$ fall into one of the following groups $G_k$:
\vspace{-2mm}
\begin{enumerate}
\item[$G_0$:=] \{$D$-sandwiched $A$, $C$-sandwiched $B$, $B$-sandwiched $C$,  $A$-sandwiched $D$\};\vspace{-2mm}
\item[$G_1$:=] \{$D$-sandwiched $A$, $C$-sandwiched $B$, $B$-sandwiched $C$,  $B$-sandwiched $D$\};\vspace{-2mm}
\item[$G_2$:=] \{$D$-sandwiched $A$, $C$-sandwiched $B$, $C$-sandwiched $C$,  $B$-sandwiched $D$\};\vspace{-2mm}
\item[$G_3$:=] \{$D$-sandwiched $A$, $D$-sandwiched $B$, $C$-sandwiched $C$,  $B$-sandwiched $D$\}.
\end{enumerate}
\item[(C2)] The derived sequence  $n(w)'$ is admissible in a diagram $j \in \{1, \dots, 7\}$;
\item[(C3)] The indices $k$ and $j$ (defined in (C1) and (C2) respectively)  are related by the formula $k=[j/2]$, where $[j/2]$ denotes the integer part of $j/2$. 
\end{itemize}
\end{defn}

\begin{defn}[Renormalized word sequence]\label{worddef}
Given a word $w \in \mathscr{A}^{\mathbb{Z}}$ and a sequence $s\in S^*$ define
a sequence of words $\{w_k\}_{k\in\mathbb{N}}$ by $w_{0}:=w$ and $w_{k} := (\pi_{s_{k-1}} \cdot w_{k-1})'$ for $k>0$.
\label{combrecursive} 
\end{defn}

\begin{defn}[Infinite coherence]\label{infcoherencedef}
Let us say that a word  $w \in \mathscr{A}^{\mathbb{Z}}$ is  \emph{infinitely coherent with respect to a sequence}
 $s \in S^*$  if it is {infinitely derivable} and, for each $k \in \mathbb{N}$, the renormalized word $w_{k}$ in Definition \ref{worddef} is coherent with respect to $(s_k,s_{k+1})$ . We say that $w$ is \emph{infinitely coherent}  if there exists a sequence $s \in S^* $ such that $w$ is infinitely coherent  with respect to $s$.
\end{defn}
\noindent 
The sequence  $s$ in Definition \ref{infcoherencedef} is not necessarily unique.  Example \ref{periodicsameseq} below shows that the first term $s_0$ can have two values. Nevertheless, given the first element $s_0$, the sequence $s$ is essentially unique, as shown by the following Lemma \ref{uniqueness}. 
\begin{lemma}\label{uniqueness} Let $w$ be infinitely coherent and  admissible in diagram $s_0$. 
\begin{itemize}
\item[(i)] Let  $w$ be \emph{non periodic}. Then there is a \emph{unique} sequence $s \in S^*$ starting with $s_0$ with respect to which $w$ is infinitely coherent.
\item[(ii)]  Let  $w$ be  \emph{periodic}. If $w$ is infinitely coherent with respect to  two sequences  $s', s''  \in S^*$ and both start with   $s'_0=s''_0=s_0$, then there exists   $\overline{k}\geq 1$ such that $s'_k = s''_k$  for all $k < \overline{k}$ and either   $s_k '=s_k''= 1$  for all $k > \overline{k}$ or $s_k '=s_k''= 7$  for all $k > \overline{k}$. 
\end{itemize}
\end{lemma}
\noindent  Lemma \ref{uniqueness} is proved in \S\ref{generationandcoherencesec}.  The ambiguity for periodic sequences described in Lemma \ref{uniqueness} is analogous to the  ambiguity of the octagon additive continued fraction expansion of terminating directions in  (\ref{2CFeven}, \ref{2CFodd}), see Lemma \ref{contfracprop2}. 

 Let us recall that by Theorem \ref{cutseqinfderiv} a cutting sequence is infinitely derivable. The following stronger statement holds.  
\begin{prop}\label{infcoherentprop}
A cutting sequence $c(\tau)$ is infinitely coherent with respect to its sequence of sectors $\{s_k(\tau)\}_{k\in \mathbb{N}}$. 
\end{prop}
\noindent The proof of Proposition \ref{infcoherentprop}  is given in \S \ref{infcoherencesec}.  Proposition \ref{infcoherentprop} 
 shows that  being infinitely coherent is a necessary condition for a sequence to be a cutting sequence. 
   In the next section we show that it is close to being sufficient, up to considering the closure of the set of cutting sequences (Theorem \ref{cuttingseqthm}).

We can now use Proposition \ref{infcoherentprop} to exhibit an example of an \emph{infinitely derivable  word} which is \emph{not a cutting sequence} (Example \ref{incoherentex}) and of a word which is infinitely coherent with respect to two sequences $s, s'$ with $s'_0\neq s_0''$ (Example \ref{periodicsameseq}).
\begin{ex}\label{periodicsameseq}
Consider the periodic word  $w= \overline{ADBCBD}$.  One can check that it is the cutting sequence of a periodic trajectory. Thus, by Proposition \ref{infcoherentprop},  it is infinitely coherent.  The word $w$ is admissible both in diagram $0$ and in diagram $4$. We can hence choose either $s'_0=0$ or $s''_0=4$ and then complete $s'_0$ to a sequence $s'$ and $s_0''$ to $s''$, so that $w$ is infinitely coherent with respect to both $s'$ and $s''$.
\end{ex}
 
\begin{ex}\label{incoherentex}
 Consider the following periodic word\footnote{This word $w$ and other similar examples can be constructed by mixing generation rules relative to different graphs in Figure \ref{generationrulesfig}. Here $w$ is obtained from $w'=\overline{CDDDCA}$ by interpolating the transitions using generation rules both from  $\mathscr{D}_6$ and  $\mathscr{D}_2 $ in Figure \ref{generationrulesfig}.} :
\bes w =  \overline{CCCBDBCCBDBCCBDBCBDADB}.      
\ees
One can check that $w$ is admissible in diagram $0$, so that $n(w)=w$, and that the derived sequence is $w'=  \overline{CDDDCA}$ which is admissible in diagram $6$. Since the successive derivatives are all equal to $\overline{DA}$ (which  is admissible up to any relabelling) $w$ is infinitely derivable. Let us remark though that the first sandwiched $C$ is $C$-sandwiched while the second sandwiched $C$  is $B$-sandwiched, so that  $w=n(w)$ does not satisfy condition $($C$1)$ in  Definition \ref{coherencedef} and thus cannot be coherent. This implies, by Proposition  \ref{infcoherentprop}, that $w$ cannot be a cutting sequence. 
\end{ex}

\subsubsection{Characterization of the closure via coherence}
Let us recall that the space $\mathscr{A}^{\mathbb{Z}}$ 
 has a natural topology which makes it a compact space (we refer e.g.~to  \cite{LM:sym}). 
The characterization of the closure of cutting sequences is the following. 
\begin{thm}[Characterization of the closure]\label{cuttingseqthm}
The closure of the space of cutting sequences coincides with the set of infinitely coherent sequences.
\end{thm}
\noindent The proof of Theorem \ref{cuttingseqthm} is given in \S \ref{proofsfindtrajectoiriessec}. An equivalent and more constructive formulation of this characterization 
 is given in Proposition \ref{generationandinfcoherence} in  \S\ref{generationsec}.

It is possible to describe explicitly the sequences which belong to the closure but are not themselves cutting sequences. The situation is very parallel to the case of Sturmian sequences (see \S\ref{torussec}). 
 Infinitely coherent sequences which are not cutting sequences have one of the two following forms.
\begin{enumerate}
\item[(i)]
If $\theta$ is terminating, they are obtained juxtaposing two periodic sequences, possibly with a letter in between them. An example of this type is given by the sequence  
\bes \dots ADADAD B CCCCCC\dots \ees
\item[(ii)]
If $\theta$ is non-terminating, they are obtained juxtaposing the two semi-infinite cutting sequences corresponding to two rays (i.e.~semi-infinite trajectories starting at a vertex), possibly with a word in between coding sides which enter the common endpoint. 
\end{enumerate}
In a forthcoming paper we give a geometric interpretation  of the closure and prove that  these two types of examples are the only possible ones. 

\subsubsection{Generation operators}\label{generationdefsec}
We will introduce a combinatorial operation on sequences which will serve as an inverse to the derivation operation. We call this operation \emph{generation}. Specifically generation allows to us invert the operation of derivation if we know the admissible diagrams (see Lemma \ref{g0invertgeneration}). We use it in \S\ref{generationsec} to give an alternative characterization of the closure of the set of cutting sequences (Proposition \ref{generationandinfcoherence}) and  to describe in \S\ref{algorithmsec} an algorithm  which produces all finite cutting sequences in a given direction $\theta$ (or in a  sector of directions). For a geometric interpretation of this  generation operation, we refer to Proposition \ref{generatetype0} in \S\ref{generationoperationsec}.

\paragraph{Interpolating words.}
Consider Figure \ref{generationrulesfig}. On each edge of the labelled diagrams $\mathscr{D}_k$, $k=1,\dots ,7$ 
 there appears a word in the letters $\{A,B,C,D\}$. Let us denote by $w^k_{L_1L_2}$ the word associated to the edge from $L_1$ to $L_2$ on $\mathscr{D}_k$. 
We refer to $w^k_{L_1L_2}$ as an \emph{interpolating word} or more precisely as an interpolating word in diagram $k$ and we call  \emph{generation rules} (or generation rules in diagram $k$) the collection of all interpolating words (in diagram $k$). 
\begin{figure}[!h]
\centering
\subfigure[$\mathscr{D}_1$]{\includegraphics[width=0.36\textwidth]{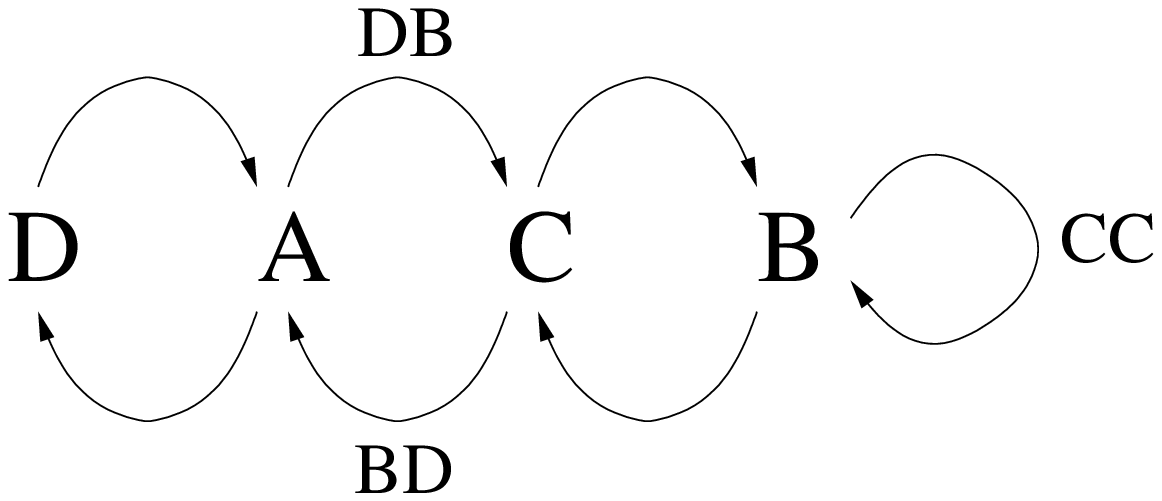}}
\hspace{9mm}
\subfigure[$\mathscr{D}_2$]{\includegraphics[width=0.36\textwidth]{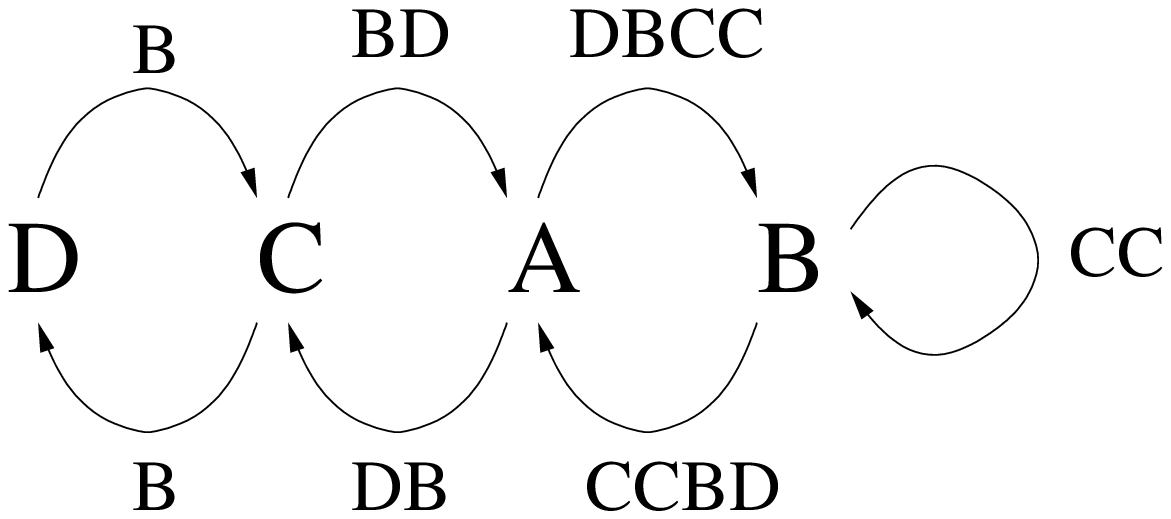}}
\subfigure[$\mathscr{D}_3$\label{generationrulesfig3}]{\includegraphics[width=0.39\textwidth]{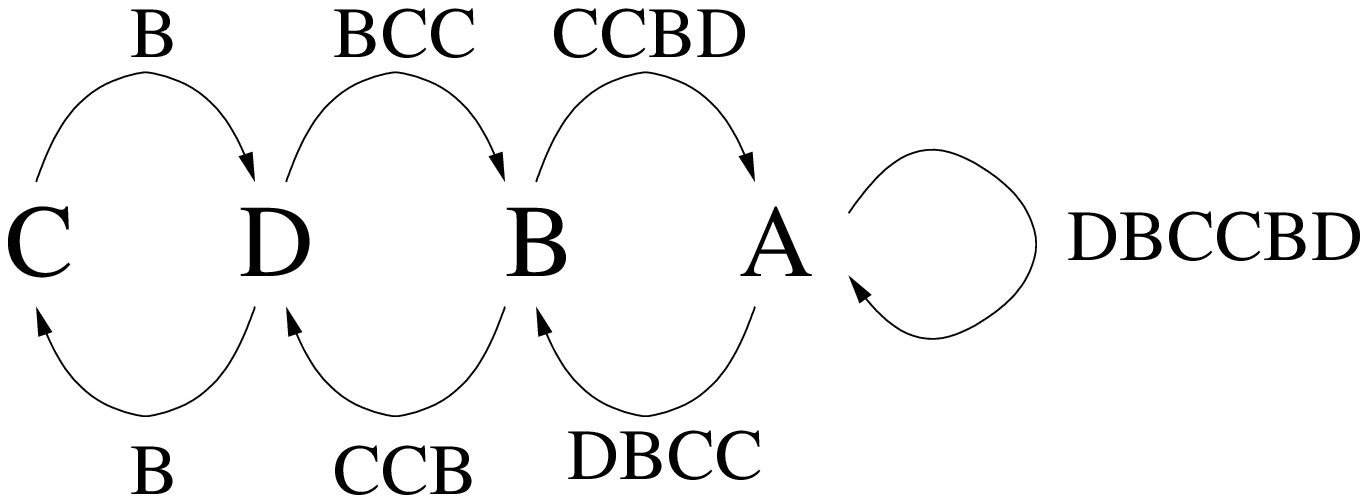}}
\hspace{9mm}
\subfigure[$\mathscr{D}_4$]{\includegraphics[width=0.36\textwidth]{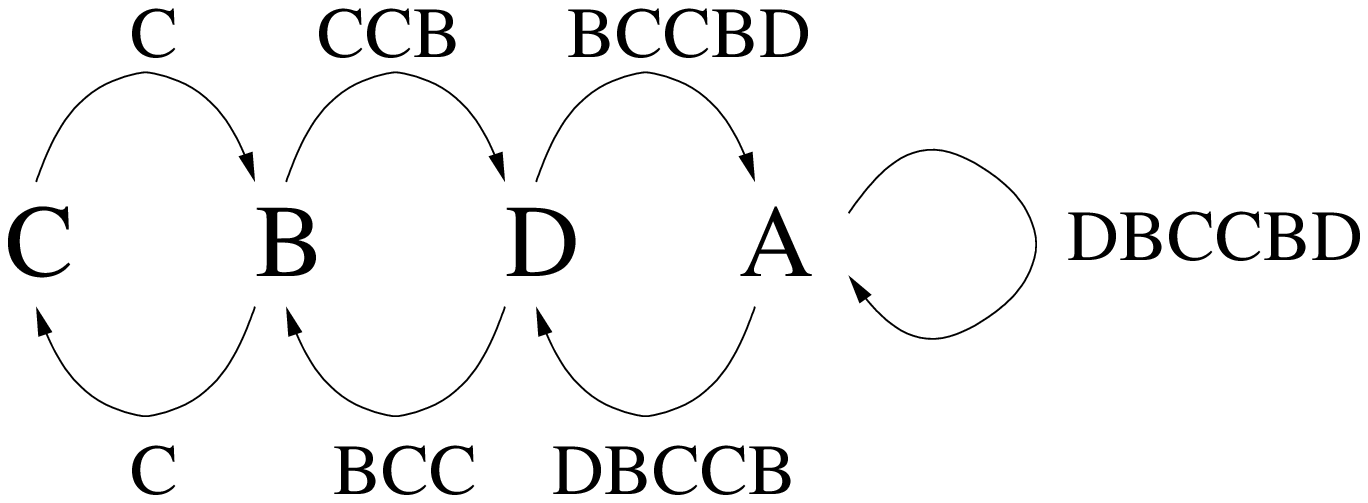}}
\subfigure[$\mathscr{D}_5$]{\includegraphics[width=0.36\textwidth]{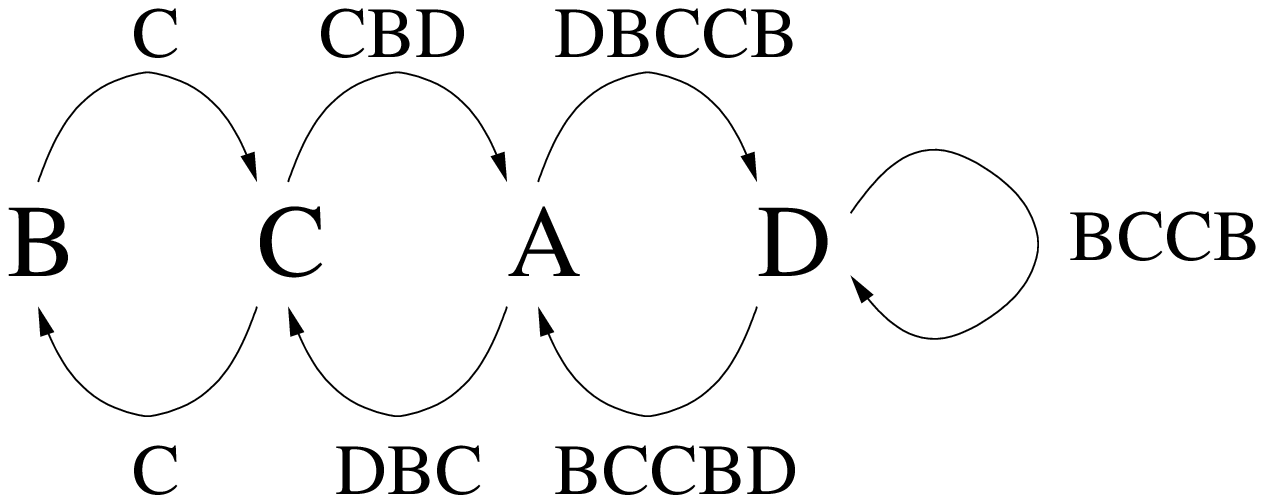}}
\hspace{9mm}
\subfigure[$\mathscr{D}_6$]{\includegraphics[width=0.36\textwidth]{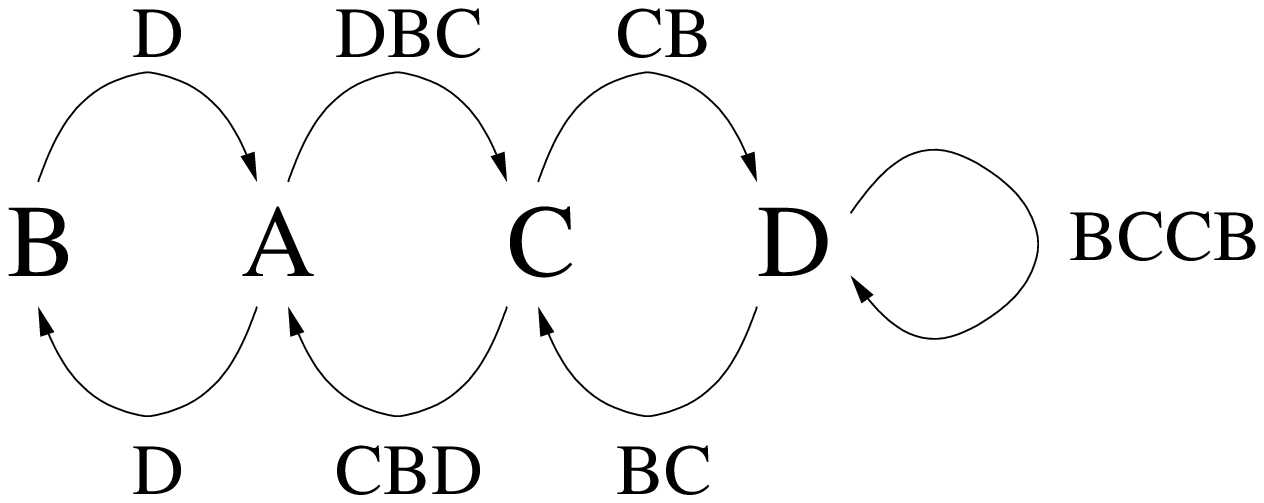}}
\subfigure[$\mathscr{D}_7$]{\includegraphics[width=0.36\textwidth]{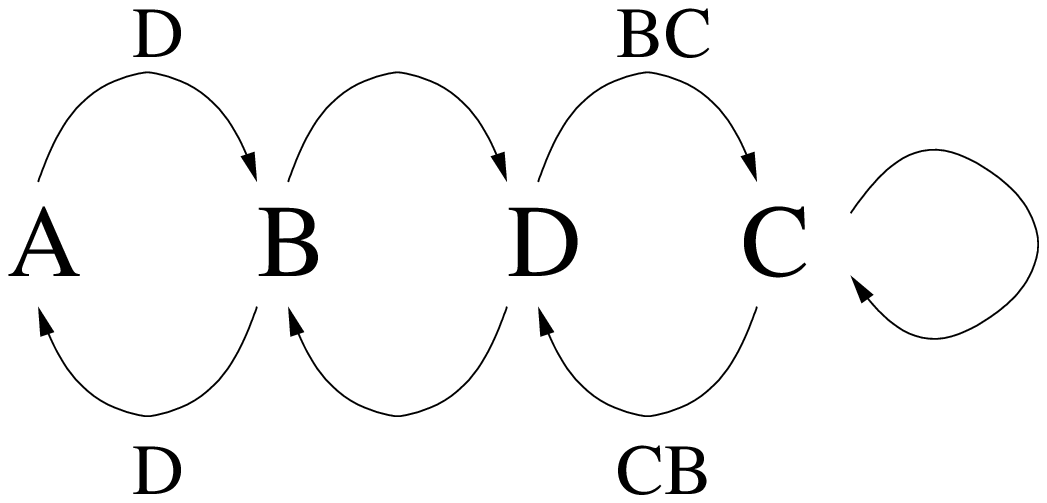}}
\caption{Generation rules to generate sequence $W$ of type  $0$. \label{generationrulesfig}}
\end{figure}

For each each $k=1,\ldots ,7$, we define an operator $\gop{k}{0}$ that acts on finite and infinite words admissible in diagram $k$ and produces words admissible in diagram $0$. 
\begin{defn}[Generation operations to diagram $0$]\label{generationopto0}
If $w$ is a finite or infinite word  admissible in diagram $k$, let $\g{k}{0}{w}$ be the word which is obtained by inserting between successive pairs of letters in  $w$ the corresponding interpolating words in diagram $k$, so that if $w = \ldots L_1L_2L_3
\ldots$ then,
\bes
\g{k}{0}{w} =\ldots L_{1} \,\, w^k_{L_{1} \,\, L_{2}}\,\,  L_2 \,\, w^k_{L_{2}L_{3}}\,\,  L_{3} 
 \ldots.
\ees
\noindent In other words, since $w$ is admissible in diagram $k$ it defines an infinite path on $\mathscr{D}_k$, the word $\g{k}{0}{w}$ is obtained by concatenating the sequences of letters associated to both the vertices and the edges that this path traverses. 
\end{defn}
\begin{ex}\label{generationex}
Let $u=CDBAABDBD$. By definition, since $u$ is a finite word  
admissible in diagram $3$, it gives a finite path through the corresponding vertices of $\mathscr{D}_3$ in Figure \ref{generationrulesfig3}. By inserting into $u$  the interpolating words in diagram $3$ given  in Figure \ref{generationrulesfig3}, 
we read off $\g{3}{0}{u}$, which is the following word (letters corresponding to vertices are written in bold font for clarity):
$$   \g{3}{0}{u}= \mathbf{C}\, B\, \mathbf{D}\, BCC \,\mathbf{B}\, CCBD\, \mathbf{A}\,DBCCBD\,\mathbf{A}\,DBCC\,\mathbf{B}\,CCB\,\mathbf{D}\,BCC\,\mathbf{B}\,CCB\,\mathbf{D}. $$ 
\end{ex}
The sequences produced in this way 
 invert derivation in the following sense.  
\begin{lemma}\label{g0invertgeneration}
If $w$ is admissible in diagram $k$ then $W:= \g{k}{0}{w} $ is admissible in diagram $0$ and satisfies the equation $W'=w$. 
\end{lemma}
\begin{proof}
It is enough to check that, for any vertex $v$ of a diagram in Figure \ref{generationrulesfig}   and any pair of edges $e_1$ and $e_2$ respectively entering and exiting $v$, if  $v$ is labelled by the letter $L$,   $w_{L_1 L}$ is the interpolating word in Figure \ref{generationrulesfig} corresponding to $e_1$  and  $w_{LL_2}$ the one associated to $e_2$,  the  concatenation $w_{L_1L} L w_{L L_2}$ 
is admissible in diagram $0$ and moreover  the only sandwiched letter in $w_{L_1L} L w_{L L_2}$   is $L$. The latter property implies that $ (\g{k}{0}{w})'=w$.  This is verified by inspecting the diagrams in Figure \ref{generationrulesfig} and we leave this to the reader.
\end{proof}
\noindent For the  connection between generation operations with the inverse branches of derivation, see also Proposition \ref{generatetype0} and Corollary \ref{7to1} in \S\ref{generationoperationsec}.
The proof of Proposition \ref{generatetype0} in \S\ref{generationoperationsec} gives an insight into how one can explicitly find the interpolating words geometrically in the case in which $w$ is a cutting sequence. 

The operator  $\g{j}{0}{w}$ produces words admissible in diagram $0$.
We now define operators that  generate words admissible in the other diagrams. These are easily obtained from $\g{j}{0}{w}$ by appropriate relabelling.

\begin{defn}\label{othergoperators}
For $j=1,\dots, 7$ and $i=0,\dots, 7$, let $\gop{i}{i}$ be the operator that acts on finite or infinite words admissible in diagram $j$, as follows:
\bes
\g{j}{i}{w}:= \pi_i^{-1}\cdot \g{j}{0}{w}.
\ees
\end{defn}
\noindent
Let us remark that the word  $\g{j}{i}{w}$ is admissible in diagram $i$ since it is obtained by relabelling a word admissible in diagram $0$ by an appropriate permutation.
Note that $\gop {j}{i} \g{k}{j}{w}$ is not defined when $j=0$ so that the composition $\gop {j}{i} \g{k}{j}{w} $  is only well defined when $j>0$.

\subsubsection{Characterization of the closure via generation}\label{generationsec}
The following Lemma shows the relation between generation and coherence.
\begin{lemma}[Coherence via generation]\label{coherencegeneration}
The word $w$ is coherent with respect to $(i,j)$ if and only if there exists a word $v$ admissible in diagram $j$ such that $w=\g{j}{i}{v}$. In this case, we have $v=n(w)'$. 
\end{lemma} 
\noindent The proof is given in \S\ref{generationandcoherencesec}.  The Lemma shows that the conditions $($C$1)$ and $($C$3)$ in the definition of coherence (Definition \ref{coherencedef}) are enough  to recover the  interpolating words. 
Thus the large collection of interpolation words shown in Figure  \ref{generationrulesfig} can be computed from the small number of coherence rules in Definition \ref{coherencedef}. Let us introduce more notation in order to reformulate  infinite coherence in terms of generation operators. 
\begin{defn}
Let us denote by $\mathscr{A}_k$ the set of all words which are admissible in diagram $k$.
Let us denote by $\mathscr{P}_k$ the set of all periodic words of period $1$ or $2$ which are admissible in diagram $k$ (see Example \ref{exgeneration}). 
\end{defn}
\begin{rem}\label{Pdescription} The set $\mathscr{P}_k$ is a finite set containing the $4$  words admissible in one of the two diagrams  $\mathscr{D}_{k-1,k}$ and  $\mathscr{D}_{k,k+1}$ (where $k$ should be considered modulo $7$) corresponding to $\pi k/8$ and $\pi (k+1)/8$ (see Figure \ref{boundarydiagrams}). 
\end{rem}
\begin{ex}\label{exgeneration}
Let $k=6$. Let us list the elements of $\mathscr{P}_6$. 
\bes
\mathscr{P}_6 = \{ \overline{BA}, \overline{AC}, \overline{CD}, \overline{D} \}.
\ees
One can check that this is the union of the cutting sequences corresponding to $6\pi /8$ and $7 \pi/8$ or equivalently of the infinite paths on the diagrams $\mathscr{D}_{5,6}$ and $\mathscr{D}_{6,7}$ in Figures \ref{D56fig}, \ref{D67fig}. 
\end{ex}
\begin{defn}
Given  $s_0 \in \{ 0,\dots, 7\}$ and $s_1, \dots, s_k$ in $\{ 1,\dots, 7\}$, let $ \mathscr{G}(s_1,\dots, s_k)$ and  $ \mathscr{P}(s_1,\dots, s_k)$ be  respectively the collection of words  of the form:
\begin{eqnarray}
\mathscr{G}(s_0,\dots, s_k) & := & \left\{ \,  \g{s_1}{s_0}{ \g{s_2}{s_1}{ \ldots \g{s_{k-1}}{s_{k-2}}{ \g{s_k}{s_{k-1}}{u}}}},  \,   \quad u \in \mathscr{A}_{s_k}  \,  \right\} ; \label{Gsetsdef} \\
\mathscr{P}(s_0,\dots, s_k) & := & \left\{ \, \g{s_1}{s_0}{ \g{s_2}{s_1}{ \ldots  \g{s_{k-1}}{s_{k-2}}{ \g{s_k}{s_{k-1}}{u}}}}, \, \quad u \in \mathscr{P}_{s_k} \,  \right\}. \label{periodiccollection}
\end{eqnarray}
\noindent Clearly, since the generation operators preserve periodicity, the words in $\mathscr{P}(s_0,\dots, s_k)$ are all periodic. One can show the following stronger fact:
\begin{lemma}\label{Pperiodic}
For each  $s \in S^{*}$ and for each $k$, all words in $\mathscr{P}(s_0,\dots, s_k)$ are \emph{cutting sequences} of periodic trajectories.
\end{lemma}
\noindent Lemma \ref{Pperiodic} is proved in \ref{proofsfindtrajectoiriessec}.

\end{defn}
\begin{ex}\label{gsetex}
Let us construct $\mathscr{P}(0,1,6)$. The set  $\mathscr{P}_6$ is given in Example \ref{exgeneration}.  In order to get $\mathscr{P}(1,6)$,   let us apply $\gop{6}{1}$, by first applying $\gop{6}{0}$ and then acting by $ \pi_1^{-1}$.  Applying $\gop{6}{0}$ to each element of $\mathscr{P}_6$, (using the labels on the edges on $\mathscr{D}_6$ in Figure \ref{generationrulesfig}),  we obtain 
\bes
\mathscr{P}(0,6) = \{  \overline{BDAD},\overline{ADBCCCBD}, \overline{CCBDBC}, \overline{DBCCB} \}
\ees
and acting by $\pi_1^{-1} = (AD)(CB)$ we get
\bes \mathscr{P}(1,6) = \{  \overline{CADA},\overline{DACBBBCA}, \overline{BBCACB}, \overline{ACBBC}  \}.
\ees
To get  $\mathscr{P}(0,1,6)$, one applies $\gop{1}{0}$, i.e.~interpolates with the generation rules on  $\mathscr{D}_1$ in Figure \ref{generationrulesfig}:
\bes \begin{split} \mathscr{P}(0,1,6)  =& \{  \overline{CBDADADB},\overline{DADBCBCCBCCBCBDA},\\ & \overline{BCCBCBDADBCBCC}, \overline{ADBCBCCBCBD}  \}.
\end{split}
\ees
\end{ex}

The following Proposition gives an equivalent characterization of infinitely coherent words  and thus, by Theorem \ref{cuttingseqthm}, an equivalent characterization of the closure of the set of cutting sequences. 
\begin{prop}[Infinite coherence via generation]\label{generationandinfcoherence}
A word $w$ is infinitely coherent with respect to $s \in  S^* $ if and only if 
\be \label{infiniteintersection}
w \in \bigcap_{k\in\mathbb{N}} \mathscr{G}_k ^s, \qquad \mathscr{G}_k^s:= \mathscr{G}(s_0,\dots, s_k). 
\ee
Thus, a word $w$ is infinitely coherent if and only if 
\be \label{unionfiniteintersection}
w \in \bigcap_{k\in\mathbb{N}} \, \, \bigcup_{\substack{ 0\leq s_0 \leq 7 \\ 1\leq s_1, \dots, s_k   \leq 7} }  \mathscr{G}(s_0,\dots, s_k). 
\ee
\end{prop}
\noindent The characterization (\ref{unionfiniteintersection}) is useful in the proof of Theorem \ref{cuttingseqthm} to see that infinitely coherent words form a closed set. The proof of Proposition \ref{generationandinfcoherence} is given in \S \ref{generationandcoherencesec}. 

\subsubsection{Approximation by periodic cutting sequences}\label{algorithmperiodicsec}
Infinitely coherent words can be approximated by periodic cutting sequences, as shown by the following Proposition \ref{finiterealization} and its Corollary \ref{approxbyperiodic}. 
For brevity, let us say that a finite word $u $ in the letters of $\mathscr{A}$ is \emph{realized} by an infinite trajectory $\tau$, or that $\tau$ \emph{realizes}  $u$, if there exists a finite piece $\tau'\subset \tau$ such that $u=c(\tau')$ is its cutting sequence.
\begin{prop}\label{finiterealization}
If $w$ is an infinitely coherent word, each finite subword of $w$ is realized by a periodic trajectory $\tau$.  In particular, if $w$ is infinitely coherent with respect to $\{s_k\}_{k\in\mathbb{N}}$,  if each finite subword of $w$ occurs as a subword of  some 
periodic cutting sequence in  
\be
 \bigcup_{k\in\mathbb{N}} \mathscr{P}_k , \qquad \mathscr{P}_k:= \mathscr{P}(s_0,\dots, s_k). \label{periodicunion}
\ee
\end{prop}
\noindent If ${S}$ is a set in $\mathscr{A}^{\mathbb{Z}}$, let denote by $\overline{S}$ its closure. 
\begin{cor}\label{approxbyperiodic} If  the word $w$ is infinitely coherent with respect to $\{s_k\}_{k\in \mathbb{N}}$, then
\be\label{infiniteintersectionperiodic}
w \in \overline{\bigcup_{k\in\mathbb{N}}  \mathscr{P}_k} , \qquad \mathscr{P}_k:= \mathscr{P}(s_0,\dots, s_k). 
\ee 
\end{cor}
\noindent Proposition \ref{finiterealization}  and its Corollary \ref{approxbyperiodic} are proved in \S\ref{proofsfindtrajectoiriessec}. 

\subsubsection{An algorithm to generate cutting sequences}\label{algorithmsec}
Let us conclude by describing how to generate all finite subwords of  cutting sequences of linear trajectories in direction $\theta$.
\begin{enumerate}
\item Given $\theta$, iterate the octagon Farey map $F$ to generate its  itinerary $\{s_k \}_{k\in\mathbb{N}}$, where $s_k$ are such that $F^k (\theta) \in \Sigma_{s_k}$, for all $k\in \mathbb{N}$;
\item For each $k$, list the elements of the set $\mathscr{P}_k:= \mathscr{P} \left( s_0, s_1,  \dots, s_k \right)$ (as in Example \ref{gsetex}); 
\item Each finite subword of a cutting sequence of a linear trajectory in direction $\theta$ appears as a subword of a word in $\mathscr{P}_k$ if $k $ is sufficiently large. 
\end{enumerate}
The above procedure can be made into an effective algorithm for finding words of fixed finite length which can occur in a cutting sequence in direction $\theta$. It is enough to iterate each step only a finite number of times, for some sufficiently large  $k$ that can be specified. 
Moreover, one can also obtain all words which occur in a sector of directions by considering all corresponding continued fraction expansions (for example for a  sector of the form $\overline{\Sigma}[s_0; s_1,\dots, s_k]$, one has to consider all expansion which share the first $k$ entries). 

The proof that the algorithm actually gives all finite subwords is obtained combaining the following Propositions. Cutting sequences are infinitely coherent with respect to their sequence of sectors (Proposition \ref{infcoherentprop}) and thus, by Proposition \ref{finiterealization}, they can be approximated through periodic cutting sequences of the form  (\ref{periodiccollection}), where the sequence $s$ coincides with the itinerary under $F$ 
by   Proposition \ref{typesequalCF}. 

\section{Basic definitions and octagon affine automorphisms group}\label{defsec}
In this section we give the definitions of translation surfaces, the $GL(2,\mathbb{R})$ action, the affine automorphism group of a translation surface and the Veech group. 
 Our approach has the novel feature that we allow orientation reversing affine automorphisms. We are lead to do this because orientation reversing affine automorphisms play and important role in the renormalization schemes we defined in \S \ref{renormalizationsec}.

\subsection{Basic definitions of the $GL(2,\mathbb{R})$ action, $\Aff(S)$ and $V(S)$.}\label{defnsec}
\paragraph{Translation surfaces and linear trajectories.}
Let us define a translation surface to be a collection of polygons $P_j\subset\R^2$ with identifications of pairs of parallel sides so that (1) sides are identified by maps which are restrictions of translations, (2) every side is identified to some other side and (3) when two sides are identified the outward pointing normals point in opposite directions. If $\sim$ denotes the equivalence relation coming from identification of sides then we define the surface $S=\bigcup P_j/\sim$. The surface  $S_O$  obtained by gluing opposite pairs of parallel sides of the octagon $O$  is an example of a translation surface. Let us stress that when we refer to a translation surface $S$, we will have in mind a particular collection of polygons in $\mathbb{R}^2$ with identifications. 

Translation surfaces can also be defined in terms of special atlases.
A \emph{translation atlas} for a (punctured) surface $X$ is a collection  $\{(U_\alpha, \phi_{\alpha})\}$ where the sets $U_\alpha$ are open sets in $X$,  and $\phi_{\alpha} : U_\alpha \rightarrow \R^2$ are charts such that the transition maps $ \phi_{\beta}\phi_{\alpha}^{-1}$ are restrictions of translations $v \rightarrow v +b$ from $\mathbb{R}^2$ to $\mathbb{R}^2$. If $S$ is a translation surface and $\Sing$ is the set of points corresponding to vertices of polygons then $X= S-\Sing$ has a translation atlas. 
The cone angle at a point in $\Sing$ is the sum of the angles at the corresponding points in the polygons $P_j$. We say that such a point is singular if the cone angle is greater than $2\pi$. For example, the surface $S_O$ has genus 2 and a single singular point with a cone angle of $6\pi$. 


If $S$ is a translation surface then the tangent spaces $T_p(S)$ are all identified with $\R^2$ by means of the coordinate charts.  Any translation invariant geometric structure on $\R^2$ can be transported to all of $S - \Sing$. A vector $v\in\R^2$ gives a \emph{parallel vector field},  a linear functional on $\R^2$ gives a \emph{parallel one form}. The metric $ds^2=dx^2+dy^2$ gives a \emph{parallel metric}. 

Due to the presence of singular points parallel vector fields do not define flows. Nevertheless we can speak of trajectories for these vector field which we call \emph{linear trajectories}.  
Since we will be considering translation surfaces which are given by specific collections of polygons in $\R^2$ we will think of trajectories as a families of line segments in $\R^2$.
We call a  trajectory which does not hit singular points  \emph{bi-infinite}. A trajectory which begins and ends at a singular point is a \emph{saddle connection}. Every periodic trajectory is contained in a maximal family of parallel periodic trajectory of the same period. If the surface is not a torus then this family fills out a \emph{cylinder} bounded by saddle connections.

\paragraph{Affine deformations and affine diffeomorphisms.}
Let $S$ be a translation surface with translation atlas $\{(U_\alpha, \phi_{\alpha})\}$.
Given  $\nu \in GL(2,\mathbb{R})$, the translation surface $S'=\nu \cdot S$ is given by the atlas  $(U_\alpha,  \nu\cdot\phi_{\alpha})_\alpha$ obtained by post-composing the charts with $\nu$. 
If $S =\bigcup P_j/\sim$  then $\nu$ acts linearly on each polygon and $\nu$ takes pairs of parallel sides to pairs of parallel sides. 
 \begin{defn}\label{actiononpolygons}
We denote by $\nu  P \subset \mathbb{R}^2$ the image of a polygon $P \subset \mathbb{R}^2$ under the linear map $\nu$.  
\end{defn}
  The surface  $\nu \cdot S$ is obtained by glueing the corresponding sides of $\nu \octcdot P_1, \dots, \nu \octcdot P_n $. There is a canonical map $\Phi_\nu$ from the surface $S$ to the surface $\nu \cdot S$ which is given by the restriction of the linear map $\nu$ to the polygons $P_1 , \dots ,  P_n$.
Clearly $\Phi_\nu$ sends infinite linear trajectories on $S$ to infinite linear trajectories on $\nu \cdot S$.
\begin{defn}
Given a linear trajectory $\tau$ on $S$ and $\Phi_{\nu}: S \rightarrow \nu\cdot S$,   
we  denote by $\PhiTraj{\nu}{\tau}$ the linear trajectory on $\nu\cdot S$ which is obtained by composing $\tau$ with $\Phi_{\nu}$.
 \end{defn}

Let $S$ and $S'$ be translation surfaces.  Consider a homeomorphism $\Psi$ from $S$ to $S'$ which takes $\Sing$ to $\Sing'$ and is a diffeomorphism outside of $\Sing$. We can identify the derivative $D\Psi_p$ with an element of $GL(2,\mathbb{R})$. We say that $\Psi$ is an \emph{affine diffeomorphism} if the derivative $D\Psi_p$ does not depend on $p$. In this case we write $D\Psi$ for $D\Psi_p$.
The canonical map $\Phi_\nu$ from $S$ to  $\nu \cdot S$ described above is an example of an affine diffeomorphism. In this case $D\Phi_\nu=\nu$.

We say that $S$ and $S'$ are \emph{affinely equivalent} if there is an affine diffeomorphism $\Psi$ between them.  
We say that $S$ and $S'$ are \emph{isometric} if they are affinely equivalent with $D\Psi\in O(2)$. We say that $S$ and $S'$ are \emph{translation equivalent} if they are affinely equivalent with $D\Psi=Id$.  If $S$ is given by identifying sides of polygons $P_j$ and $S'$ is given by identifying sides of polygons $P'_k$  then a translation equivalence $\Upsilon$ from $S$ to $S'$ can be given by a ``\emph{cutting and pasting}" map. That is to say we can subdivide the polygons $P_j$ into smaller polygons and define a map $\Upsilon$ so that  the restriction of $\Upsilon$ to each of these smaller polygons is a translation and the image of $\Upsilon$ is the collection of polygons $P'_k$. 


An affine diffeomorphism from $S$ to itself is an \emph{affine automorphism}. The collection of affine diffeomorphisms is a group which we denote by $\Aff(S)$. The collection of isometries of $S$ is a finite subgroup of $\Aff(S)$ and  the collection of translation equivalences is a subgroup of the group of  isometries. 
If $S$ is given as a collection of polygons with identifications then we can realize an affine automorphism of $S$ with derivative $\nu$ as a composition of a map $\Psi_\nu:S\to\nu\cdot S$ with a translation equivalence, or cutting and pasting map, $\Upsilon:\nu\cdot S \to S$. 

\paragraph{Veech group and lattice surfaces.}
The Veech homomorphism is the homomorphism $\Psi\mapsto D\Psi$  from $\Aff(S)$ to $GL(2,\R)$. The image of this homomorphism lies in the subgroup of matrices with determinant $\pm1$ which we write as $SL_{\pm}(2,\mathbb{R})$. We call \emph{Veech group} and we denote by $V(S)$ the image of $\Aff(S)$ under the Veech homomorphism.
It is common to restrict to orientation preserving affine diffeomorphisms in defining the Veech group. We write $V^+(S)$ for the image of the group of \emph{orientation preserving} affine automorphisms in $SL(2,\R)$.  Since we will make essential use of orientation reversing affine automorphisms we will use the term Veech group for the larger group $V(S)$. Note that the term Veech group is used by some authors to refer to the image of the the group of orientation preserving affine automorphisms in the projective group $PSL(2,\R)$. We denote the image of $V(S)$ in $PGL(2,\R)$ by $V_P(S)$. We denote the image of $V^+(S)$ in $PSL(2,\R)$ by $V^+_P(S)$.

The kernel of the Veech homomorphism is a the finite group of translation equivalences of $S$.  This kernel is trivial if and only if an affine automorphism of $S$ is determined by its derivative. If $S$ is given by identifying polygons and the the kernel of the Veech homomorphism is trivial then the map $\Upsilon$ above is uniquely determined.

A translation surface $S$ is called a \emph{lattice surface} if  $V^+(S)$ is a lattice in $SL(2, \mathbb{R})$. This is equivalent to saying that $V(S)$ is a lattice in $SL_{\pm}(2, \mathbb{R})$. The torus $T^2=\mathbb{R}^2 / \mathbb{Z}^2$ is an example of a lattice surface whose  Veech group is $GL(2, \mathbb{Z})$. Veech proved  more generally that all  translation surfaces obtained from regular polygons are lattice surfaces (see \S \ref{veechsec}).
Lattice surfaces satisfy the \emph{Veech dichotomy} (see Veech \cite{Ve:tei}, Vorobets \cite{Vor:pla}) which says that if we consider a direction $\theta$ then one of the following two possibilities holds: either there is a saddle connection in direction $\theta$ and the surface decomposes as a finite union of cylinders each of which is a union of a family of closed geodesics in direction $\theta$, or each trajectory in direction $\theta$ is dense and uniformly distributed.

\subsection{The Veech group  and the affine automorphism group of $S_O$} \label{veechsec}
The surface $S_O$ obtained by identifying opposite sides of the octagon $O$ is closely related to the surface that Veech discusses in \cite{Ve:tei} but is not the same. Veech considers the double octagon surface obtained by starting with two copies of the octagon and identifying each side on the first copy with a parallel side on the second copy. This surface has a translation equivalence of order two which interchanges the two octagons and acts freely. The quotient of this surface by this translation equivalence is the surface $S_O$. These two surfaces have the same Veech group (see \cite{HS:inv}, Lemma J).  Veech gives elements in $V^+_P$ of the double octagon which generate a lattice group. A priori the Veech group could be larger than the group generated by these elements but Earle and Gardner in \cite{CG:tei} show that these same elements are generators for $V^+_P(S_O)$.
We will summarize these results here. It will be useful for us to know the group $\Aff(S_O)$ as well as the Veech group, since  affine automorphisms act on trajectories and this action plays a key role in the construction of our renormalization schemes.

\paragraph{Some elements of $\Aff(S_O)$.}
Since we allow orientation reversing transformations, the entire isometry group $\Die_8$ of the octagon $O$ is contained in $\Aff(S_O)$. For $\eta\in D_8$ let $\Psi_\eta:S_O\to S_O$ denote the corresponding affine automorphism. We have $D\Psi_\eta=\eta$. Consider in particular  the reflection $\Psi_\alpha $ of the octagon at the horizontal axes and  the reflection $\Psi_\beta $ in the tilted line which forms an angle $\pi/8$ with the horizontal axis. These are given by the two matrices
\be \label{alphabetadef}
\alpha := \begin{pmatrix} 1 & 0 \\  0 & -1  \end{pmatrix} \qquad  \beta := \begin{pmatrix} \frac{\sqrt{2}}{2} & \frac{\sqrt{2}}{2} \\ - \frac{\sqrt{2}}{2} & \frac{\sqrt{2}}{2}  \end{pmatrix}   = \rho_{-\frac{\pi}{8}}  \cdot \alpha \cdot \rho_{\frac{\pi}{8}}
 \ee 
  where $\rho_\theta$  is the matrix representing counterclockwise rotation by the angle $\theta$. 
The elements $\nu_j$, $j=0,\dots, 7$ defined in (\ref{nujdef}) can be expressed as words in $\alpha$ and $\beta$. 
We remark also that $\alpha \beta =  \rho_{\frac{\pi}{4}} $, which is an orientation preserving element and hence  $\Psi_\alpha\Psi_\beta\in V^+(S_O)$. Note that the $(\alpha \beta)^4=-Id$ so that the element $(\Psi_\alpha\Psi_\beta)^4$ is in the kernel of the homomorphism from $\Aff(S_O)$ to $V_P(S_O)$. The image of $\Die_8$ in $V_P(S_O)$  is isomorphic to $\Die_4$.

An affine automorphism of a translation surface is \emph{parabolic} if its derivative is a shear which is to say that  it has two equal eigenvalues but is not equal to the identity. Parabolic automorphisms of translation surfaces are closely related to cylinder decompositions. Every parabolic automorphism gives a cylinder decomposition and if we have a cylinder decomposition into cylinders for which the inverse moduli are commensurable then there is a parabolic automorphism which acts as a multiple of a Dehn twist in each cylinder. 

\begin{figure}[!h]
\centering

\subfigure[Horizonatal decomposition]{\label{cylinders1fig}\includegraphics[width=0.37\textwidth]{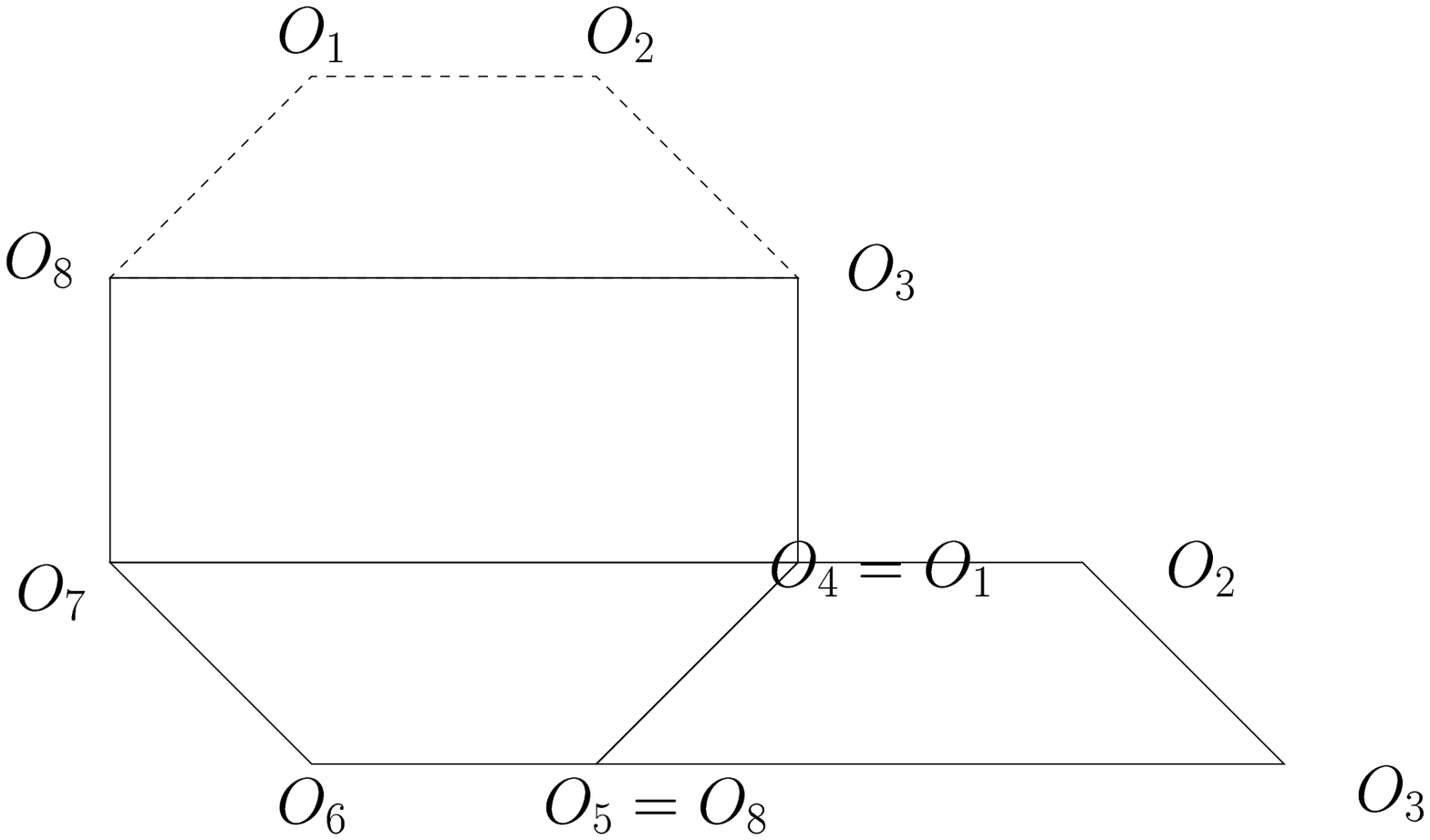}}
\hspace{5mm}
\subfigure[Decomposition in direction $\pi/8$]{\label{cylinders2fig}\includegraphics[width=0.42\textwidth]{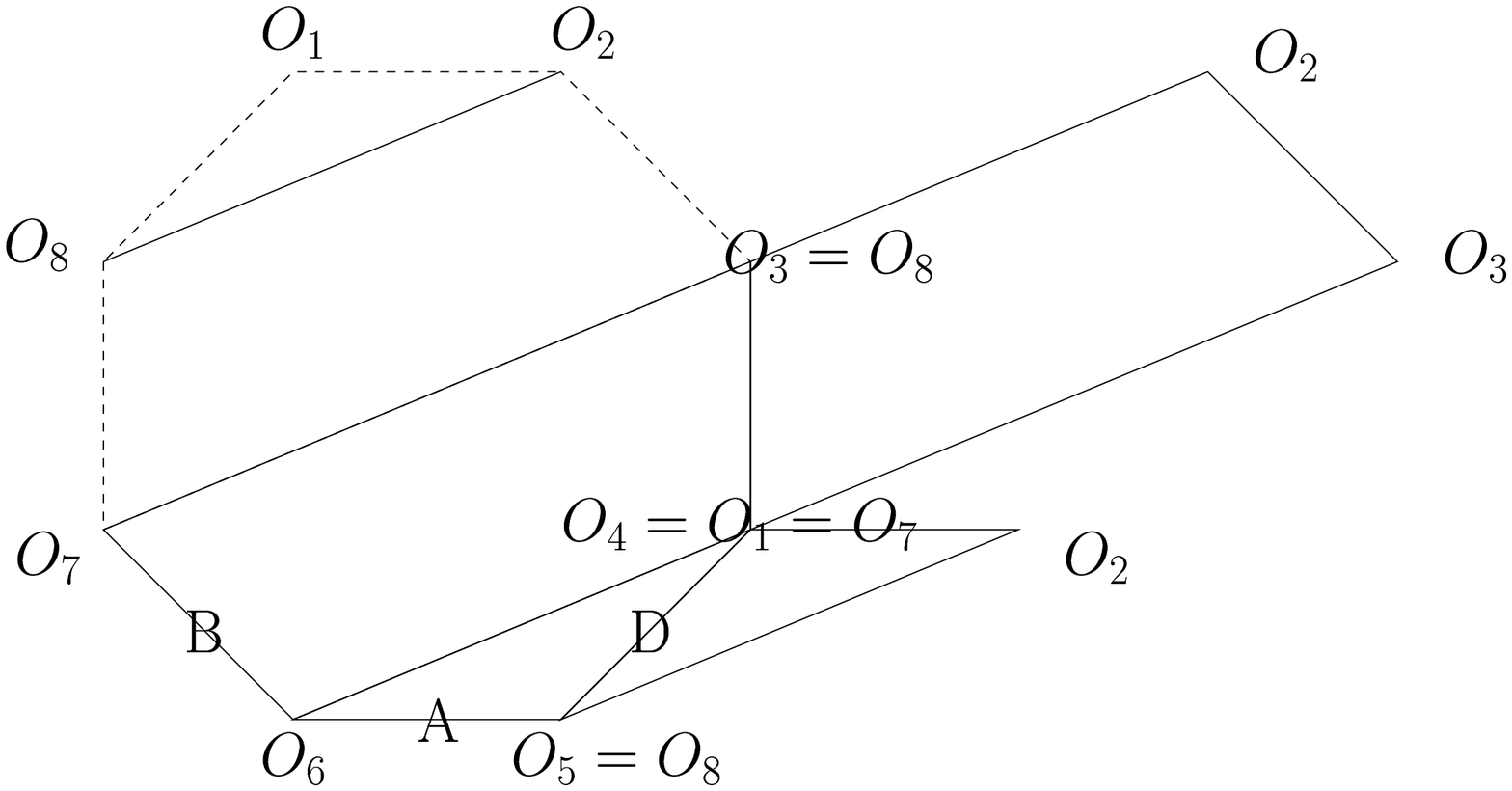}}
\caption{Cylinder decompositions of $S_O$.\label{cylindersdecomp}}
\end{figure}
The surface $S_O$ has a horizontal cylinder decomposition which we can see as follows.
 Cut the octagon along the horizontal diagonals  $\overline{O_8 O_3}$ and $\overline{O_7 O_4}$;  since, recalling the identifications,  $\overline{O_8 O_1}$ is glued to  $\overline{O_5 O_4}$ and  $\overline{O_8 O_7}$ is glued to  $\overline{O_3 O_4}$, we get two parallelograms, which give a cylinder decomposition for the horizontal flow, as in Figure \ref{cylinders1fig}. 
 As remarked by Veech,
the inverse modulus (i.e.~ratio of width divided by height) of one cylinder is $\mu=(1+ \sqrt{2}) $ and the inverse modulus of the other $2\mu$. Hence, the affine diffeomorphism $\Psi_{\sigma}$ which is given in one cylinder by the a Dehn twist which acts by  fixing horizontal lines,  sending  vertical lines  to lines of slope $\mu$ and fixing the singular points and by the square of the same Dehn twist  in the other cylinder, is globally well defined. The corresponding matrix in the Veech group is
\bes 
\sigma  := \begin{pmatrix} 1 &  2(1+ \sqrt{2}) \\ 0 &  1  \end{pmatrix}. \ees 
\noindent The image $O':= \sigma \octcdot O \subset {\mathbb{R}^2} $ is the affine octagon shown in Figure \ref{shearedoctagon1}; $O'$ can be mapped to the original octagon $O$ by cutting it into polygonal pieces, as shown in Figure \ref{shearedoctagon1}, and rearranging these pieces without rotating them to form $O$ (in Figure \ref{shearedoctagon1} the pieces in $O$ and $\sigma \octcdot O $ have been numbered to show this correspondence). Let us denote  by $\Upsilon_o : O' \rightarrow O$ this \emph{cut and paste} map. 
Clearly $D \Upsilon_o = Id$. The automorphism  $\Psi_{\sigma}$ is given by the composition $\Upsilon_o  \sigma$. 
\begin{figure}[!h]
\centering
\includegraphics[width=.9\textwidth]{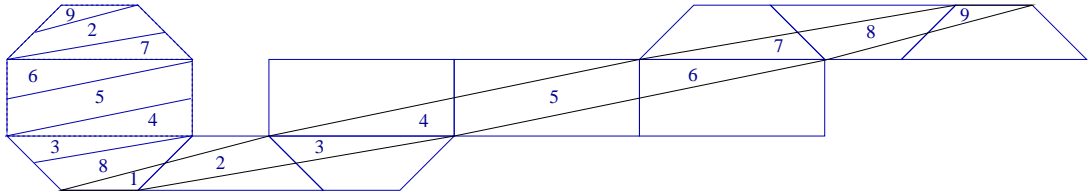}
\caption{The affine octagon $O'$ and the cut and paste map $\Upsilon_o$.\label{shearedoctagon1}}
\end{figure}

\noindent Similarly, also the direction $\pi/8$ has a  decomposition into cylinders, which  is shown in Figure \ref{cylinders2fig}. 

Earle and Gardner (following Veech) proved the following. 
 \begin{prop}\label{Veechthm} The Veech group $V^+_P(S_O)$ is a lattice and it is generated by the rotation $\rho_{\frac{\pi}{4}}$ and by the shear $\sigma$.
\end{prop}

\paragraph{The affine reflection $\Psi_\gamma$.} \label{gammadescription} We will describe an orientation reversing element which will play a crucial role in \S \ref{derivationsec}. 
Let us consider the following matrix:
\be \label{gammadef} 
\gamma :=  \begin{pmatrix} -1 &  2(1+\sqrt{2})  \\ 0  & 1
 \end{pmatrix} . 
\ee
Let us remark that $\gamma$ is an \emph{involution}, i.e. $\gamma^2=id$ or $\gamma=\gamma^{-1}$.
One can  check that  $ \gamma \, \nu_7 = \sigma$, where $\nu_7$ is the reflection at the vertical axes, see (\ref{nujdef}). Thus,
 the action of $\gamma$ on $O$ is obtained by first reflecting it with respect to the vertical axis (this sends $O$ to $O$, but reflecting the orientation), then shearing it through $\sigma$.  The image  $O':= \gamma \octcdot O$ is the same as in Figure \ref{shearedoctagon1}, but the orientation is opposite (the action of $\gamma$ sends the sides of $O$ labelled by $B,C, D$ to the sides of $O'$ labelled by $B',C',D'$ in Figure \ref{cutseqshearedoctagonfig}). Since $\gamma$ is an involution, $\gamma \octcdot O'=O$. 
 If we compose  $\gamma : O \rightarrow O'$ with the cut and paste map  $\Upsilon_{o}: O' \rightarrow O$ defined above (which cut and pastes the pieces in Figure \ref{shearedoctagon1}), 
 we get an affine automorphism $\Psi_{\gamma}:= \Upsilon_{o} \gamma  $ of $S_O$.  
 

Let us give an alternative description of $\Psi_{\gamma}$  that shows that it is a \emph{hidden symmetry} in the following sense. 
The finite order elements of the Veech group arising from isometries of the octagon are isometries of $S_O$. 
Veech describes finite order elements of the Veech group which are not conjugate to isometries in the Veech group as hidden symmetries. Such elements are isometries with respect to some metric, hence their derivatives are conjugate to elements of $O(2)$ in $GL(2,\R)$, but they are not conjugate by elements of $GL(2,\R)$ which are derivatives of elements in  $V(S_O)$.  In the group $V^+_P(S_O)$ there are no such hidden symmetries, but $\Psi_{\gamma}$ is such  an element in $V_P(S_O)$ and  one can see it as follows.   

\begin{figure}[!h]
\centering
\subfigure[The affine octagon $\nu \octcdot O$.]{\label{halfshearedoctagon}
\includegraphics[width=0.5\textwidth]{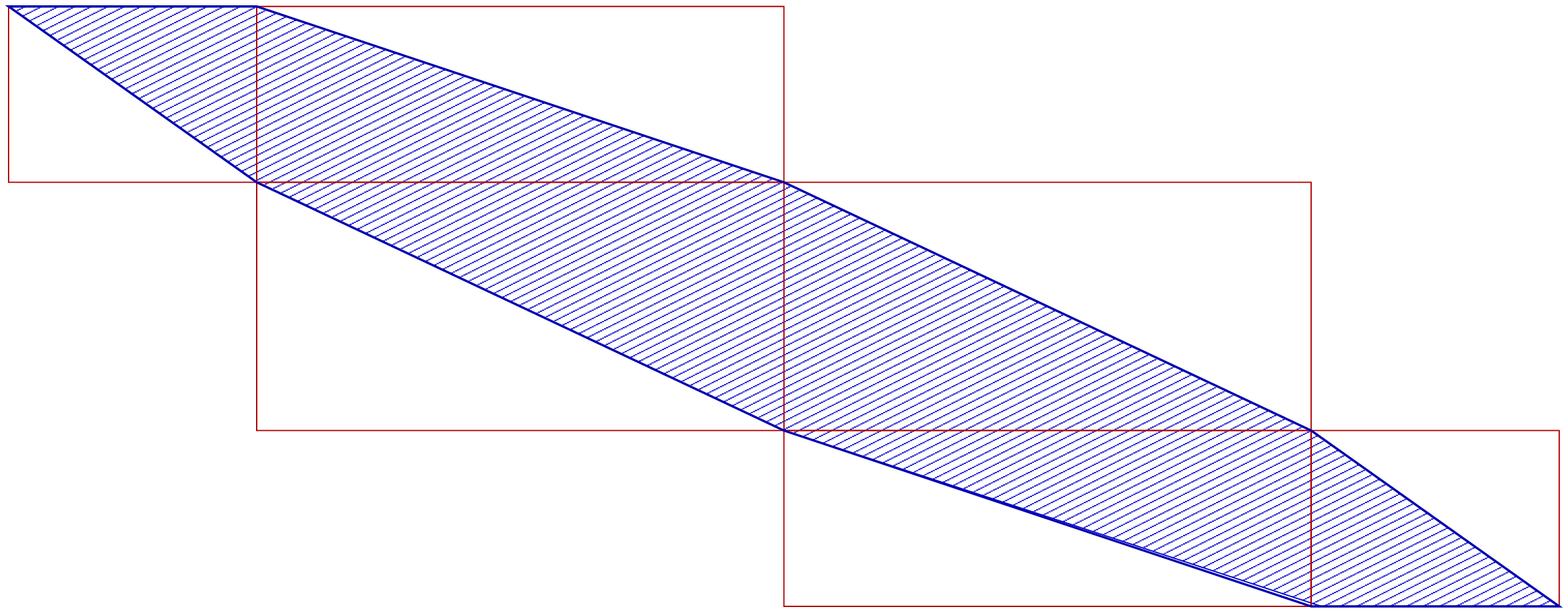}}
\hspace{1mm}
\subfigure[L-shape and symmetry $f_v$]{\label{reflectedLshape}
\includegraphics[width=0.33\textwidth]{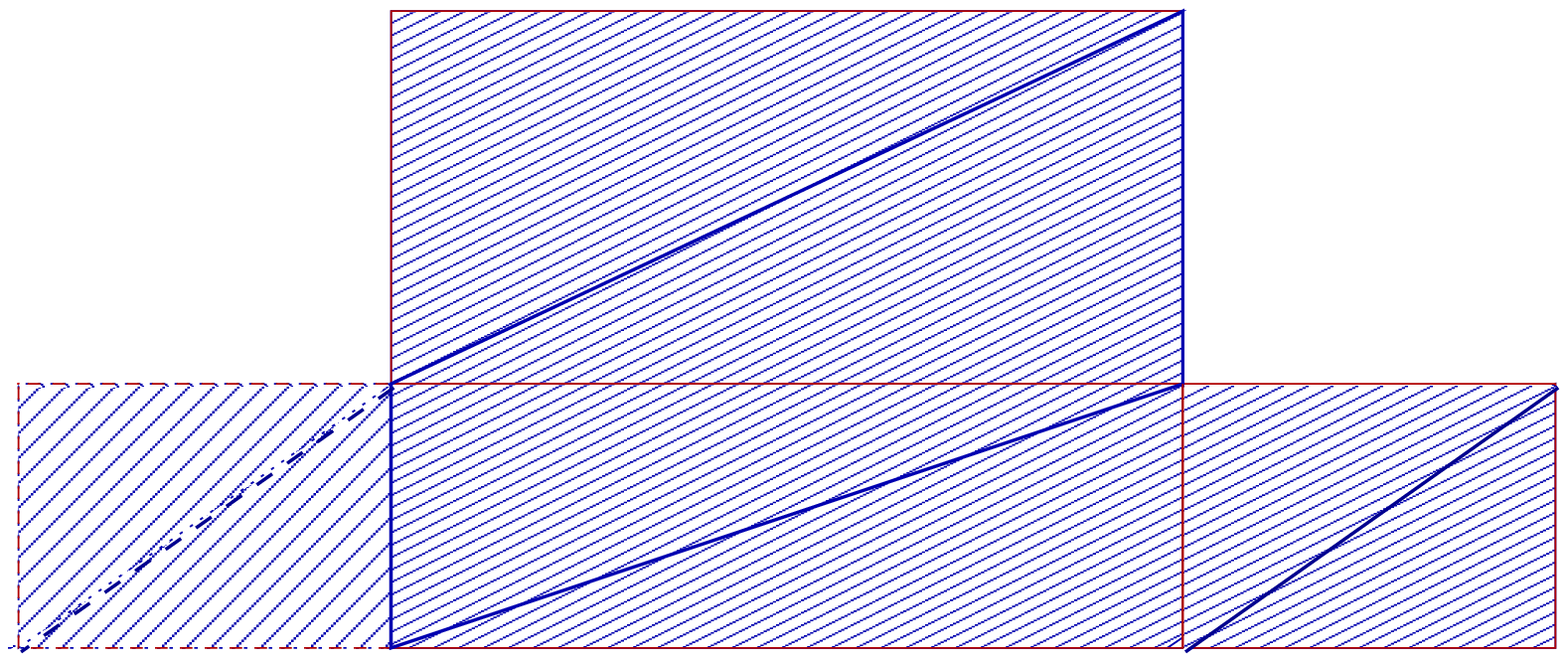}}
\caption{The hidden symmetry $\Psi_\gamma$.}
\end{figure}
We have seen that there both $\theta=0$ and $\theta= \pi/8$ are directions which correspond to cylinder decompositions in $S_O$ (Figure \ref{cylindersdecomp}). We wish to alter the metric (preserving the affine structure) to make these two directions perpendicular. We can do this by applying the linear transformation
\be\label{halfshearnudef}
\nu =\begin{pmatrix} 1 &  -1-\sqrt{2}  \\ 0  & 1
 \end{pmatrix}.
\ee
This transformation preserves the horizontal direction and has the effect of straightening the direction of the second cylinder decomposition so that it is vertical, as shown in Figure \ref{halfshearedoctagon}. These horizontal and vertical cylinders partition $\nu\cdot S_O$ into rectangles (shown in Figure \ref{halfshearedoctagon}). One can cut and paste $\nu\octcdot O $ to get three rectangles in an L shape. Now consider the effect of the linear transformation that reflects in the vertical line. If we apply this to the three rectangles we get three rectangles of the same shape which we can translate back into their original positions (see figure \ref{reflectedLshape}). In particular this vertical reflection $f_v$, composed with a cut and paste map, gives us an isometry (in the sense explained in \S \ref{defnsec}) of the surface $\nu\cdot S_O$. Since the affine automorphism group of $\nu\cdot S_O$ is conjugate to the affine automorphism group of $S_O$ by the action of $\nu$, this isometry corresponds to a hidden symmetry, $\Psi_\gamma$, of the original surface $S_O$. The derivative of the hidden symmetry is the affine reflection $\gamma=\nu^{-1}f_v\nu$ (this relation is well illustrated by Figure \ref{rectanglesedecompfig}).

To describe $\Aff(S_O)$ and  $V(S_O )$ we need the following Lemma.
\begin{lemma}\label{trivialkernel} The  kernel of the Veech homomorphism from $\Aff(S_O)$ to $GL(2,\R)$ is trivial. 
\end{lemma}
\begin{proof}
Let $f$ be an affine automorphism of $S_O$ with $Df=Id$. The map $f$ acts isometrically on $S_O$. The shortest saddle connections in the surface $S_O$ are those corresponding to the boundary of the octagon $O$. Since $f$ is an isometry it takes a saddle connection of minimal length to a saddle connection of minimal length so it preserves the boundary of $O$ and hence it must take $O$ to itself. We have seen however that the isometry group of the octagon is $\Die_8$ so $f\in\Die_8$. On the other hand the only element of $\Die_8$ with derivative equal to the identity is the identity so $F$ is the identity.
\end{proof}
\begin{cor}\label{Psigammainverse}
The affine diffeomorphism  $\Psi_\gamma$ is an involution.
\end{cor}
\begin{proof}
Since $\Psi_{\gamma}=\gamma \Upsilon_0$ and $\gamma$ is an involution, $\Psi_{\gamma}^{-1}=\Upsilon_o^{-1}\gamma$. Thus, since $D\Upsilon_o=id$, we have  $D(\Psi_{\gamma} \Psi_{\gamma}^{-1})=\gamma\, D\Upsilon_o\,  D\Upsilon_o^{-1}\, \gamma = id$. By Lemma \ref{trivialkernel}, this implies that  $\Psi_{\gamma} \Psi_{\gamma}^{-1}=Id$.
\end{proof}

\begin{lemma}[Veech group of the octagon] The Veech group $V(S_O )$ is generated by the linear parts  $\alpha$, $\beta$ and $\gamma$ of the  affine diffeomorphisms $\Psi_\alpha$, $\Psi_\beta$ and $\Psi_\gamma$. The affine automorphism group $\Aff(S_O )$ is generated by $\Psi_\alpha$, $\Psi_\beta$ and $\Psi_\gamma$. 
\end{lemma}
\begin{proof}
Let $\Gamma$ be the subgroup of $GL(2,\R)$ generated by $\alpha$, $\beta$ and $\gamma$. 
Let $\Gamma_P$ be the image of $\Gamma$ in $PGL(2,\R)$ and $\Gamma_P^+$ be the intersection of $\Gamma_P$ with $PSL(2,\R)$. Since $\Gamma\subset V$ we also have $\Gamma_P\subset V_P$ and $\Gamma_P^+\subset V_P^+$. According to \cite{CG:tei} $\Gamma_P^+=V_P^+$. Using this we have the following relation for indices of subgroups: $(\Gamma^+_P:\Gamma_P)(\Gamma_P:V_P)=(V_P^+:V_P)$. Since the subgroup of orientation preserving elements is the kernel the determinant homomorphism $\nu \to \det ( \nu ) $ into the two element group $\{\pm 1\}$ the first and last terms can only be equal to one or two. Since $\Gamma_P$ contains orientation reversing elements we see that 
$(\Gamma^+_P:\Gamma_P)=2$. We conclude that $(V_P^+:V_P)=2$ and $(\Gamma_P:V_P)=1$. So $\Gamma_P=V_P$. 

Now the homomorphisms $\Gamma\to \Gamma_P$ and $V\to V_P$ both have kernels which have order one or two. Since $-Id\in \Gamma$ both kernels have order two. So we have an inclusion of each term of the short exact sequence $\{\pm Id\}\to \Gamma\to \Gamma_P$  into the short exact sequence $\{\pm Id\}\to V\to V_P$. On the first and third terms these inclusions are isomorphisms. It follows that the inclusion of the second term is also an isomorphism so $\Gamma=V$.

Since the kernel of the homomorphism from $\Aff(S_O)$ to $GL(2,\R)$ is trivial by Lemma \ref{trivialkernel},  $\Aff(S_O )$ is isomorphic to its image $V(S_O)$. This shows that  $\Aff(S_O)$ is generated by $\Psi_\alpha$, $\Psi_\beta$ and $\Psi_\gamma$ and concludes the proof.
\end{proof}

\section{Derivation and renormalization schemes}\label{derivationproofssec}
\subsection{Derived sequence as cutting sequence}\label{derivationsec}
In this section we prove Proposition \ref{derivationiscuttingseqprop}. We prove first the following special case. 
\begin{prop}[Renormalization in $\Sigma_0$]\label{derivationiscuttingseqpropsec0}
 Let  $\tau$ be a trajectory in direction $\theta \in \overline{\Sigma}_0$ and let $c(\tau)$ be the cutting sequence of $\tau$. Let  $\tau'= \Psi_{\gamma} \trajcdot \tau$, where 
  $\Psi_{\gamma}$ is described in \S\ref{veechsec}. Then $c(\tau')=c(\tau)'$. 
\end{prop}
\noindent Proposition \ref{derivationiscuttingseqprop} can be proved by reducing it to Proposition \ref{derivationiscuttingseqpropsec0} as follows.
\begin{proofof}{Proposition}{derivationiscuttingseqprop}
Let $w= c(\tau)$ be the cutting sequence of a trajectory $\tau$ in direction $\theta\in \Sigma_k$. Let us show that $w'$ is the cutting sequence of the trajectory $\nu_k^{-1}\trajcdot \Psi_{\gamma} \nu_k \trajcdot \tau$. Since we can apply Proposition \ref{derivationiscuttingseqpropsec0} to the normal form $n(w)=\nu_k\trajcdot \tau$, using also (\ref{relabellingeq}), we have that $ c(\Psi_{\gamma} \nu_k \trajcdot \tau) = c(\nu_k \trajcdot \tau)'= (\pi_k \cdot c(\tau))'$. Since by Remark \ref{commute} derivation commutes with relabelling, $(\pi_k \cdot c(\tau))' = \pi_k \cdot c(\tau)' $. Using again (\ref{relabellingeq}) and combining the previous equalities  we get 
\bes
c(\nu_k^{-1}\trajcdot \Psi_{\gamma} \nu_k \trajcdot \tau ) = \pi_k^{-1} \cdot c( \Psi_{\gamma} \nu_k \trajcdot \tau )=\pi_k^{-1} \cdot  (\pi_k \cdot c(\tau))' =  \pi_k^{-1} \cdot  \pi_k \cdot c(\tau)' = w', 
\ees
which is what we wanted to prove.
\end{proofof}

 Derivation can be seen as the combinatorial counterpart of a geometric operation, which sends a linear trajectory in the octagon into a new linear trajectory in an octagon (see Proposition  \ref{derivationiscuttingseqpropsec0} and the proof of Proposition \ref{derivationiscuttingseqprop}). We note that under this operation the parametrization of the linear trajectory  changes. This is the geometric counterpart to the fact that derivation acts by deleting symbols.

\begin{figure}
\centering
\subfigure[Auxiliary diagonals]{\label{diagonalsoctagonfig}
\includegraphics[width=0.33\textwidth]{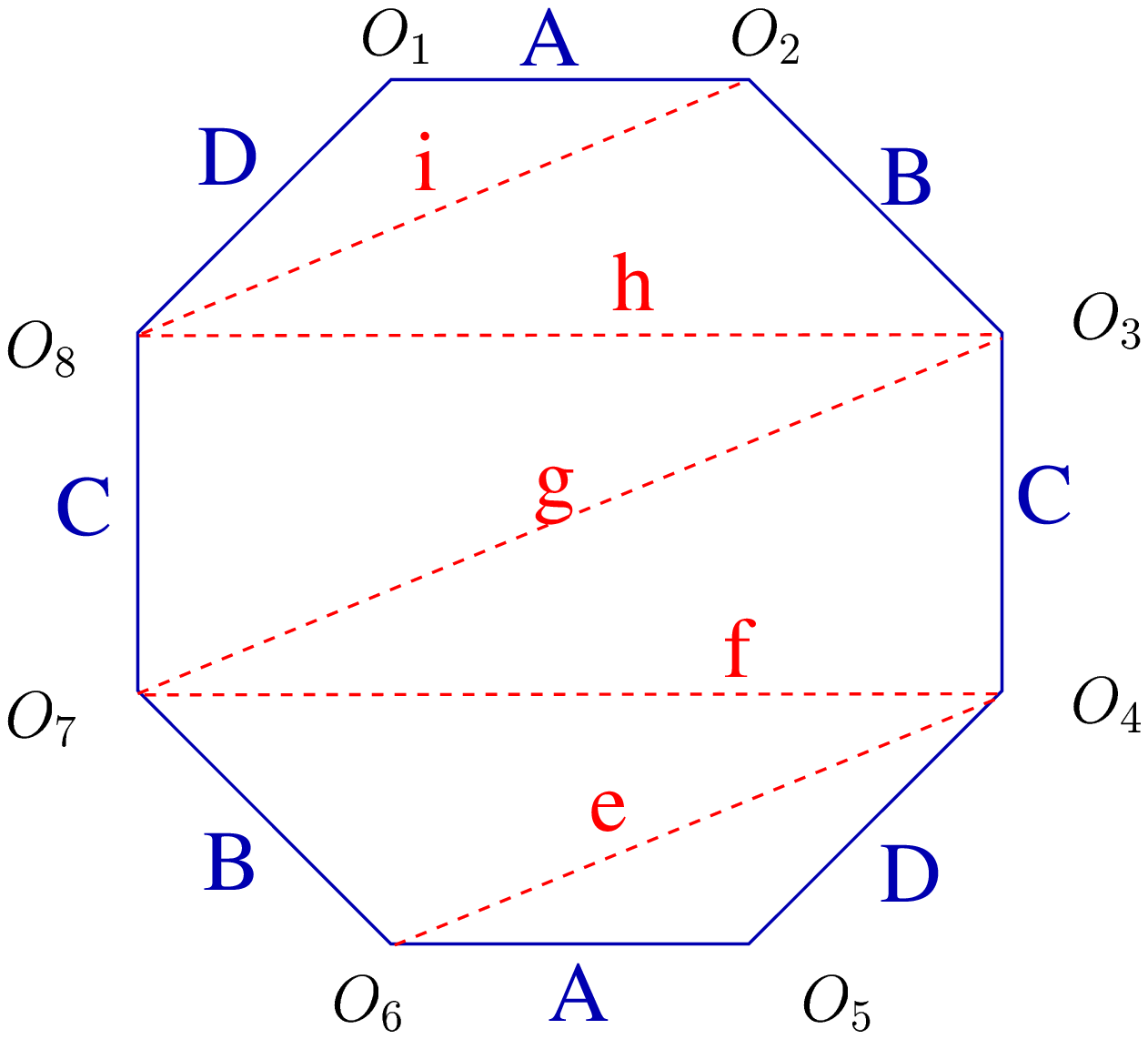}}
\hspace{5mm}
\subfigure[Augmented diagram $\widetilde{\mathscr{D}}_0$]{\label{augmentedD0fig}
\includegraphics[width=0.36\textwidth]{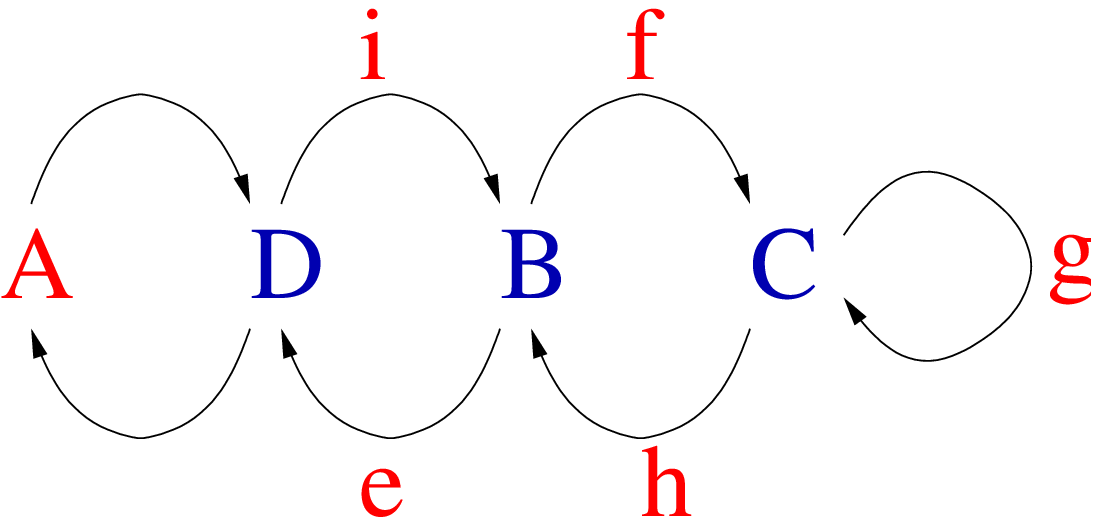}}
\caption{Construction of the augmented cutting sequence for $0\leq \theta \leq \pi/8$.}
\end{figure}

\begin{proofof}{Proposition}{derivationiscuttingseqpropsec0}
Let us add  the octagon diagonals  in directions $\theta= 0 $ and $\theta= \pi/8$, as in Figure \ref{diagonalsoctagonfig}, i.e.\ the diagonals $\overline{O_6 O _4}$, $\overline{O_7 O _4}$, $\overline{O_7 O _3}$, $\overline{O_8 O _3}$, $\overline{O_8 O _2}$, as auxiliary edges and let us label them by the  letters $e$, $f$, $g$, $h$, $i$ respectively.   
Let $\tilde{c}(\tau) \in \{A,B,C,D,e,f,g,h,i\}^{\mathbb{Z}}$ be the cutting sequence of $\tau$ with respect to the octagon sides \emph{and} the auxiliary edges, which we will call the  \emph{augmented sequence}\footnote{Compare with the proof of Proposition \ref{toruscutseq} for the square, where the auxiliary diagonal in direction $\pi/4$, labelled by $c$, was added.}. 

The augmented sequence $\tilde{c}(\tau)$ is completely determined by the cutting sequence $c(\tau)$. For example, since $0\leq \theta \leq  \pi/8$, whenever there is a $BD$ transition, the trajectory crosses from side $\overline{O_6 O_7}$ to side $\overline{O_4 O_5}$ (see Figure \ref{diagonalsoctagonfig}).  In between the trajectory is forced to cross $e$. Similarly when there is a $BC$ transition, $\tau$ is forced to cross $f$. Reasoning in this way, one can obtain the augmented diagram $\widetilde{\mathscr{D}}_0$ in Figure \ref{augmentedD0fig}, where each edge, which corresponds to a transition, has been labeled by the letter of the auxiliary side crossed by the trajectory during that transition. We remark that the two edges with no labels correspond to transitions ($DA$ and $AD$) during which no auxiliary diagonal is crossed. 

Recall that ${c}(\tau)$ corresponds to an infinite path on   $\widetilde{\mathscr{D}}_0$ and  $\widetilde{c}(\tau)$ is simply the sequence of labels of both vertices and edges along this path; from the combinatorial structure of $\widetilde{\mathscr{D}}_0$ in Figure \ref{augmentedD0fig}, one can see that the subsequence of the letters $e,f,g,h,i$  together with the letter\footnote{The letter $A$ plays a special role since it corresponds to an octagon side in direction $\theta=0$ in the sector $\Sigma_0$.}  $A$ completely determine the path.
 Let us therefore switch the role of edges and vertexes: consider the graph  $\widetilde{\mathscr{D}}_0'$ shown in Figure \ref{dualDfig} which is obtained by switching the role of all edges and vertices, with the exception of the vertex $A$ (and hence ``almost''  dual to $\widetilde{\mathscr{D}}_0$). Hence, let us keep only the letters $\{A,e,f,g,h,i\}$ in $\tilde{c}(\tau)$ (which we call \emph{auxiliary letters}) and get an auxiliary word $\hat{c}(\tau) \in \{A,e,f,g,h,i\}^{\mathbb{Z}} $ which corresponds to  the sequence of vertices of an infinite path on $\widetilde{\mathscr{D}}_0'$. 

Consider now the affine octagon $O' := \gamma \octcdot O$ in Figure \ref{cutseqshearedoctagonfig}, obtained by applying the affine reflection $\gamma$  in (\ref{gammadef}) (see \S \ref{gammadescription}). 
\begin{figure}
\centering
\includegraphics[width=1\textwidth]{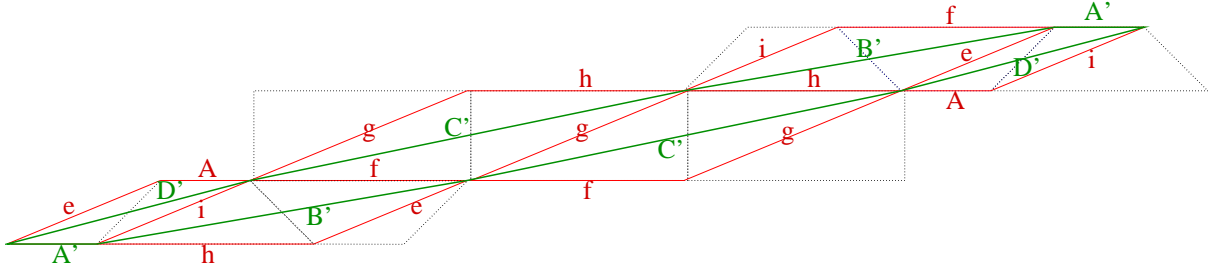}
\caption{The affine octagon $O':= \gamma \octcdot O$ and the auxiliary diagonals.\label{cutseqshearedoctagonfig}}
\end{figure}
Let us recall that $O' = \Upsilon_o^{-1} O$, where $ \Upsilon_o : O' \rightarrow O$ is the cut and paste map defined in \S \ref{veechsec} referring to  Figure \ref{shearedoctagon1}.  The image $\Upsilon_o^{-1} \tau$ of the trajectory $\tau$ is clearly a trajectory on $O'$ in the same direction $\theta$. Let $\overline{c}(\tau) \in \{ A, B', C', D'\}^{\mathbb{Z}}$ be the cutting sequence of $\Upsilon_0^{-1} \tau$ with respect to the sides of $O'$, where we
denote by primed letters $ B', C', D'$ the sides of $O'$ which are the image under $\gamma$ of the sides of $O$ labelled by $B$, $C$, $D$ respectively (this labelling of the sides of $O'$ is shown in Figure \ref{cutseqshearedoctagonfig}). 
\begin{figure}
\centering
\subfigure[Cutting sequences for $S_0$ on $\widetilde{\mathscr{D}}_0'$]{\label{dualDfig}
\includegraphics[width=0.42\textwidth]{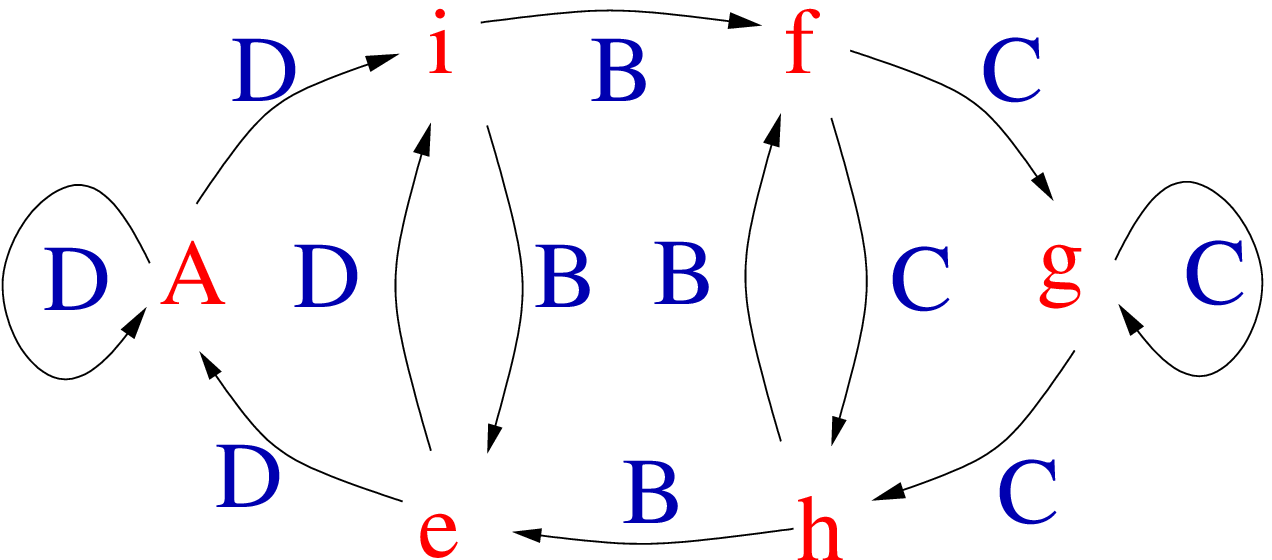}}
\hspace{2mm}
\subfigure[Cutting sequences for $\gamma \cdot S_0$ on ${\mathscr{D}}_0'$]{\label{dualDnewoctagonfig}
\includegraphics[width=0.42\textwidth]{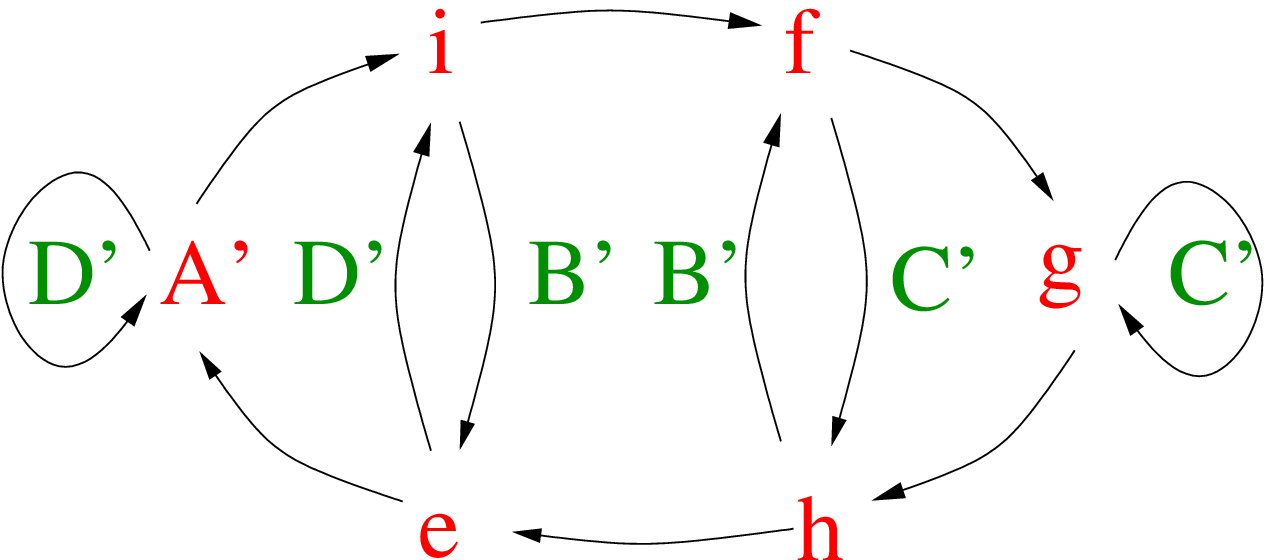}}
\caption{Comparison of dual augmented transition graphs for $S_O$ and $\gamma \cdot S_O$.}
\end{figure}   
Let us show that  also the cutting  sequence $\overline{c}(\tau)$ is completely determined by the auxiliary sequence $\hat{c}(\tau)$ of letters $\{A,e,f,g,h,i\} $. 
As it can be seen in Figure \ref{cutseqshearedoctagonfig}, each of the sides of $O'$ is a diagonal of a parallelogram (with sides in direction $0$ or $\pi/8$) formed by auxiliary edges together with  the side labelled by $A$. Thus, using the assumption that $0\leq \theta \leq \pi/8$ and reasoning as before, the transitions in $\hat{c}(\tau)$ (i.e. pairs of consecutive letters $\{A,e,f,g,h,i\}$)  determine whether the corresponding diagonals are crossed. For example, if the trajectory hits the auxiliary side $e$ and then $i$, it has to cross $D'$ (see Figure \ref{cutseqshearedoctagonfig}) but during an  $f g$ transition, no sides of $O'$ are crossed.  The graph ${\mathscr{D}}_0'$ in Figure \ref{dualDnewoctagonfig} is obtained by relabelling   the edges of  $\widetilde{\mathscr{D}}_0'$, labelling  each edge by the letter in $\{ A, B', C', D'\}$ of the side of $O'$ crossed during the auxiliary letter transition corresponding to that edge. If no edge of $O'$ is crossed, the edge does not carry any label (see Figure \ref{dualDnewoctagonfig}).  
A convenient picture for visualizing the previous assertions can be obtained by applying the shear $\nu$ in (\ref{halfshearnudef}) to $O$ and $O'$. This has the effect of mapping all parallelograms made by auxiliary sides to rectangles, as shown  in Figure \ref{rectanglesedecompfig}. Figure \ref{rectanglesedecompfig} also shows that $O$ is sheared to the left  by $\nu$ and  $ O'$ is sheared back, so that  $\nu \octcdot O$ and $\nu \octcdot O'$ become exactly symmetric to each other with respect to the vertical axes.  Both $\nu\octcdot O$ and $\nu\octcdot O'$  can be cut and pasted into a union of three rectangles\footnote{
We remark that the pairs of diagonals of these rectangles are $\{B,B'\}$, $\{C,C'\}$ and $\{D,D'\}$; the affine octagon $O'$ is indeed obtained from $O$ by three simultaneous \emph{Delaunay switches} which consists in changing a rectangle diagonal with the opposite one.} in an L-shape. Furthermore, in Figure \ref{rectanglesedecompfig}, the trajectory  $\tau$ is sent to the trajectory $\nu \trajcdot \tau $ whose direction $\theta_{\nu}$ satisfies $0 \leq \theta_\nu \leq \pi/2$.  Thus, all the edges labels in Figure \ref{dualDnewoctagonfig}  can be easily checked by looking at rectangles and reasoning as in the case of the square (see \S \ref{torussec}).  

\begin{figure}
\centering
\includegraphics[width=1\textwidth]{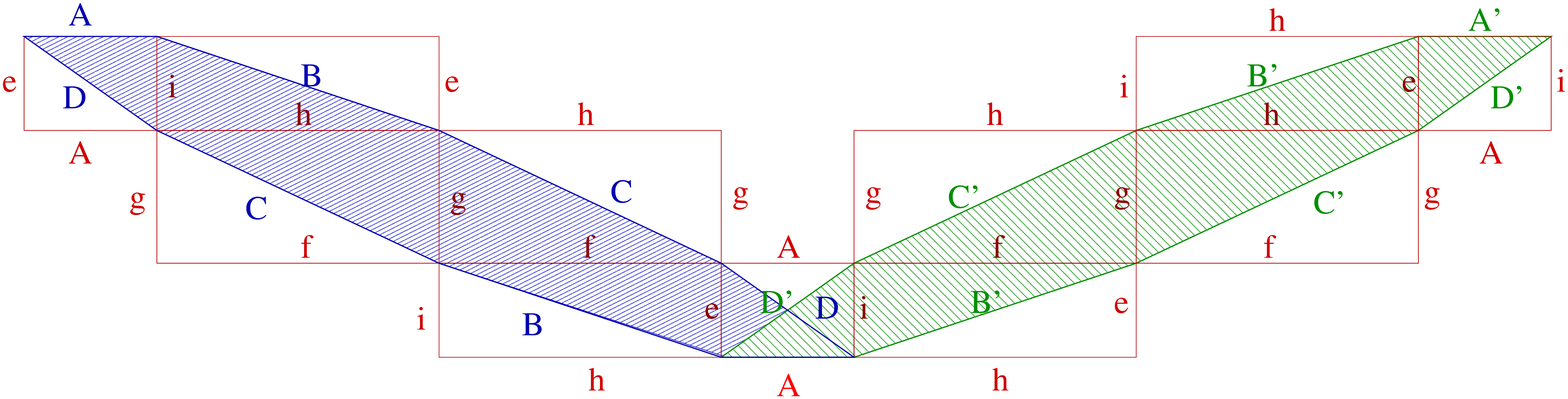}
\caption{Sheared copies of $O$ and $O'$ and rectangles decompositions.\label{rectanglesedecompfig}}
\end{figure}   
 
Let us now show that the cutting sequence $\overline{c}(\tau)$ is exactly the derived sequence of $c(\tau)$ with $B$, $C$ and $D$ replaced by $B'$, $C'$ and $D'$. As we have just proved, both cutting sequences  are determined by the sequence $\hat{c}(\tau)$ of auxiliary edges crossed by $\tau$: $c(\tau)$ is obtained reading off the labels of the edges on $\widetilde{\mathscr{D}}_0'$ of the infinite path  given by  $\hat{c}(\tau)$, while $\overline{c}(\tau)$ is obtained reading off the labels of the edges on  $\mathscr{{D}}_0'$ corresponding to the same path. Hence it is enough to compare the two different labellings of edges in  $\widetilde{\mathscr{D}}_0'$ and  $\mathscr{{D}}_0'$ (i.e. Figures \ref{dualDfig} and \ref{dualDnewoctagonfig}). Observe that the labelled edges in $\mathscr{{D}}_0'$ are a subset of those in $\widetilde{\mathscr{D}}_0'$ and are labelled by the corresponding primed letters. One can check that each of these labelled edges corresponds to sandwiched letters in $\widetilde{\mathscr{D}}_0'$. Some oriented edges of $\widetilde{\mathscr{D}}_0'$ are preceeded or followed by more than one edge, but when an oriented edge of $\widetilde{\mathscr{D}}_0'$ is preceded  by two edges, 
they are both labelled by the same letter and the same is true when it is followed by two edges. When an oriented edge is labelled by a sandwiched letter, all its incoming and outgoing oriented edges carry the same label. For example, a path that uses the $D$ edge from $e$ to $i$ in $\widetilde{\mathscr{D}}_0'$, came from an edge labelled $B$ and can be followed by two edges, both labelled by $B$, hence the corresponding $D$ is $B$-sandwiched. Conversely, one can also check that each of the edges which do not carry a label in ${\mathscr{D}}_0'$ correspond to letters which are not sandwiched. For example, if a path goes through the $B$ edge from $i$ to $f$ in $\widetilde{\mathscr{D}}_0'$, there are two possible  previous edges, both labelled by $D$   and two possible following ones,  both are labelled by $C$, hence the $B$ is not sandwiched. 
Thus, we have proved that the cutting sequence $\overline{c}(\tau)$  of the trajectory $\Upsilon_o^{-1}  \tau$ (which is a trajectory in direction $\theta$ in $O'$) with respect to the labels of the affine octagon $O'$ is obtained from $w$ by erasing exactly the letters which are not sandwiched and writing $B'$, $C'$, $D'$ for $B,C,D$ respectively. Hence, $\overline{c}(\tau)$   is exactly the derived sequence $w'$, with $B,C,D$ primed.  

In order to conclude the proof of the Proposition, it is enough to 
  apply  $\gamma$, so that $O'$ is sent to $O$, the sides labelled by $B', C', D'$ are sent to  $B, C, D$ respectively (by definition of $B',C',D'$) and $\Upsilon_o^{-1}\tau$ is sent to the trajectory $\tau'= \gamma \Upsilon_o^{-1 } \tau$. Thus, $w'$ is the cutting sequence $c(\tau')$ of a linear trajectory $\tau'$ in $O$. Finally, since $\Psi_{\gamma} = \Upsilon_o \gamma = \gamma \Upsilon_o^{-1}$ by Corollary \ref{Psigammainverse}, we have $\tau' = \Psi_{\gamma}\tau$. 
\end{proofof}

\subsection{Direction of renormalized trajectories and itineraries}
In this section we give the proof  by induction of Proposition \ref{typesequalCF}. Let us first prove two Lemmas which we use in  the  inductive step.
\begin{lemma}\label{sec0directionrenormalization}
If  the direction $\theta $ of the trajectory $\tau$ is in $\Sigma_0$, then the direction $\theta'$ of the trajectory  $\tau'= \Psi_\gamma \tau$  is given by 
\be\label{theta'}
\theta' = F_0(\theta)= F(\theta),  
\ee
where $F$ is the octagon Farey map in \S \ref{Fareymapsec}. 
\end{lemma}
\begin{proof}
 Since $\Psi_{\gamma}  = \gamma \Upsilon_{O}^{-1} $ (by definition of $\Psi_\gamma$ and Remark \ref{Psigammainverse}),  $\tau'$ is obtained from $\Upsilon_O^{-1} \tau$ acting by $\gamma$. The cut and pasted trajectory $\Upsilon_O^{-1} \trajcdot \tau$ has the same direction $\theta$ of $\tau$. Thus, considering 
 the action of  $\gamma$ induced on inverses of slopes, we get that the direction $ \theta'$ of $\tau'$ is given by $ \theta' =\gamma[\theta]$, where the action is the one defined in (\ref{actionslopes}). This concludes the proof, since $  \gamma[\theta] = F_0(\theta) $ by definition (\ref{Fareymapdef}) of the Farey map $F$.  
\end{proof}

\begin{lemma}[Inductive step] \label{inductivestep} Let $\tau$ be a trajectory in direction $\theta \in \Sigma_{i}$, $0\leq i \leq 7$, and $n(\tau)$ its normal form (given by Definition \ref{normalizedtrajdef}). Then the direction $\theta'$  of the derived trajectory  $n(\tau)'$ (see Definition \ref{derivedtrajdef}) is given by
\bes
\theta' = F_{i}(\theta).
\ees
\end{lemma}
\begin{proof}
Since by assumption $s(\tau)=i$, the normal form, by Definition \ref{normalizedtrajdef}, is given by $n(\tau)=\nu_i \trajcdot \tau$.  Thus, using the notation (\ref{actionslopes}) introduced in \S \ref{Fareymapsec}, the direction of $n(\tau)$ is $\nu_{i} [\theta] \in \Sigma_0$. Hence, by Lemma \ref{sec0directionrenormalization}, the direction $\theta'$ of  $n(\tau)'$ is  $F_0( \nu_{i} [\theta] )$. By definition (\ref{Fareymapdef}) of $F$,  $F_0 (\nu_{i} [\theta] )= F_{i} (\theta) $. 
\end{proof}

\begin{proofof}{Proposition}{typesequalCF}
Let us show by induction that $\theta_k=F^k(\theta)$. For $k=1$, since $\theta \in \Sigma_{s(\tau)}$ and  $\theta_1$ is by definition the direction of $\tau_1 = n(\tau)'$, by Lemma \ref{inductivestep} we have $\theta_1 = F_{s(\tau)}(\theta) = F(\theta)$. In the same way, for all $k\in \mathbb{N}$, since $\theta_{k} \in \Sigma_{s_{k}(\tau)}$ (recall the Definition \ref{sectorsdef} of the sequence of sectors)  and  $\theta_{k+1}$ is  the direction of $\tau_{k+1} = n(\tau_{k})'$ (recall Definition \ref{renormalizationtrajectoriesdef}), Lemma \ref{inductivestep} gives that $\theta_{k+1} = F_{s_k(\tau)}(\theta_k) = F(\theta_k)$ and, since by inductive assumption $\theta_k = F^k(\theta)$, we get $\theta_{k+1}= F^{k+1}(\theta)$. The second part of Proposition \ref{typesequalCF} is simply a rephrasing of Definition \ref{sectorsdef} of the sequence of sectors and Definition \ref{itinerarydef} of itinerary. 
\end{proofof}

\subsection{Properties of the octagon additive continued fraction expansion }\label{contfracpropproof}
We include here for completeness the proofs of Lemmas \ref{contfracprop1} and \ref{contfracprop2}, which are obtained by applying to our map $F$ the classical arguments about coding of piecewise expanding maps.
 \begin{proofof}{Lemma}{contfracprop1}
Let $s_0 \in \{0,\dots, 7\}$ and $s_k\in \{1,\dots, 7\}$ for all $k\geq 1$. 
Let us prove first that the intersection $\bigcap_k \overline{\Sigma}[s_0,\dots,s_k]$  is \emph{non-empty}. 
Remarking that $F_{s_k}^{-1}[0,\pi]\subset [\pi/8,\pi]$ and for all branches $F_i$, $i=1,\dots, 7$, $F_i^{-1}[\pi/8,\pi]\subset [\pi/8,\pi]$, one can see that all intersections of the form $F_{s_0 }^{-1} F_{s_1 }^{-1} \cdots  F_{s_k}^{-1} [0, \pi]$ (which by definition give $\overline{\Sigma}[s_0,\dots,s_k]$) are non-empty. 
 Moreover, as $k$ increases, they give a sequence of nested compact sets. This implies that the intersection (\ref{directionasintersection}) is non-empty. 

Let us now know show that the intersection $\bigcap_k \overline{\Sigma}[s_0,\dots,s_k]$   
 consists of a \emph{unique point}. 
The inverse branches  $F_i^{-1}$ are contracting but not uniformly contracting. One can verify 
that they are uniformly contracting  outside of a fixed neighborhood  of the points $\{\pi/8, \pi \}$, which are the only points where the derivative is $\pm 1$. If $s_{k}\notin \{0,1,7\}$, $F_{s_k}^{-1}[0,\pi] \subset [2\pi/8,7\pi/8]$. Moreover,  since  $F_7^{-1}([0,\pi])= \overline{\Sigma_7}$ and $F_1^{-1}([0,\pi])=\overline{ \Sigma}_1$, both $F_1^{-1}F_7^{-1} [0,\pi]$ and $F_7^{-1}F_1^{-1} [0,\pi]$ are contained in $[\underline{\theta}, \overline{\theta}] $ where $ \underline{\theta}:= F_1^{-1}(7\pi/8)> \pi/8 $ and $ \overline{\theta}:= F_7^{-1}(2\pi/8) < \pi$, so in all this cases the composition of the inverse branches is uniformly contracting.

Hence, if the sequence $\{s_k\}_{k\in \mathbb{N}} \in S^*$ contains either infinitely many $s_k\notin \{ 1,7 \}$ (recall that $s_k\neq 0$ if $k>0$, so this gives infinitely many $s_k\notin \{ 0,1,7 \}$) or infinitely pairs $\{s_k,s_{k+1}\}=\{1,7\}$, the intersection   (\ref{directionasintersection})  consists of a unique point. 
The only possibility that we are left to consider is that there exists a $k_0$ such that 
 $s_k=1$ for all $k> k_0$ or $s_k=7$ for all $k> k_0$. In this case, the image of the intersection (\ref{directionasintersection}) under $F_{s_{k_0}} \dots F_{s_{1}}F_{s_{0}}$ must consist of  either the fixed point $\pi/8$ or the fixed point $\pi$ (the only directions whose sequence of sectors is constantly $1$ or $7$) and hence again the intersection consists of a unique point (the preimage of a fixed point).
\end{proofof}

\begin{proofof}{Lemma}{contfracprop2}
Let us prove $(i)$. Given $\theta$, let $\{s_k\}_{k\in\mathbb{N}}$ be its itinerary under the Farey map. By definition of itinerary,   $\theta_k:= F^k(\theta ) \in \Sigma_{s_k}$, so that for each $k\in \mathbb{N}$ we can write $\theta_{k}=F^{-1}_{s_{k}}(\theta_{k+1})$ . Thus, $\theta =\theta_0=  F^{-1}_{s_0} F^{-1}_{s_{1}}\ldots F^{-1}_{s_k} \theta_{k+1}$. This shows that $\theta$ belongs to the intersection (\ref{directionasintersection}). Moreover, since the sequence $s_0,s_1,s_2,\dots$ is an itinerary, it satisfies the constraints in  Definition \ref{CFexpansion} by Remark \ref{itineraryproperties}. Thus, $\theta = [s_0; s_1,s_2,\dots]$. 

To prove $(ii)$, assume that $\theta$ has two expansions $[s_0'; s_1', \dots]$ and $[s_0''; s_1'', \dots]$. Let $k$ be the minimum $l$ such that $s_l'\neq s_l''$. Thus, by Definition \ref{CFexpansion}, $\overline{\theta}:=F_{s'_{k}}\ldots F_{s'_0}(\theta) =  F_{s''_{k}}\ldots F_{s''_{0}}(\theta) $ belongs to both $F_{s_{k}'}^{-1}[0,\pi]$ and to $F_{s_{k}''}^{-1}[0,\pi]$. This shows that, up to inverting the role of $s_k'$ and $s_k''$, $s_{k'}=s_k$ and $s_{k}''=s_k+1$ for some $s_k \in \{1,\dots, 6\}$, since otherwise $F_{s_k'}^{-1}[0,\pi] $ and $F_{s_k''}^{-1}[0,\pi] $ are disjoint  (and by Definition \ref{CFexpansion} $s_{k+1}\neq 0$, $s_k\neq 0$).  Thus,  $\overline{\theta} \in F_{s_k}^{-1}[0,\pi] \cap F_{s_k+1}^{-1}[0,\pi] = \pi (s_k+1)/8 $. Thus,  for all $l\in \mathbb{N}$, $F^l (\overline{\theta} )$ is either equal  to $\pi$ (if $s_{k}$ is odd) or $\pi/8$ (if $s_{k}$ is even). Since $\pi$ belongs only to  $F_{7}^{-1}[0,\pi]$ and $\pi/8$ belongs only to  $F_{1}^{-1}[0,\pi]$ or  $F_{0}^{-1}[0,\pi] $ and we exclude  by Definition \ref{CFexpansion} that $s_k'$ or $s_k''$ are $0$, this shows that either $s_l'=s_l' = 7$ for all $l>k$ (if $s_{k}-1$ is even)  or  $s_l'=s_l ''= 7$ for all $l>k$ (if $s_{k}-1$ is odd). 
\end{proofof}

\subsection{Direction recognition}\label{directionsthmproof}
We first prove Proposition \ref{admsectorsfornonper} and then Theorem \ref{directionsthm}. 
\begin{proofof}{Proposition}{admsectorsfornonper}
Assume that $w$ is non-periodic.  Let us first show that $w=w_0$ is admissible in an unique diagram. We know that $w$ is the cutting sequence of some $\tau$ in an unknown  direction $\theta$.
 Let $s_0= s(\tau)$ so that $w$ is admissible in diagram $s_0$. A priori $w$ could be admissible in some other diagram too (see Example \ref{periodicsameseq}) and we want to rule out this possibility. 
 We are going to show that all transitions which are allowed in sector $s_0$ actually occur, so that by Lemma \ref{all_transitions}  $w$ is admissible \emph{only} in diagram $s_0$.

Since $w$ is non-periodic, the trajectory $\tau$ cannot be periodic. The Veech dichotomy (see \S\ref{defnsec}) this implies that $\tau$ is dense in $O$.
Let $L_1L_2$ be a transition allowed in  $\mathscr{D}_{s_0}$. This means that we can choose inside the octagon $O$ a segment in direction $\theta$ that connects an interior point on a side labelled by $L_1$ with an interior point on a side labelled $L_2$. Since $\tau$ is dense, it comes arbitrarily close to the segment. Since by construction $\tau$ and the segment are parallel, this shows that $c(\tau)$ contains the transition $L_1L_2$. Repeating the argument for all transitions in $\mathscr{D}_{s_0}$, we get that $w$ gives a path on $\mathscr{D}_{s_0}$ which goes through all arrows.  This shows by Lemma \ref{all_transitions} and Remark \ref{uniquenesstyperk} that $s_{0}(\tau)$ is uniquely determined by the sequence $w$ and equal to $s_0$.

Let us proceed by induction. Assume that we showed that $s_0, \dots ,s_k$ are uniquely defined and that we have used them to define the words $w_{k+1}$ as in (\ref{renormalizedseq}). Let $\{\tau_k\}_{k\in\mathbb{N}}$ be the sequence of normalized trajectories in Definition \ref{renormalizationtrajectoriesdef}. We have $c(\tau_{k+1}) = w_{k+1}$ by Lemma \ref{cutseqrelationlemma}.  
Since $\tau_{k+1}$ is obtained from a  dense trajectory $\tau$ by applying a sequence of elements of the Veech group and affine diffeomorphisms are in particular homeomorphisms, it will also be the case that $\tau_{k+1}$ is a dense trajectory. Hence, we can repeat the same argument used for $w_0$ to show that $w_{k+1}$ is admissible  only in  diagram $s_{k+1}:=s(\tau_{k+1})$ and hence $s_{k+1}$ is uniquely defined and can be recovered only from the knowledge of $w_k$. 
\end{proofof}

\begin{proofof}{Theorem}{directionsthm}
Let $w=c(\tau)$ be non-periodic. By Proposition \ref{admsectorsfornonper}, the combinatorial renormalization algorithm is well defined for all $k\in \mathbb{N}$ and each $w_k$ is admissible in a unique diagram $d_k(\tau)$. Since by Lemma \ref{cutseqrelationlemma}, $c(\tau_k)=w_k$ and $\tau_k$ is in sector $s_k(\tau)$, we also have that $w_k$ is admissible in diagram $s_k(\tau)$ and this shows that $d_k(\tau)=s_k(\tau)$ for all $k\in \mathbb{N}$. 

Let $\theta$ be the direction of $\tau$. By Proposition \ref{typesequalCF}, the itinerary of $\theta$ under $F$ is hence $d_0(\tau)$, $d_1(w)$, $d_2(w), \dots  $ and this, by Lemma \ref{contfracprop2} $(i)$, shows  that $\theta$ is given by the octagon additive continued fraction expansion $[d_0(w); d_1(w), d_2(w), \dots ]$ and is unique (see  Lemma \ref{contfracprop1}, recalling Definition \ref{CFexpansion}).
\end{proofof}

\subsection{Terminating directions and periodic sequences}\label{3periodicproof}
Let us prove the  characterization of terminating directions in terms  of periodic trajectories.
\begin{proofof}{Proposition}{3periodic}
The implication $(i) \Rightarrow (ii)$ is obvious. Let us show $(ii)\Rightarrow (iii)$. Let $w=c(\tau)$ be a periodic cutting sequence. Let us remark that the operation of derivation strictly decreases the period unless all letters are sandwiched. Thus, in a finite number iterations (say $k$) of the renormalization scheme in \S\ref{trajrenormsec}, we obtain a trajectory $\tau_k$ such that  $w_k = c(\tau_k)$ contains only sandwiched letters. Let  $w_k = \dots L_{j-1} L_j L_{j+1} \dots $. To say that $L_j$ is sandwiched means that $L_{j-1}=L_{j+1}$.  Thus, since all letters are sandwiched, $w_k$ has period one or two. One can see from the diagrams in Figure \ref{diagrams} and in Figure \ref{boundarydiagrams}  that all  admissible sequences of period one or two (thus in particular $w_k$)  are admissible in a  diagram $\mathscr{D}_{j,j+1}$  in Figure \ref{boundarydiagrams}. This forces   the direction $\theta_k$ of $\tau_k$ to be a multiple $j \pi/8$ for some $j\in \mathbb{N}$. Let us recall that by Proposition \ref{typesequalCF} $\theta_l = F^l(\theta) $ for all $l\in\mathbb{N}$. Thus, for all $l >0$, $F^{l+k}(\theta) = F^l(\theta_k) = F^l(\pi j /8)$ is either always $\pi/8$ or always $\pi$ (according to the parity of $j$) and the itinerary of $\theta$ is respectively either eventually $1$ or eventually $7$.  By Lemma \ref{contfracprop2} $(i)$, this shows  that the direction is terminating (Definition \ref{terminatingdef}). 

Let us conclude by showing  $(iii) \Rightarrow (i)$. If $\theta$ is terminating, there exists $k$ such that either $F^k(\theta) = \pi$ or $F^k(\theta) = \pi/8$. Consider the corresponding renormalized trajectory $\tau_k$, which, by Proposition \ref{typesequalCF},   has direction $\theta_k = F^k (\theta)$ equal to $\pi$ or $\pi/8$. Since  trajectories in these two directions are  periodic (see Remark \ref{periodicmultipleofpi8}), $\tau_k$  is a periodic trajectory. Since the property of being periodic is preserved by the affine automorphisms $\gamma \nu_i $, $0\leq i \leq 7$, which act in the renormalization scheme, also $\tau$ is periodic. 
\end{proofof}

\section{Generation, coherence and characterization of the closure}\label{generationproofsec}
\subsection{Infinite coherence of cutting sequences}\label{infcoherencesec}
In this section we prove that cutting sequences are infinitely coherent (Proposition \ref{infcoherentprop}). The basic step is given by the following Lemma which is a consequence of a  further analysis of the geometric renormalization operation introduced in the proof of Proposition \ref{derivationiscuttingseqpropsec0}.  

\begin{lemma}\label{coherencecutseq1step}
Let $w=c(\tau)$ be a cutting sequence and let $i=s(\tau)$ and $j=s(n(\tau)')$. Then $w$ is coherent with respect to the pair $(i,j)$. 
\end{lemma}
\begin{proof}
Let us verify  conditions $ ($C$0) , ($C$1), ($C$2),$ and $($C$3)$ in the Definition \ref{coherencedef} of coherence.  
Since $s(\tau)=i$, $w$ is admissible in diagram $i$ and $($C$0)$ holds. Moreover, since  $n(w)'=c(n(\tau)' )$ (see for example Lemma \ref{cutseqrelationlemma} for $k=1$) and we assume that $s(n(\tau)')=j$, $n(w)'$ is admissible in diagram $j$. By Lemma \ref{typesequenceslemma} $(i)$,  the sector $s(n(\tau)')=s_1(\tau)$ is not $0$, so $j \neq 0$ and $($C$2)$ also holds.  Since the conditions $($C$1)$ and $($C$4)$ are expressed through the normalized sequence $n(w)$, for notational convenience let us assume from now on that $\theta \in \Sigma_0$ so that $n(\tau)=\tau$ and $n(w)=w$. 

In the  proof of Proposition \ref{derivationiscuttingseqpropsec0} we showed that the sequence of sandwiched letters in $w$, up to replacing the the letters $B$, $C$ and $D$ by their primed counterparts, gives the cutting sequence of the cut and pasted trajectory $\Upsilon_0^{-1}\tau$ in $O'= \gamma\octcdot O$ with respect to the labels on the sides of $O'$ in Figure \ref{cutseqshearedoctagonfig}. 
Let $\overline{c}({\tau}) \in \{A,B',C',D',e,f,g,h,i\}^{\mathbb{Z}}$ be the augmented cutting sequence of $\Upsilon_0^{-1}\tau$ with respect to both sides of $O'$ and  auxiliary diagonals in Figure \ref{cutseqshearedoctagonfig}. As we showed in  Proposition \ref{derivationiscuttingseqpropsec0}, $\overline{c}(\tau) $ is obtained reading off the labels of both edges and vertices of an infinite path on diagram $\mathscr{D}_0'$ in Figure \ref{dualDfig}. 

Let $\theta_{B'} , \theta_{C'} $ and $\theta_{D'}$ be the directions of the sides of $O'$ labelled by $B', C', D'$ respectively. Consider for each primed letter, for example $B'$, the parallelogram in Figure \ref{cutseqshearedoctagonfig} (or equivalently  the rectangle in Figure \ref{rectanglesedecompfig}) which has $B'$ as diagonal and whose sides are labeled by the auxiliary letters. If $\theta < \theta_{B'} $ (and hence in Figure \ref{rectanglesedecompfig} the slope of the trajectory is less  than the slope of the diagonal), $B'$ is preceded and followed in $\overline{c}({\tau})$  by $i$ and $e$ respectively. This identifies the corresponding $B'$ arrow in diagram ${\mathscr{D}}_0'$ in Figure \ref{dualDnewoctagonfig}, the one from $i$ to $e$. Since this arrow in diagram $\widetilde{\mathscr{D}}_0'$ in Figure \ref{dualDfig} is preceded and followed by arrows with the label $D$,  this shows that if $\theta < \theta_{B'} $ then all sandwiched $B$'s in $w$ are $D$-sandwiched. Similarly, if $\theta > \theta_{B'} $, the corresponding auxiliary letter transitions is $h f$ (see Figure \ref{cutseqshearedoctagonfig} or \ref{rectanglesedecompfig}) and hence, looking at the $B$-labelled arrow from $h$ to $f$ in
diagram $\widetilde{\mathscr{D}}_0'$ in Figure \ref{dualDfig}, $B$ in this case is always $C$-sandwiched. Reasoning in the same way, if $\theta < \theta_{C'} $ then $C$ is $C$-sandwiched and $C$ is $B$-sandwiched if $\theta > \theta_{C'}$. We see that $D$ is $B$-sandwiched if $\theta < \theta_{D'} $ and $D$ is $A$-sandwiched if $\theta > \theta_{D'}$.
We also observe that diagram $\widetilde{\mathscr{D}}_0'$ in Figure \ref{dualDfig} shows that $A$ is always $D$-sandwiched.  

Since $0< \theta_{B'} < \theta_{C'} < \theta_{D'}< \pi/8$ (see Figure \ref{cutseqshearedoctagonfig}), these inequalities force certain implications  between the sandwiched letters combinations that can occur in a given cutting sequence. For example if $0\leq \theta < \theta_{B'}$ then it follows that $\theta < \theta_{C'} < \theta_{D'}$. Thus we see from the above considerations that all $B$'s are $D$-sandwiched, all $C$'s are $C$-sandwiched and all $D$'s are $B$-sandwiched. This, recalling that $A$ is always $D$-sandwiched, gives the combinations of letters in group $G_3$ in Definition \ref{coherencedef}. The other groups correspond to the other  cases,
as summarized in the Table \ref{coherencetable}. This concludes the proof of condition $($C$1)$ of Definition \ref{coherencedef}.
\begin{table} \centering
\begin{tabular}{c| c||c|c|c|| c}
$\theta$ & $\theta'$ & B & C & D &  \\
\hline\hline \vspace{1mm}
$0\leq \theta \leq  \theta_{B'}$ & $\theta' \in \overline{\Sigma}_6  \cup \overline{\Sigma}_7$ & $D$-sandwiched & $C$-sandwiched & $B$-sandwiched & $G_3$ \\ 
\hline  \vspace{1mm}
$\theta_{B'} \leq \theta \leq  \theta_{C'} $& $\theta' \in \overline{\Sigma}_4 \cup \overline{\Sigma}_5$  & $C$-sandwiched & $C$-sandwiched & $B$-sandwiched & $G_2$ \\ 
\hline \vspace{1mm}
$ \theta_{C'} \leq \theta \leq  \theta_{D'}$ & $\theta' \in \overline{\Sigma}_2 \cup \overline{\Sigma}_3$ & $C$-sandwiched & $B$-sandwiched & $B$-sandwiched & $G_1$ \\ 
\hline \vspace{1mm}
$ \theta_{D'} \leq \theta \leq \pi/8 $&
$\theta' \in \overline{\Sigma}_1$ &  $C$-sandwiched & $B$-sandwiched& $A$-sandwiched &$ G_0 $ \\ 
\end{tabular}
\caption{Coherence of sandwiched letters according to sectors in $[0,\pi/8]$.\label{coherencetable}}
\end{table}

Let us verify condition $($C$3)$ of Definition \ref{coherencedef}. 
We observe that the diagram $j$ in which $n(w)'=w'$ is admissible is determined by the sector $\Sigma_j$ to which the direction $\theta'$ of $\tau'= n(\tau)'$ belongs. On the other hand, as we have just seen, the index $k$ of group $G_k$ in $($C$1)$ is determined by the inequalities between $\theta$ and $ \theta_{B'}, \theta_{C'}, \theta_{D'}$. The angles $\theta$ and $\theta'$ are related by $\theta'= \gamma [\theta]$ (see Proposition \ref{typesequalCF}). 
One can check (for example using that $O=\gamma \octcdot O'$ and analyzing the directions of corresponding sides) that:

\be\label{imagesangles}
\gamma [ 0]=\pi, \quad    \gamma [ \theta_{B'}] = {3\pi/ 4},  \quad      \gamma [ \theta_{C'}]=  {\pi}/{2},  \quad    \gamma [ \theta_{D'}]=  {\pi}/{4},  \quad     \gamma [ \pi/8] = {\pi}/{8}. 
\ee  
Since $\gamma$ is orientation reversing, this gives the following correspondence,  illustrated in Figure \ref{gammasectorsfig}:
\be
\label{imagessectors}
\begin{array}{lcllcl}
 \gamma\left[\phantom{\frac{}{}} [0, \theta_{B'}] \phantom{\frac{}{}} \right]&=& \overline{\Sigma}_6\cup \overline{\Sigma}_7 ,  &
  \gamma \left[\phantom{\frac{}{}} [\theta_{B'}, \theta_{C'}] \phantom{\frac{}{}} \right]& = & \overline{\Sigma}_4 \cup \overline{\Sigma}_5,   \\ 
   \gamma\left[ \phantom{\frac{}{}} [\theta_{C'} , \theta_{D'}] \phantom{\frac{}{}} \right]& = & \overline{\Sigma}_2 \cup \overline{\Sigma}_3, \phantom{\frac{\pi}{8}}   & \gamma\left[ \left[ \theta_{D'} ,  \frac{\pi}{8}\right] \phantom{\frac{}{}}\right]&= & \overline{\Sigma}_1 .
\end{array}\ee
\noindent From  the first equality, i.e.~$\gamma\left[\phantom{\frac{}{}} [0, \theta_{B'}] \phantom{\frac{}{}} \right]= \overline{\Sigma}_6\cup \overline{\Sigma}_7$, we have that $\theta \in [0,\theta_{B'}]$, and hence $k=3$, if and only if $\theta' \in \overline{\Sigma}_6 \cup \overline{\Sigma}_7$, and hence $j=6$ or $j=7$. In this case,  $k=[j/2]$. The other  cases are similar and  are summarized in Table \ref{coherencetable}, where  one can check that the relation $k=[j/2]$ always holds. This concludes  the proof of $($C$3)$ and of the Lemma. 
\begin{figure}[!h]
\centering
\includegraphics[width=.9\textwidth]{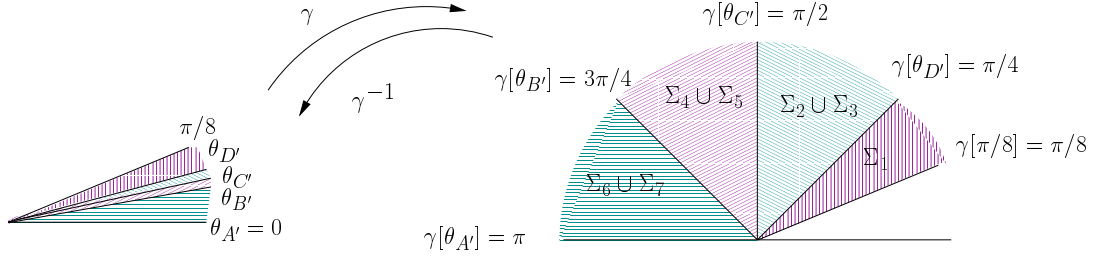}
\caption{The action of $\gamma=\gamma^{-1}$  on sectors of directions.\label{gammasectorsfig}}
\end{figure}
\end{proof}

\begin{proofof}{Proposition}{infcoherentprop} 
Let $w=c(\tau)$  be a a cutting sequence  and 
let $\{ s_k \}_{k\in\mathbb{N}}$ be the sequence of sectors of $\tau$ (i.e.~$s_k=s_k(\tau)$, see Definition \ref{sectorsdef}).  Let us show by induction that the $k^{th}$ renormalized word $w_{k}$ (in Definition \ref{combrecursive}) is coherent with respect to $(s_k,s_{k+1})$ for each $k \in \mathbb{N}$. 
Let $\tau_k$ be the sequence of renormalized trajectories (see Definition \ref{renormalizationtrajectoriesdef}) and let us first verify\footnote{We remark that we cannot use Lemma \ref{cutseqrelationlemma} here since we are considering  periodic words for which the renormalization scheme in \S\ref{combrenormsec} is not well-defined. Moreover, if $w$ is not periodic, we have to verify that the sequence of renormalized words $w_k$ in Definition \ref{combrecursive} (renormalized with respect to an abstract sequence $s$), coincides with the words $w_k$ of the
renormalization scheme in \S\ref{combrenormsec} when $s$ is the sequence of sectors.} by induction on $k$ that $c(\tau_k)= w_k$. Since  $s_k$ is the sector of $\tau_k$, $n(\tau_k)=\nu_{s_k} \trajcdot \tau_k$ and by definition  $\tau_{k+1} = n(\tau_k)'$, so that $\tau_{k+1} = (\nu_{s_k} \trajcdot \tau_k)'$. 
Thus, if we know by inductive assumption that  $w_k = c(\tau_k)$, using (\ref{derivedseqeq}), (\ref{relabellingeq}) and Definition \ref{combrecursive}, we get 
\bes c (\tau_{k+1}) = c\left( (\nu_{s_k} \trajcdot \tau_k)'\right)= c(\nu_{s_k} \trajcdot \tau_k )' = \left(\pi_{s_k} \cdot c(\tau_k)\right)' = (\pi_{s_k} \cdot w_k )' =w_{k+1} . \ees
Thus, for each $k\in \mathbb{N}$, we have $w_k = c(\tau_k) $ and  since $s(\tau_{k})=s_{k}$ and  $s(\tau_{k+1})=s(n(\tau_k)')= s_{k+1}$, we can apply   Lemma \ref{coherencecutseq1step} to $w_k = c(\tau_k) $, with $i =s_k$ and $j=s_{k+1}$. Lemma \ref{coherencecutseq1step} gives that $w_k $ is coherent with respect to $(s_k,s_{k+1})$, concluding the proof.   
\end{proofof}

\subsection{Finding generation rules and inverting derivation}\label{generationoperationsec}
Just as the proof of Proposition \ref{derivationiscuttingseqprop} gave a geometric counterpart to derivation,  the following Proposition \ref{generatetype0} gives a geometric counterpart to the operation of generation. Moreover, the proof of Proposition \ref{generatetype0} shows how the interpolating words and the generating rules in Figure \ref{generationrulesfig} can be explicitly found from the analysis of the renormalization operation described in  Proposition \ref{derivationiscuttingseqpropsec0}.
\begin{prop}\label{generatetype0}
Let  $w=c(\tau)$ be the cutting sequence of a trajectory $\tau$ and assume that $s(\tau)=k \neq 0$. Let $W$  be the cutting sequence of the trajectory  $\Psi_\gamma \tau$. Then  $W = \g{k}{0}{w} $. In particular, $W$ is admissible in diagram $0$ and $W'=w$. 
\end{prop}
\noindent Proposition \ref{generatetype0} shows that when $w$ is a cutting sequence, $\g{k}{0}{w}$ is also a cutting sequence and in particular that it is the cutting sequence of the trajectory  $\Psi_\gamma \tau$.  More generally:
\begin{cor}\label{gopgeometrically}
Let $w$ be as in Proposition \ref{generatetype0}. The generated word $\g{k}{i}{w}$ is the cutting sequence of the trajectory $\nu_i^{-1} \trajcdot \Psi_{\gamma}\tau$. 
\end{cor}
\noindent We also remark that the operation of derivation is not uniquely invertible, since we have the following. 
\begin{cor}\label{7to1}
Given a cutting sequence $w = c(\tau)$ where $s(\tau)=0$, there exists cutting sequences $W_1$, $W_2, \dots, W_7$ such that  $W_j$ is admissible in diagram $j$  and  $W_j'=W$. 
\end{cor}
\noindent The Corollaries \ref{gopgeometrically}  and \ref{7to1} are both proved in this section, after the proof of Proposition \ref{generatetype0}.

\begin{proofof}{Proposition}{generatetype0}
Let us set $\tilde{\tau}:= \Psi_\gamma  \tau $, so we that, since $ \Psi_\gamma  =  \Psi_\gamma ^{-1}$ by Corollary \ref{Psigammainverse}, we can think of  $\tau$ as $\Psi_{\gamma} \tilde{\tau} $. Let use denote by $\tilde{\theta}$ be the direction of $ \tilde{\tau} $   and by  $W$ be the cutting sequence of $\tilde{\tau}$ with respect to $O$. Since by assumption $\theta \in [\pi/8,\pi]$  and since $\tilde{\tau} = \Psi_\gamma
\tau$ implies that $\tilde{\theta} =\gamma[\theta]$, it follows that $\tilde{\theta}\in \gamma \left[ [\pi/8,\pi]  \right] = \overline{\Sigma}_0$. Thus, by Proposition \ref{derivationiscuttingseqpropsec0}, $W' = w$  and moreover $w$, if we prime the letters $B,C,D$, is also the cutting sequence of the trajectory $\Upsilon_o^{-1} \tilde{\tau}$ ontained by cut and paste (which is trajectory in $O'$ in the same direction $\tilde\theta$)  with respect to the sides of $O'=\gamma O$ labelled as  in Figure \ref{cutseqshearedoctagonfig}. We will show that we can recover the sequence $W$ from the knowledge of $w$ using generation.  

Consider  the auxiliary diagonals labelled by $\{ e,f,g,h,i\}$  inside $O'$ (refer to  Figure \ref{cutseqshearedoctagonfig}) and let $\overline{c}(\tilde{\tau})\in \{ A,B',C',D',e,f,g,h,i\}^{\mathbb{Z}}$ be the augmented cutting sequence of $\Upsilon_o^{-1} \tilde{\tau}$  with respect to both the sides of $O'$ and the auxiliary diagonals. 
Let $0= \theta_{A'} < \theta_{B'} < \theta_{C'} < \theta_{D'}< \pi/8$ denote the directions of the sides of $O'$ labelled by $A, B', C', D'$ respectively (see Figure \ref{cutseqshearedoctagonfig}).  From Figure \ref{cutseqshearedoctagonfig}, reasoning as in the proof of Lemma \ref{coherencecutseq1step}, we see  that
if $\tilde\theta<  \theta_{B'}$, then $B'$ is preceded and followed in $\overline{c}(\tilde{\tau})$ by $i$ and $e$ respectively, while if  $\tilde\theta >  \theta_{B'}$,  $B'$ is always preceded and followed  by $h$ and $f$ respectively. Similar conclusions can be obtained from Figure \ref{cutseqshearedoctagonfig} for $C'$ and $D'$. We  summarize the conclusions  in the last two columns of Table \ref{auxiliaryletterstable}. From  (\ref{imagessectors}), recalling that $\gamma$ is an involution,  
 we have the following correspondence of directions under the action of $\gamma$ (see Figure \ref{gammasectorsfig}):
\bes
\begin{array}{lcllcl}
 \gamma\left[\overline{\Sigma}_1\right]&=& \left[  \theta_{D'} ,  {\pi}/{8}\right], &
  \gamma\left[\overline{\Sigma}_2 \cup \overline{\Sigma}_3 \right] &=& \left[  \theta_{C'} , \theta_{D'} \right],  \\ 
  \gamma\left[\overline{\Sigma}_4 \cup \overline{\Sigma}_5 \right]&=& \left[ \theta_{B'},  
\theta_{C'}  \right], & \gamma\left[ \overline{\Sigma}_6\cup \overline{\Sigma}_7\right] &=& \left[ 0, \theta_{B'}  \right].
\end{array}
\ees
\noindent Using this correspondence and knowing the sector\footnote{Let us remark that it is enough to know less than the sector $\overline{\Sigma}_k$ (or equivalently the type of $w$),~i.e. it is enough to know only in which of the four larger sectors $\overline{\Sigma}_1$, $\overline{\Sigma}_{2}\cup \overline{\Sigma}_3$, $\overline{\Sigma}_{4}\cup \overline{\Sigma}_5$ or  $\overline{\Sigma}_{6}\cup \overline{\Sigma}_7$ the direction $\theta$ belongs.} of $\theta$, we can recover the order relation  between  $\tilde\theta$ and  $\theta_{A'}, \theta_{B'}, \theta_{C'}, \theta_{D'}$ and thus  we can recover the auxiliary letters that precede and follow in $\overline{c}(\tilde{\tau})$ each letter of $w$ (with primed letters $B',C',D'$ for $B,C,D$),  as summarized in the Table \ref{auxiliaryletterstable}.
\begin{table} \centering
\begin{tabular}{c| l| l||c|c|c }
 Fig.~\ref{graphs4sectorsfig}& $\theta$-sector & $\theta'$-sector & B$'$ & C$'$ & D$'$ \\
\hline\hline  
$\mathscr{G}_0$& $\theta \in \overline{\Sigma}_1$ &  $ \theta_{D'} < \tilde\theta \leq \pi/8 $ & h B$'$ f  & f C$'$ h  & AD$'$A \\
\hline 
$\mathscr{G}_1$& $\theta \in \overline{\Sigma}_2 \cup \overline{\Sigma}_3$ & $ \theta_{C'} \leq  \tilde\theta \leq  \theta_{D'}$& h B$'$ f& f C$'$ h & e D$'$ i \\
\hline  
$\mathscr{G}_2$&
$\theta \in \overline{\Sigma}_4 \cup \overline{\Sigma}_5$  & $\theta_{B'} \leq \tilde\theta \leq \theta_{C'} $ & h B$'$ f& g C$'$ g & e D$'$ i \\
\hline  
$\mathscr{G}_3$&
$\theta \in \overline{\Sigma}_6 \cup \overline{\Sigma}_7$ &$0\leq \tilde\theta \leq \theta_{B'}$ &  i B$'$ e & g C$'$ g  & e D$'$ i \\
\end{tabular}
\caption{Auxiliary letters preceding and following each letter of $\overline{c}(\tilde{\tau})$.\label{auxiliaryletterstable}}
\end{table}

To conclude the proof, we want to describe the generation rules that allow us to recover the word $W$ from $w$.
 We know from the proof of Proposition \ref{derivationiscuttingseqpropsec0} that $W$  is obtained reading the labels in diagram $\widetilde{\mathscr{D}}_0'$ (Figure \ref{dualDfig}) of the  path given by   $\overline{c}(\tilde{\tau})$   on $\mathscr{D}_0'$ (Figure \ref{dualDnewoctagonfig}). 
\begin{figure}[!h]
\centering
\subfigure[$\mathscr{G}_0$ for $\theta \in \overline{\Sigma}_1$\label{G0fig}]{\includegraphics[width=0.42\textwidth]{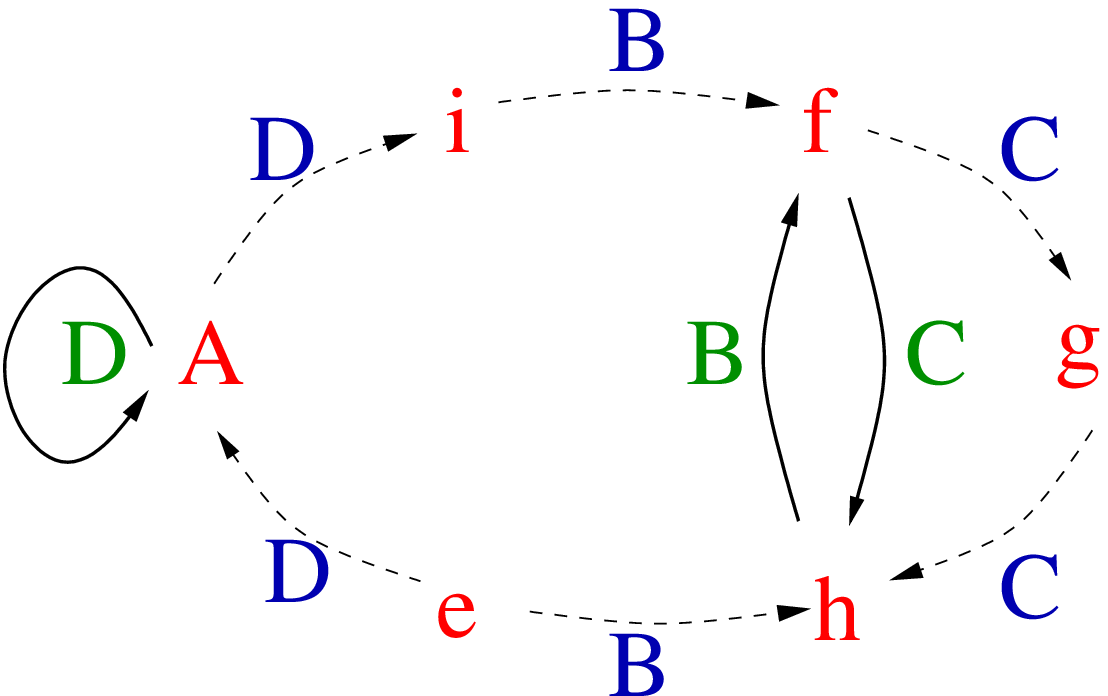}}
\hspace{3mm}
\subfigure[$\mathscr{G}_1$ \label{G1fig} for $\theta \in \overline{\Sigma}_2 \cup \overline{\Sigma}_3$ ]{\includegraphics[width=0.42\textwidth]{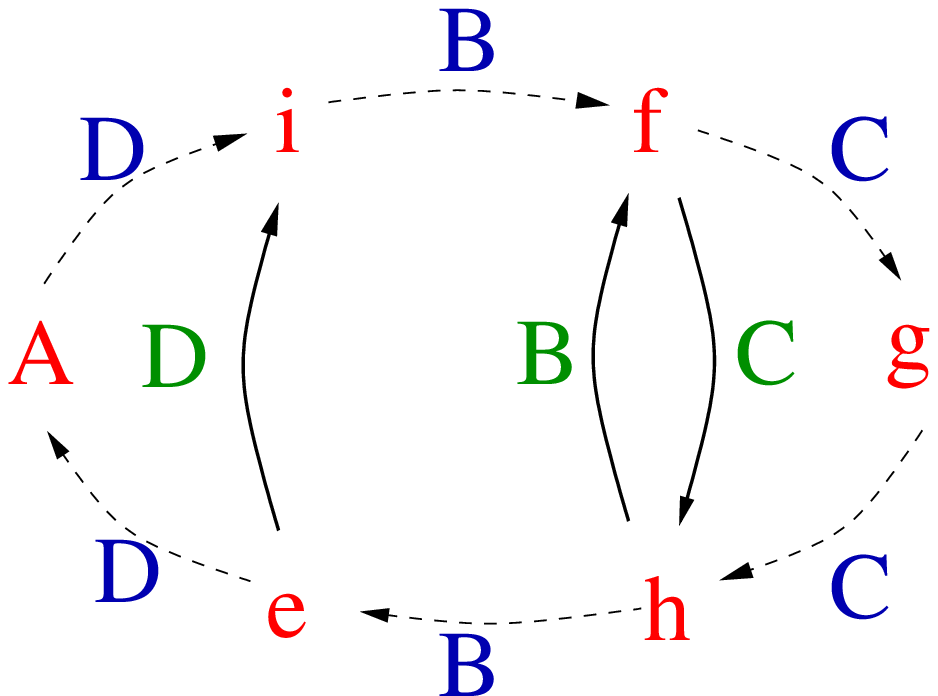}}
\subfigure[$\mathscr{G}_2$ \label{G2fig} for $\theta \in \overline{\Sigma}_4 \cup \overline{\Sigma}_5$]{\includegraphics[width=0.42\textwidth]{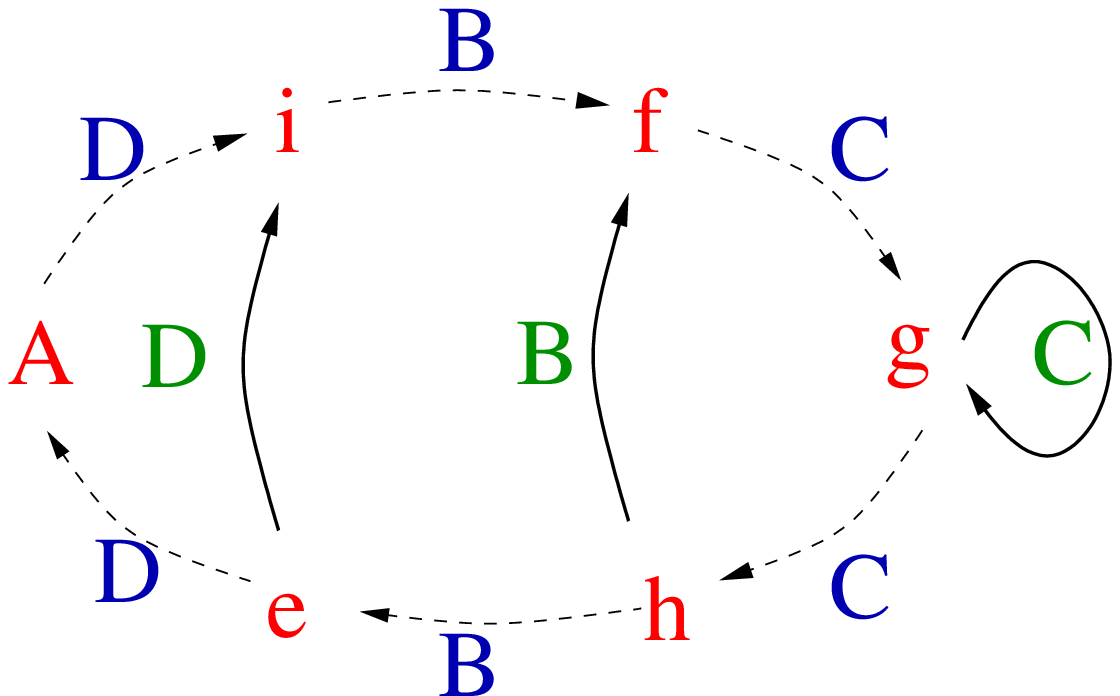}}
\hspace{3mm}
\subfigure[$\mathscr{G}_3$ for $\theta \in \overline{\Sigma}_6 \cup \overline{\Sigma}_7$]{\includegraphics[width=0.42\textwidth]{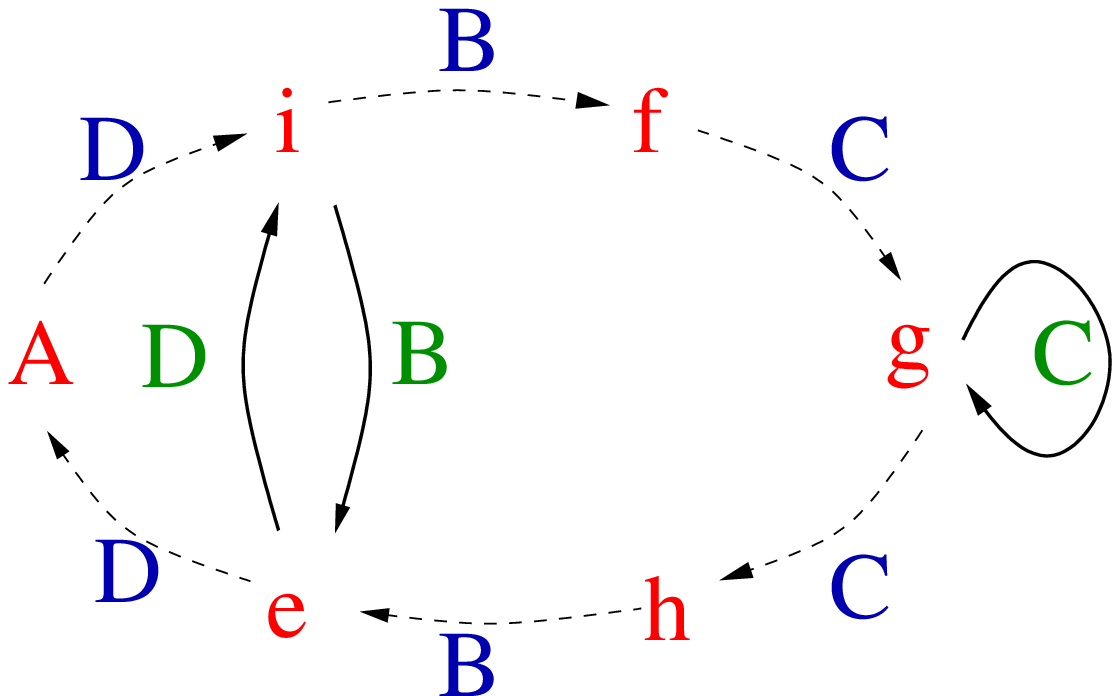}}
\caption{Graphs used to recover the generation rules. \label{graphs4sectorsfig}}
\end{figure}
Observe that in $\mathscr{D}_0'$ there are multiple arrows labelled with the same letter (for example both the edges from $i$ to $e$ and from $h$ to $f$ are labelled by $B'$). Which arrow is used depends only on $\theta$ and can be determined by using Table \ref{auxiliaryletterstable}. For example,   if $\theta<  \theta_{B'}$, looking at the last line of Table \ref{auxiliaryletterstable} we know that  $B'$ is preceded by $i$ and followed by $e$ and  this shows that $B'$-arrow used is the one from $i$ to $e$. Hence, according to the four  sectors in Table \ref{auxiliaryletterstable}, the set of  arrows whose labels give  the letters $B,C,D$ in $W$ is the set of solid-line arrows in one  of the four graphs $\mathscr{G}_0$,   $\mathscr{G}_1$, $\mathscr{G}_2$ or $\mathscr{G}_3$ in Figure \ref{graphs4sectorsfig}. 

One can check that the graphs $\mathscr{G}_0$,   $\mathscr{G}_1$, $\mathscr{G}_2$ or $\mathscr{G}_3$
have the following property. Given any two solid arrows,  there is a unique way of going from the ending vertex of one to the entering vertex of the other using only dashed arrows and not going through the $A$ vertex. Similarly, given $A$ and a solid arrow, there is a unique way of going from $A$ to the entering vertex of the arrow or from the exiting vertex of the arrow to the $A$ vertex without using solid arrows. 
Thus, given two letters $L_1, L_2 \in \{A,B,C,D\}$, let $v_{L_i}$, $i=1,2$, be the vertex labelled by $A$ if $L_i=A$ and let otherwise $v_{L_1}$ be the exiting vertex of the solid arrow labelled  $L_1$ and  $v_{L_2}$ be the entering vertex of the solid arrow labelled by $L_2$. 
Let  $p_{L_1 L_2}$ be the unique path  which goes from $v_{L_1}$ to $v_{L_2} $ and does not use solid arrows nor goes through the $A$ vertex.   
Reading off the labels of the sides of $p_{L_1 L_2}$,  we obtain precisely the generation words $w_{L_1L_2}$  in Figure \ref{generationrulesfig}. In order to find the labels of one the edges of $\mathscr{D}_k$,  we use the appropriate graph for $\Sigma_k$, which (according to the correspondence summarized in Table \ref{auxiliaryletterstable}) is $\mathscr{G}_{[k/2]}$. For example, to interpolate between the solid arrow labelled $D$   and the solid arrow labelled $B$ in $\mathscr{G}_1$ (see Figure \ref{G1fig}), one has to go through dashed arrows labelled in order $BCC$ and, since $[3/2]=1$, this gives exactly the label of the $DB$-edge in  $\mathscr{D}_3$  in Figure \ref{generationrulesfig}(c). We leave to the  reader the verification of  the other interpolation words.  
\end{proofof}

\begin{proofof}{Corollary}{gopgeometrically}
By Proposition \ref{generatetype0}, $\g{k}{0}{w}$ is the cutting sequence of $\Psi_{\gamma}\tau$. Since $\g{k}{i}{w} = \pi_i^{-1 } \g{k}{0}{w}$ (see  Definition \ref{othergoperators}), the Corollary follows from (\ref{relabellingeq}). 
\end{proofof}

\begin{proofof}{Corollary}{7to1}
For $k=1, \dots, 7$, consider the cutting sequences $\pi_k^{-1}\cdot w$  of the trajectories $\nu_k^{-1} \trajcdot \tau$, whose sectors are all different from $0$,  since $s(\tau)=0$ and $\nu_k^{-1}$ maps $\Sigma_0 $ to $\Sigma_k$. Thus,  for each  $k=1, \dots, 7$, we can apply Proposition \ref{generatetype0} to $\pi_k^{-1}\cdot w$ and obtain a cutting sequence $\widetilde{W}_k$  of type $0$ such that $(\widetilde{W}_k)' = \pi_k^{-1}\cdot w$. 
Let us set $W_k: = \pi_k \cdot \widetilde{W}_k$.  Thus, for each $k=1, \dots, 7$, we have that $W_k$ is  of type $k$ by construction and that  $(W_k)'= w$ since derivation commutes with the action of permutations (see Remark \ref{commute}).  
\end{proofof}

\subsection{Generating rules and coherence}\label{generationandcoherencesec}
In this section we prove Lemma \ref{coherencegeneration} and Proposition \ref{generationandinfcoherence} and  show that the coherence conditions allow us  to fully recover the generation rules. Let us remark  that in Proposition \ref{generatetype0} we showed how to recover generation rules for words which are cutting sequences and thus we could use geometric considerations. Here we assume only that $w$ is an infinitely coherent word.   In this section we also give the proof of Lemma \ref{uniqueness}, which is based on Lemma \ref{coherencegeneration}.

\begin{proofof}{Lemma}{coherencegeneration}
Let us first show that if $w=\g{j}{i}{v}$,  where $v\in \mathscr{A}_j$,  then $v=n(w)'$. Indeed, since $w$ is of type $i$, by Definition \ref{othergoperators} and by Lemma \ref{g0invertgeneration}, we have that   $n(w)'= (\pi_i \cdot \g{j}{i}{v})' =  (\g{j}{0}{v})'= v  $. 

Let us now show that if $w$ can be written as $\g{j}{i}{v}$ where $v\in \mathscr{A}_j$, then $w$ is coherent with respect with $(i,j)$. 
 Let us verify the conditions in Definition \ref{coherencedef}.  
Clearly $w$ is of type $i$ and, as we just proved,    $n(w)'=  v  $ where $v$ is by assumption is admissible in diagram $j$, so that $($C$0)$ and   $($C$2)$ in Definition \ref{coherencedef} of coherence are  satisfied. Let us remark that $n(w)=  \g{0}{j}{v} $ and that, by Lemma \ref{g0invertgeneration} and Definition \ref{generationopto0}, the sandwiched letters in $n(w)$ are exactly the letters of $v$, hence are given by the labels of the vertices of the path described by $n(w)$ on  $\mathscr{D}_j$ in Figure \ref{generationrulesfig}. 
 Thus, to check $($C$1)$ and $($C$3)$, it is enough to verify that, for each $j=1,\dots, 7$,  the generation rules in diagram $\mathscr{D}_j$ in  Figure \ref{generationrulesfig} are such that the letters corresponding to vertices are sandwiched exactly as in the group $G_{[j/2]}$ in Definition \ref{coherencedef}.  

Now we prove that if $w$ is coherent with respect with $(i,j)$, then, if we set $W:=n(w)$ and $v:= W'=n(w)'$,  $w$  can be written in the form $\g{i}{j}{v}$ where $v = n(w)'\in \mathscr{A}_j$.    Given  $w$ coherent with respect to $(i,j)$,  We know by $($C$2)$ in Definition \ref{coherencedef} that $v$ is of type $j$.  
Let us show that $W = \g{j}{0}{v}$. Since $W=n(w) = \pi_i \cdot w$ because $w$ is of type $i$ (by $($C$0)$ in Definition \ref{coherencedef}), this is enough to show that  $w =  \pi_i^{-1} \cdot \g{j}{0}{v}   =  \g{j}{i}{v} $ (recall Definition \ref{othergoperators}) and hence conclude the proof. 

Since $W$ is admissible in diagram $0$ and its derived sequence is $v$, 
 let us start by listing all possible ways of interpolating a transition in $v$ by a path on  $\mathscr{D}_0$ which does not contain any other sandwiched letter.  
 Because of the structure of  $\mathscr{D}_0$ (refer to Figure \ref{D0}), the requirement of avoiding sandwiched letters leaves only the possibilities listed in Table \ref{transitionsinterpolations} (where the letters which belong to  $v$ are written in bold fonts to distinguish them from the  interpolating letters). Let us  justify for example the entry corresponding to the transition $AC$, which shows all ways of interpolating between a sandwiched $A$ and a sandwiched $C$. A path in $\mathscr{D}_0$ which starts from  $A$, goes to $D$: if after $D$ there were another $A$ then $D$ would be sandwiched, hence there has to be a $B$ and similarly after $B$ there is a $C$ (otherwise the $B$ would be sandwiched). If after $C$ there is another $B$ then $C$ is sandwiched and thus the interpolation word is complete and is $DB$. There could also be another $C$. In this case, there has to  be an additional $C$, so that the first $C$ is $C$-sandwiched: $ADBCCB$ cannot be completed in any other way, since $ADBCCBC$, $ADBCCBDB$ and $ADBCCBDAD$ all contain another sandwiched letter ($B$, $D$ or $A$ respectively).  

\begin{table}
\centering
\begin{tabular}{c||l|l|l|l |}& {\A} & \B &\C & \D \\
\hline\hline
\multirow{2}{*}{\A} &	 \multirow{2}{*}{\A DBCCBD\A} & \A D\B & \A DB\C & 	\A \D \\
& 	& \A DBCC\B & \A DBC\C 		& \A DBCCB\D\\
\hline	
\multirow{2}{*}{\B} &	 \B D\A& \multirow{2}{*}{\B CC\B} 	& \B \C & \B \D \\
& \B CCBD\A 	&   & \B C\C & \B CCB\D \\
\hline
\multirow{2}{*}{\C} & \C BD\A& 	\C \B	  & \multirow{2}{*}{\C \C}&  \C B\D \\
& \C CBD\A&		\C C\B	&	& \C CB\D \\
\hline
\multirow{2}{*}{\D}	& \D \A&	\D \B  & \D B\C    & \multirow{2}{*}{	\D BCCB\D }\\
& \D BCCBD\A 	& \D BCC\B      	& \D BC\C		& \\
\hline 
 \end{tabular}
\caption{Possible interpolations of type $0$ of sandwiched transitions.\label{transitionsinterpolations}}
\end{table}

Since $ v= W'= n(w)'$ and $v$ is admissible in diagram $j$, conditions $($C$1)$ and $($C$3)$ in Definition \ref{coherencedef} imply that all sandwiched letters in $W$ belong to $G_{[j/2]}$. This requirement determines the interpolation, since when in an entry in Table  \ref{transitionsinterpolations} there appear two possible interpolations, they correspond to two different ways of sandwiching a  letter: for example, there are two possible  interpolations of the transition $AC$, but $ADBC$ can be used only if $C$ is $B$-sandwiched and hence belongs to the groups $G_0$ or $G_1$, while  $ADBCC$ only if $C$ is $C$-sandwiched and hence in $G_2$ or $G_3$. This leads once again to the generation rules in Figure \ref{generationrulesfig} and thus (recalling Definition \ref{generationopto0}) shows that $W=\g{j}{0}{v}$.  
\end{proofof}

\begin{proofof}{Proposition}{generationandinfcoherence}
Let us first show that, given any $s_0, \dots, s_{k} \in \{0, \dots, 7\}$ and $w=w_0$ and letting $w_l$, for $l=1, \dots, k$,  be the renormalized words in Definition \ref{combrecursive}, $w_l$  is $(s_l,s_{l+1})$-coherent for all $l=0, \dots, k-1$ if and only if   
\be\label{finiteintersection}
w = \g{s_1}{s_0}{ \g{s_2}{s_1}{ \ldots  \g{s_{k-1}}{s_{k-2}}{ \g{s_{k}}{s_{k-1}}{w_k}}}}. 
\ee
If (\ref{finiteintersection}) holds, using that if $w=\g{j}{i}{v}$  then $v=n(w)'$ (see Lemma \ref{coherencegeneration}) and Definition \ref{combrecursive}, one can see by induction on $l=0,\dots , k-1$ that  $w_l = \g{s_{l+1}}{s_{l}}{ \ldots  \g{s_{k}}{s_{k-1}}{w_k}}$. Thus,  $w_l$  is $(s_l,s_{l+1})$-coherent by  Lemma \ref{coherencegeneration}. Conversely, let us assume that $w_l$  is $(s_l,s_{l+1})$-coherent for all $l=0, \dots, k-1$ and prove (\ref{finiteintersection}) by induction on $l$. The inductive assumption gives $w = \g{s_{1}}{s_{0}}{ \ldots  \g{s_{l}}{s_{l-1}}{w_l}}$. Since $w_{l}$ is $(s_{l},s_{l+1})$-coherent, by Lemma \ref{coherencegeneration}, $w_{l} =   \g{s_{l+1}}{s_{l}}{ (n(w_l)')}$ and by Definition \ref{combrecursive}, this gives $w = \g{s_{1}}{s_{0}}{ \ldots  \g{s_{l+1}}{s_{l}}{w_{l+1}}}$. 

It follows that $w$ is infinitely coherent with respect with $\{s_k\}_{k\in \mathbb{N}}$ iff for each $k\in \mathbb{N}$ $w$ can be expressed as (\ref{finiteintersection}), which equivalently means that $w \in \mathscr{G}(s_0, \dots, s_k)$. Hence $w$ is infinitely coherent with respect with $\{s_k\}_{k\in \mathbb{N}}$ iff  it belongs to the intersection (\ref{infiniteintersection}). 

Let us now prove the equivalence with (\ref{unionfiniteintersection}). If $w$ is infinitely coherent, we just showed that it belongs to the intersection (\ref{infiniteintersection}) for some $s\in S^{*}$, thus it is in (\ref{unionfiniteintersection}). Conversely, let  $w$  be in (\ref{unionfiniteintersection}). Then for each $k\in \mathbb{N}$ there exists $s_0^k,s_1^k, \dots, s_k^k$ such that  $w \in \mathscr{G}(s^k_0,\dots, s^k_k)$, but, fixed $i$,  $s_i^k$ could a priori depend on $k$. Since $s_0^k $ can assume only finitely many values, there exist $\overline{s}_0\in \{0,\dots, 7\}$ and a subsequence $\{k_j\}_{j\in \mathbb{N}}$ such that  $s_0^{k_j} =\overline{s}_0$ for all $j$. The proof of Lemma \ref{uniqueness} shows that, for a fixed $j$, the sequence $s_i^{k_j}$, $i=0, \dots, k_j$, is determined by $s_0^{k_j} =\overline{s}_0$ up to the ambiguity described in Lemma \ref{uniqueness} $(ii)$. Thus, passing to a further subsequence, there exists $\overline{s} = \{ \overline{s}_k\}_{k\in\mathbb{N}}\in S^*$ such that $w \in  \mathscr{G}(\overline{s}_0,\dots, \overline{s}_{k_j})$ for all $j\in\mathbb{N}$. Since the sets $\mathscr{G}(\overline{s}_0,\dots, \overline{s}_k)$ are nested, $w$ is in their intersection over all $k$.  Thus, by the characterization proved in the first part of the argument,  $w$ is infinitely coherent with respect to $\overline{s}$.  
\end{proofof}

\begin{proofof}{Lemma}{uniqueness}
Let us first prove that if $w$ is coherent with respect to  $(i,j)$, then $j$ is uniquely determined by $i$ unless $n(w)'$ is periodic of period $1$ or $2$. Equivalently, by Lemma \ref{coherencegeneration}, we have to prove that if $w= \g{j}{i}{v}$
for some $v \in \mathscr{A}_j$, $1\leq j\leq 7$, then $j$ is uniquely determined by $i$ unless $v$ is periodic of period $1$ or $2$.
Observe that if the word $v=n(w)'$  is not periodic of period $1$ or $2$, it contains at least two distinct transitions of the form $L_1L_2$ and $L_2L_3$, where
 $\{ L_1,L_2\} \neq \{ L_2, L_3\}$. From $w= \g{i}{j}v$, we can read off the corresponding interpolating words $w_{L_1, L_2}$ and $w_{L_2, L_3}$.
 
From Figure \ref{generationrulesfig},  one can check that 
given two interpolating words $w = w_{L_1 , L_2}$ and $w'= w_{L_1', L_2'}$ (possibly empty words), if $\{L_1,L_2\} \neq \{ L'_1,L'_2\}$, there is a unique diagram $\mathscr{D}_j$ in which  both $w$ and $w'$ appear as interpolating words. We leave this verification to the reader. Clearly, the same assertion remains true if we apply $\pi_i^{-1}$ to all the interpolating words in Figure \ref{generationrulesfig} (which means that we compare all interpolating words than correspond to the operators $\gop{j}{i}$ for $j=1,\dots, 7$). 
Thus,  if $v$ is not periodic of period $1$ or $2$ and, as already remarked, we can read off from $w= \g{j}{i}{v}$  two   interpolating words $w_{L_1 , L_2}, w_{L_2, L_3}$, $\{L_1,L_2\} \neq \{ L_2,L_3\}$, occurring in $\g{i}{j}{v}$,   this shows that $j$ is determined. 



Let us now consider a word $w$ which is infinitely coherent with respect to $\{s_k \}_{k\in\mathbb{N}} \in S^*$.  
By Definition \ref{infcoherencedef}, each $w_k$ is  coherent with respect to $(s_k, s_{k+1})$. Thus, by induction, given $s_0$ each $s_{k}$   is uniquely  determined, unless $w_{k+1}= n(w_k)'$ is periodic of period $1$ or $2$. In this case, $w$  is also periodic since it is obtained by applying the composed operator $\gop{s_1}{s_0} \dots \g{s_{k+1}}{s_{k}}{w_{k+1}}$ which takes  periodic words to periodic words. Thus, if $w$ is not periodic, all $s_k$ with $k\geq 1 $ are uniquely determined by $s_0$, proving  $(i)$. 

 If on the other hand there exists a $\overline{k}$ such that $w_{\overline{k}}$ is   periodic of period  $1$ or $2$, assume that $\overline{k}$ is the first such $k$. Observe first that, reasoning as above, $s_1, \dots, s_{\overline{k}-1}$ are determined by $s_0$. 
Since periodic words of period $1$ or $2$ are fixed points under derivation, recalling Definition \ref{worddef} we have  $w_{\overline{k}+1} = ( \pi_{s_{\overline{k}}} \cdot w_{\overline{k}} )' =  {\pi_{s_{\overline{k}}}} \cdot w_{\overline{k}} $.   Thus $w_{\overline{k}+1}$ is admissible in diagram $0$ and  is periodic of period $1$ or $2$.  One can check from Figure \ref{diagrams} that this implies that either   $ w_{\overline{k}+1} \in  \{\overline{AD}, \overline{BC}\}$ (Case $a$) or  $ w_{\overline{k}+1} \in \{\overline{DB}, \overline{C}\}$ (Case $b$). 
 Since $ w_{\overline{k}+1}$ is coherent with respect to $(s_{\overline{k}+1}, s_{\overline{k}+2} )$, it is  admissible in diagram $s_{\overline{k}+1}$ (by  $($C$0)$ in Definition \ref{coherencedef}).  
   In Case $a$, we see from Figure \ref{diagrams}  that $ w_{\overline{k}+1} $ is admissible in diagram $0$ and in diagram $1$.  Since  $s_{\overline{k}+1} = 0$ is excluded because $\{s_k \}_{k\in\mathbb{N}} \in S^*$, this shows that $s_{\overline{k}+1} = 1$. Moreover, in Case $a$ $w_{\overline{k}+1}$ is also a fixed point for the action of $\pi_1$, so that  $ w_{\overline{k}+2}= ({\pi_1} \cdot w_{\overline{k}+1})' =  w_{\overline{k}+1}   $. Repeating the same considerations by induction, we get that 
   $s_k=1$  for each $k > \overline{k}$. In Case $b$, one can see that $ w_{\overline{k}+1}$ is admissible in diagram $0$ and in diagram $7$ and, reasoning in the same way, conclude that 
  $s_k= 7$ for each $k > \overline{k}$.  This  concludes the proof of  $(ii)$.
\end{proofof}

\subsection{Approximation by periodic cutting sequences and closure}\label{proofsfindtrajectoiriessec}
In this section we first prove that infinitely coherent words can be approximated by periodic cutting sequences, giving the proofs of Lemma \ref{Pperiodic}, Proposition \ref{finiterealization} and Corollary \ref{approxbyperiodic}. We then give the proof of Theorem \ref{cuttingseqthm}. 
\begin{proofof}{Lemma}{Pperiodic}
Let $w =  \g{s_1}{s_0}{  \ldots   \g{s_k}{s_{k-1}}{u}}$, where $u\in \mathscr{P}_{s_k}$. By Remark \ref{Pdescription}, $u$ is admissible in either $\mathscr{D}_{s_k-1,s_k}$ or $\mathscr{D}_{s_k,s_{k+1}}$ in Figure \ref{boundarydiagrams} (where the indexes should be intended modulo $7$). Thus, one can find a periodic trajectory $\tau_{s_k}$ in direction $\pi s_k /8$ or $\pi (s_k+1)\pi/8$ such that $u =c(\tau_{s_k})$ (see  Remark \ref{periodicmultipleofpi8}). Let us recursively define, for $j=0, \dots k-1$,  $\tau_{s_{k-(j+1)}}:= \nu_{s_{k-(j+1)}}^{-1} \trajcdot \Psi_{\gamma}\tau_{s_{k-j}}$.  Using  Lemma \ref{gopgeometrically} and induction, one can show  that 
$  \g{s_{k-(j+1)}}{s_{k-j}}{  \ldots   \g{s_k}{s_{k-1}}{u}}$ is the cutting sequence of the trajectory $\tau_{s_{k-(j+1)}}$. For $j=k-1$,  this shows that $w $ is the cutting sequence of the periodic trajectory $ \tau_{s_0}$. Since moreover each map of the form  $\nu^{-1}\Psi_{\gamma}$ sends periodic trajectories to periodic trajectories, $\tau_0$ is a periodic trajectory. The proof is concluded by setting $\tau:= \tau_0$.  
\end{proofof}

\begin{proofof}{Proposition}{finiterealization}
Let $u$ be a finite subword of a word $w$ infinitely coherent  with respect to $\{s_k\}_{k\in \mathbb{N}} \in S^*$. Let $w_k$ be
 the words in Definition \ref{combrecursive}.  Let $u_k$ be the subword of $w_k$ (possibly empty) of letters of $u$ which survive as letters in $w_k$ (see Example \ref{examplereconstruction} below).  

If $u_{k}$ is not empty and not all letters of $u_{k}$   are sandwiched,  $u_{k+1}$ is strictly shorter than $u_{k}$ (since the action of permutations do not change the property of being sandwiched and non-sandwiched letters are erased by derivation). Thus, either  there exists a $\overline{k}$ such that all letters of   $u_{\overline{k}}$ are sandwiched and hence $u_{k}$ has the same number of letters as $u_{\overline{k}}$ for all $k\geq \overline{k}$ 
 (let us call this situation Case $(i)$), or  there  exists $\overline{k}$  such that $\overline{k}+1$ is the minimum  ${k}\in \mathbb{N}$ for which $u_{{k}}$ is empty (Case $(ii)$).

Let us show that in both cases  $u_{\overline{k}}$ is realized by some periodic trajectory $\tau_{\overline{k}}$. In Case $(i)$, since all letters of  $u_{\overline{k}}$ are sandwiched, either $u_{\overline{k}} = LLL \cdots L$ or $u_{\overline{k}} =L_1L_2L_1L_2 \cdots L_{i}$ ($i=1,2$). In this case, we define $\tau_{\overline{k}}$ to be  the periodic trajectory  realizing the transition $LL$ or $L_1L_2$, which has direction $l\pi/8 $ for some $0\leq l < 8$ (see Figure \ref{boundarydiagrams} and Remark \ref{periodicmultipleofpi8}).  In Case $(ii)$, consider the normal form  sequence ${\pi_{s_{\overline{k}}}} \cdot {w_{\overline{k}}} $ and let $L_1$ be the last sandwiched letter of  ${\pi_{s_{\overline{k}}}} \cdot {w_{\overline{k}}} $  before the occurrence of the subword ${ \pi_{s_{\overline{k}}} } \cdot{ u_{\overline{k}}}$ and $L_2$ be the first sandwiched letter after $ {\pi_{s_{\overline{k}}}} \cdot {u_{{\overline{k}}}}$ (see Example \ref{examplereconstruction} below). Since by construction there are no other sandwiched letters in between,    $w_{k+1}= ( \pi_{s_k} \cdot w_{{\overline{k}}})'$ contains the transition $L_1L_2$ and is admissible in diagram $s_{\overline{k}+1} \neq 0$ by definition of coherence. Let  $\overline{\tau}$ be a periodic trajectory admissible in diagram $s_{\overline{k}+1}$ realizing $L_1L_2$.  By Lemma \ref{coherencegeneration},  $ {\pi_{s_{\overline{k}}}} \cdot {w_{{\overline{k}}}} = \g{s_{\overline{k}+1}}{0} { w_{k+1} }$ and the subword of ${\pi_{s_{\overline{k}}}} \cdot { w_{{\overline{k}}}} $ between $L_1$ and $L_2$ is obtained by interpolating $L_1L_2$ by the  generation word $w^{s_{\overline{k}+1}}_{L_1L_2}$  on diagram $\mathscr{D}_{s_{\overline{k}+1}}$ in Figure \ref{generationrulesfig}. Consider the trajectory $\Psi_{\gamma} \trajcdot  \overline{\tau}$: by Proposition \ref{generatetype0}, $c( \Psi_{\gamma} \trajcdot  \overline{\tau})= \g{s_{\overline{k}+1}}{0}{c(\overline{\tau}) }$ and since $c(\overline{\tau})$ contains the transition $L_1L_2$, $c( \Psi_{\gamma} \trajcdot  \overline{\tau})$ contains $w^{s_{\overline{k}+1}}_{L_1L_2}$ and in particular the subword ${\pi_{s_{\overline{k}}}} \cdot {u_{{\overline{k}}}}$. In other words,  $ \Psi_{\gamma} \trajcdot  \overline{\tau}$  realizes ${\pi_{s_{\overline{k}}}}\cdot { u_{{\overline{k}}}} $  and,  by (\ref{relabellingeq}), if we define $\tau_{\overline{k}}:=\nu_{s_{\overline{k}}}^{-1}\Psi_{\gamma} \trajcdot  \overline{\tau} $, then $\tau_{\overline{k}}$ realizes  $u_{\overline{k}}$. 

Let us show by induction that if $\tau_{k}$ realizes  $u_k$ for some $k\geq 1$, then  $\tau_{k-1}:=  \nu_{s_{k-1}}^{-1}\Psi_{\gamma} \trajcdot \tau_k $ realizes  $u_{k-1}$. By Definition \ref{infcoherencedef}, $w_{k-1}$ is coherent with respect to $(s_{k-1}, s_k)$ and by Lemma \ref{coherencegeneration}, this means equivalently that $w_{k-1}= \g{s_{k}}{s_{k-1}}{v}$ for some $v$ admissible in diagram $s_k$. Reasoning as in the proof of Proposition \ref{generationandinfcoherence}, we can see that $v=w_{k}$ and $w_{k-1}= \g{s_{k}}{s_{k-1}}{w_k}$. This implies that $u_{k-1}=  \g{s_{k}}{s_k-1}{u_k}$. 
Make the  inductive assumption that  $u_{k}$ is realized by $\tau_k$, or in other words that $c(\tau_k)$ contains the finite word $u_k$. Let us remark that $w_k$ and thus $u_k$ are admissible on diagram $s_k\neq 0$ by definition of coherence. Consider $\Psi_{\gamma} \trajcdot \tau_k $. By  Proposition \ref{generatetype0}, $c( \Psi_{\gamma} \trajcdot \tau_k  ) = \g{s_k}{0}{c(\tau_k)}$, thus  $c( \Psi_{\gamma} \trajcdot \tau_k  )$ contains in particular the interpolated word $ \g{s_k}{0}{u_k}$. By (\ref{relabellingeq}), $\tau_{k-1}=  \nu_{s_{k-1}}^{-1}\Psi_{\gamma} \trajcdot \tau_k $ contains  $\pi_{k-1}^{-1} \cdot \g{s_k}{0}{u_k} = \g{s_{k}}{s_{k-1}}{u_k}$ (Definition \ref{othergoperators}) and since we showed above that $u_{k-1}=  \g{s_{k}}{s_{k-1}}{u_k}$, this shows that $\tau_{k-1}$ realizes $u_{k-1}$. For $k=1$, we get a trajectory $\tau:= \tau_0$ that realizes $u_0=u$. 

Let us now prove  (\ref{periodicunion}). Notice that, both in Case $(i)$ and in Case $(ii)$, the trajectory $\tau_{\overline{k}}$ is periodic of period $1$ or $2$ and admissible in diagram $s_{\overline{k}}$, thus $c( \tau_{\overline{k}}) \in \mathscr{P}_{s_{\overline{k}}}$. We showed that $u$ is a subword of the word $c(\tau)$. Since  $c(\tau)$ is obtained from $c(\tau_{\overline{k}}) \in \mathscr{P}_{s_{\overline{k}}}$ by applying the composition  of operators $\gop{s_1}{s_{0}} \dots \gop{s_{j}}{s_{j-1}}\dots \gop{s_{\overline{k}}}{s_{\overline{k}-1}} $,  this shows that $c(\tau)$ is in  $\mathscr{P}(s_0, \dots, s_{\overline{k}})$  and thus it belongs to the union in (\ref{periodicunion}), concluding the proof. 
\end{proofof}

\begin{ex}\label{examplereconstruction}
Consider the word $w$ in Example \ref{typesex} and let  $u=u_0= AAADBDBCBD$ be its subword. Then $u_{1}= CCDBAB$,  $u_{2}= A$ and $\overline{k}=2$. Moreover $\pi_2\cdot u_{2}= \B$ is a subword of $\pi_2 \cdot w_{2} =  n(w_{2})= \dots  CBCC\B DAD \cdots $ (where $\pi_2\cdot u_{2}$ is in bold font), hence the preceeding sandwiched letter is $L_1=B$ and the following sandwiched letter is    $L_2=A$.
\end{ex}

\begin{proofof}{Corollary}{approxbyperiodic}
Since for each arbitrarily long finite subword $u$ of $w$ the second part of Proposition \ref{finiterealization} gives an element of some $\mathscr{P}_k$ which contains $u$, the Corollary follows simply by  the definition of topology on $\mathscr{A}^\mathbb{Z}$ (see \cite{LM:sym}).   
\end{proofof}

\begin{proofof}{Theorem}{cuttingseqthm} 
Let $IC$   denote the set of infinitely coherent words in $\mathscr{A}^{\mathbb{Z}}$, $CS$ 
  the set of cutting sequences and  $\overline{CS}$ be its closure.   
In order to show that $\overline{CS} = IC$, one has to show that $CS \subset IC$, that $IC$ is closed and that $CS$ is dense in  $IC$.
Proposition \ref{infcoherentprop} states that cutting sequences are infinitely coherent with respect to their sequences of sectors  $\{s_k\}_{k\in\mathbb{Z}}$, which are in $S^*$ by Corollary \ref{typesequenceslemma}, thus the inclusion  $CS \subset IC$ holds. 
 By Proposition \ref{generationandinfcoherence},  $IC$ can be equivalently expressed as the countable intersection of finite unions  (\ref{unionfiniteintersection}).   
 Since the set $\mathscr{A}_k$ of admissible words in diagram $k$  is a subshift of finite type, $\mathscr{A}_k$ is closed (see for example Chapter $6$, \cite{LM:sym}). Moreover, one can check that the operators $\gop{j}{i}$ are Lipschiz, since if $u, v \in \mathscr{A_j}$  have a common subword, the interpolated words $\g{j}{i}{u}$ and  $\g{j}{i}{v}$ have an even longer common subword. Thus, the sets $\mathscr{G}(s_0, \dots, s_k)$ in (\ref{unionfiniteintersection}), being the image of a closed set under a continuous map, are closed. 
  Since, 
 in (\ref{unionfiniteintersection}) for each  $k$ one considers a finite union of closed sets, 
  $IC$  is a countable intersection of closed sets  and thus is closed.  
Finally, Corollary \ref{approxbyperiodic} combined with Lemma \ref{Pperiodic} show that cutting sequences of periodic trajectories  (thus in particular cutting sequences) are dense in  $IC$. 
\end{proofof}

\section{Regular $2n$-gons}\label{2ngonssec}
In this section we briefly sketch how our results for the regular octagon generalize to a  regular polygon with an even number of sides.  Let  $O_{2n}$ be a regular polygon with $2n$ sides, for $n\ge 3$. The case $n=2$, which corresponds to the square, was described in section \S\ref{torussec} and follows a somewhat different pattern from the other polygons. The sides of $O_{2n}$  can be identified by glueing pairs of opposite parallel sides. Let us denote the pairs, in clockwise order, by the letters $L_1, \dots, L_n$. We will denote this alphabet by $\mathscr{A}_n=\{ L_1, \dots, L_n\}$. Let $S_{2n}$ denote the surface obtained from $O_{2n}$ by making these identifications. Let $\Sing \subset S_{2n}$ be the set of points corresponding to vertices. The formulation of the results will differ slightly according to whether $n=2k$ is even or $n=2k+1$ is odd. If $n=2k$ is even then $S_{2n}$ is a surface of genus  $k$ where $\Sing$ consists of a single point with cone angle $2\pi(n-1)$. If $n=2k+1$ is odd then $S_{2n}$ is a surface of genus $k$ where $\Sing$ consists of two points both with cone angle $2\pi k=\pi(n-1)$. Note that when $n=3$ the surface $S_6$ is the torus and the set $\Sing$ consists of two non-singular points. When we refer to the Veech group in this case we mean the Veech group of $S_6$ consisting of toral automorphisms which take $\Sing$ to itself. This will be a subgroup of index two in the Veech group of the torus.

\subsection{Veech group and transition diagrams}
\paragraph{Isometries.}
The group $\Die_{2n}$  of isometries of $O_{2n}$ is generated by the reflection $\alpha_{n}$ in the horizontal axis  and the reflection $\beta_{n}$ in the slanted line forming and angle $\pi/{2n}$ with the horizontal. 
Let $\nu_0, \dots, \nu_{2n-1}$ be the elements of $\Die_{2n}$ given by\footnote{One can verify that if $n=4$ these elements coincide with the ones defined in (\ref{nujdef}).} $\nu_i = (\alpha\beta)^k $ if $i=2k$ and $\nu_i = (\beta \alpha)^k \beta$ if $i=2k+1$. The element $\nu_i$  brings the sector $\overline{\Sigma}^{(n)}_{i}:=[\pi i / 2n , \pi (i+1)/2n]$ to  $\overline{\Sigma}^{(n)}_{0}:=[0, \pi/2n]$. Let $\pi_0, \dots, \pi_{2n-1}$ be the corresponding induced permutations (see Definition \ref{inducedpermutation}) on the elements of $\mathscr{A}_n$.

\paragraph{Veech group.}
As in the case of the octagon, an additional element in the Veech group of $O_{2n}$ can be obtained acting by simultaneous Dehn twists in a cylinder decomposition. If $n=2k+1$, $O_{2n}$ can be cut and pasted  into $k$ horizontal cylinders, all of which have inverse modulus $\mu_n:=2\cot (\pi/2n) $.  If $n=2k$, $O_{2n}$ can be decomposed into $k$ cylinders, $k-1$ of modulus $\mu_n$ and one of modulus $\mu_n/2$.
One can obtain a globally well defined affine diffeomorphism $\Psi_{\sigma_n}$ whose  derivative is given by 
\bes 
D\Psi_{\sigma_n }   : = \sigma_n  : =  \begin{pmatrix} 1 &  2\cot(\pi/2n)  \\ 0 &  1  \end{pmatrix}. \ees  
This diffeomorphism is a shear which acts as a Dehn twist on the cylinders of modulus $\mu_n$ and acts as the square of a Dehn twist on the cylinder of modulus $\mu_n/2$.
Let $\gamma_n = \sigma_n \nu_{n}$ be the corresponding affine reflection with axes in the horizontal and $\pi/2n$ directions. Let $\Psi_{\gamma_n}$ be the corresponding affine automorphism, obtained by composing $\gamma_n$ with a cut and paste map $\Upsilon$ (which is again uniquely determined since one can prove as in Lemma \ref{trivialkernel} that the kernel of the Veech homomorphism is trivial).  Veech proved that $V^+(O_{2n})$ is generated by $\sigma_n$ and $\rho_{\pi/2n}$. Reasoning as in \S\ref{veechsec}, 
one can prove that $V(O_{2n})$ is generated by $\alpha_n, \beta_n, \gamma_n$.

\paragraph{Transitions diagrams.}
The structure of the transition diagram $\mathscr{D}_0^n$ (according to the parity) for a trajectory in $\overline{\Sigma}_0^n=[0, \pi/2n]$ is given by the diagram in Figure \ref{Dkoddeven}. 
The other transition diagrams $\mathscr{D}_1^n, \dots,\mathscr{D}_{2n-1}^n $ are obtained by permuting the letters by $\pi_{k}$, $k=1, \dots, 2n-1$. 
\begin{figure}[h]
\centering
\subfigure[The diagram $\mathscr{D}_0^n$ for $O_{2n}$ when $n=2k$ is even.\label{Dkevenfig}]{\includegraphics[width=0.9\textwidth]{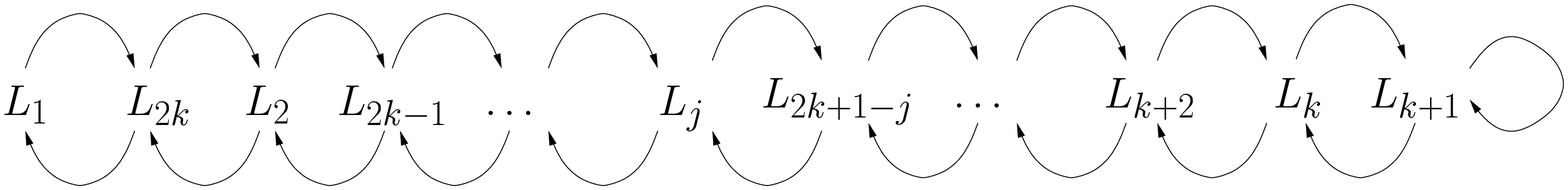}}\\
\hspace{9mm}

\subfigure[The diagram $\mathscr{D}_0^n$ for $O_{2n}$ when $n=2k+1$ is odd.\label{Dkoddfig}]{\includegraphics[width=0.9\textwidth]{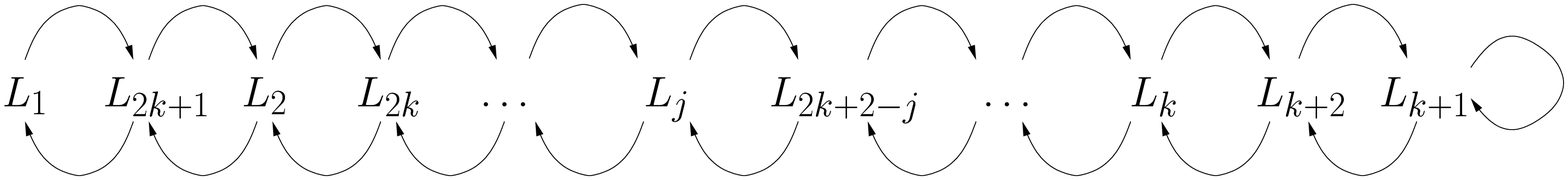}}
\hspace{6mm}
\caption{The basic transition diagram for a regular $2n$-gon, $n=2k$ or $n=2k+1$.\label{Dkoddeven}}
\end{figure}

\subsection{Derived sequences and $2n$-gon Farey map} 
Let us say that a sequence is \emph{admissible} if it corresponds to an infinite path on one of the graphs $\mathscr{D}_0^n, \dots, \mathscr{D}_{2n-1}^n$.  We derive a sequence as we did before by deleting all letters which are not sandwiched. The necessary condition for cutting sequences in terms of derivability in Theorem \ref{cutseqinfderiv} holds for any $2n$-gon.
 
\begin{thm} An infinite cutting sequence in a regular 2n-gon  is infinitely derivable.
\end{thm}

Similarly, the direction recognition Theorem \ref{directionsthm} extends to $O_{2n}$, where the sequence of diagrams is now a sequence in $\{0, \dots, 2n-1\}^{\mathbb{N}}$ and the corresponding $2n$-gon Farey map is a piecewise expanding map with $2n$ branches given by 
\begin{equation*}
 F^{(n)}(\theta)  := \cot^{-1} \left( (\gamma_n \nu_i) \left[  \cot \theta \right] \right)  \quad \mathrm{if} \,  \theta \in \overline{\Sigma}^n_i, \quad i=0, \dots, 2n-1. 
\end{equation*} 

We have an additive continued fraction algorithm as before. Any $\theta\in[0,\pi]$ corresponds to an infinite sequence $[s_0; s_1,\ldots ]$ with $s_0\in\{0,\ldots , n-1\}$ and  $s_k\in\{1, \ldots, n-1\}$ for $k>0$, which is given by its Farey itinerary.

\begin{thm} An infinite non-periodic cutting sequence determines its direction uniquely via a Farey itinerary.
\end{thm}

Let us say that a direction is \emph{terminating} if its itinerary with respect to the $2n$-gon Farey-map is eventually $1$ or eventually $n-1$. One can show (as in Proposition \ref{3periodic}) that a cutting sequence $c(\tau)$ if periodic if and only if the direction of $\tau$ is terminating.

\begin{rem}\label{holonomy} The holonomy field or trace field for $S_{2n}$ is the field $\Q[\cos{(\pi/n)}]$ (see \cite{HS:inf}). In general the holonomy field does not necessarily correspond to the slopes of terminating directions (as in the octagon, see Theorem \ref{Calta}) if $n$ is different from 3 or 4 (see \cite{Leu:ube}, \cite{Leu:ube2}).
\end{rem}

\subsubsection{Complexity of cutting sequences}\label{wordcomplexitysec} Let us recall the argument which shows that all the words that we are considering have linear complexity.
\begin{prop}\label{growth} 
If we fixed a  direction, the number of words of length $\ell$ appearing as subwords of bi-infinite cutting sequences in that direction is bounded above by $(n-1)\ell+1$. If the direction is non-terminating,  then the number of words of length $\ell$ that appear is $(n-1)\ell+1$. 
\end{prop}
\begin{proof} This follows from the fact that we can realize the flow in a fixed direction as a special flow over an interval exchange transformation on $n$ intervals, so that cutting sequences coincide with symbolic coding of orbits of the interval exchange transformation with respect to the natural partition into subintervals. The number of words of length $\ell$ is given by the number of intervals in the $\ell$-th iterate of this interval exchange transformation and this is bounded above by $(n-1)\ell+1$ (see \cite {CFS:erg}). If the direction is non-terminating, cutting sequence are not periodic and the Veech dichotomy implies that there are  no saddle connections. This means that the interval exchange transformation satisfies the Keane condition and in this case we have equality (see again \cite {CFS:erg}).
\end{proof}

\subsection{Coherence and closure}
The characterization of cutting sequences can also be generalized to $2n$-gons. 
Let  $\pi_n^{\rightarrow}$ be the permutation of $\mathscr{A}_n$ induced  by the vertical reflection (see Definition \ref{inducedpermutation}) and $\pi_n^{\uparrow}$ be the one induced by a reflection in an axis forming a angle of $\pi/2n$ with respect to the vertical (i.e.~an angle $\pi/2 + \pi/2n$ with respect to the horizontal). We have
\begin{equation}\label{pi_top_right_def}
\begin{array}{ll}
\pi_n^{\rightarrow}  = (L_2 L_{n})(L_3 L_{n-1}) \dots ( L_j L_{n+2-j}) \dots (L_{k} L_{k+2})(L_{k+1}) , \qquad & \mathrm{if } \quad n=2k; \\
 \pi_n^{\rightarrow}  = (L_2 L_{n})(L_3 L_{n-1})   \dots ( L_j L_{n+2-j}) \dots (L_{k+1} L_{k+2}), \qquad & \mathrm{if }\quad n=2k+1;  \\
 
 \pi_n^{\uparrow}=  (L_1 L_{n}) (L_2 L_{n-1}) \dots (L_j L_{n+1-j})\dots  ( L_{k+1}L_k),   \qquad & \mathrm{if}\quad n=2k; \\
 \pi_n^{\uparrow}= (L_1 L_{n}) (L_2 L_{n-1})\dots (L_j L_{n+1-j}) \dots (L_{k} L_{k+2})(L_{k+1}), \qquad & \mathrm{if }\quad   n=2k+1;
 \end{array}
 \end{equation}

\paragraph{Coherence}
\noindent The coherence Definition \ref{coherencedef} is generalized as follows (see also \S\ref{oncoherencesec} below  for more insight on this Definition).
\begin{defn}[Coherence for $O_{2n}$]\label{coherencedef2ngons} A word $w \in \mathscr{A}_n^{\mathbb{Z}}$ is coherent if there exists a pair $(i,j)$ of indexes  in $\{0,\dots, 2n-1\}$ such that :  
\begin{itemize}
\item[(C0)$'$] The word $w$ is admissible in diagram $i$;
\end{itemize}
If we normalize $w$ by setting $n(w):=\pi_i \cdot w$, then:
\begin{itemize}
\item[(C1)$'$] The sandwiched letters which occur in $n(w)$ fall into one of the groups  $G^n_0, \dots, G^n_{n-1}$, where
\bes  G_l:=\left\{ \begin{array}{ll} L_j \ \mathrm{ is}\   \pi_n^{\uparrow}(L_j)-\mathrm{sandwiched}\  \mathrm{for}\ 1\leq j \leq n-l, \\   L_j\ \mathrm{ is}\  \pi_n^{\rightarrow}(L_j)-\mathrm{sandwiched} \ \mathrm{ for}\  n-l < l \leq n \end{array} \right\};
 \ees 
\item[(C2)$'$] The derived sequence  $n(w)'$ is admissible in a diagram $j \in \{1, \dots, 7\}$;
\item[(C3)$'$] The indices $k$ and $j$ (defined in (C1)$'$ and (C2)$'$ respectively)  are related by the formula $k=[j/2]$. 
\end{itemize}
\end{defn}
The Definition of infinite coherence is then analogous to Definition \ref{infcoherencedef}.
Using this definition of coherence,   Proposition \ref{infcoherentprop} on infinite coherence of cutting sequences holds and one can give a characterization of the closure of cutting sequences (see Theorem \ref{closure2ngons} below).


\paragraph{Generation rules.}
In view of the structure of the diagram $\mathscr{D}_0^n$, one can write two Tables (according to the parity of $n$) analogous to the Table \ref{transitionsinterpolations}, whose entries $L_iL_i$ and $L_iL_j$ for $i<j$ have the form shown in Tables \ref{interpolations2neven} and \ref{interpolations2nodd} (compare with Figure \ref{Dkoddeven}, each entry corresponds to a path which does not contain any sandwiched letter).
\begin{table}
 \begin{tabular}{c ||l  }
$ \LL_i \LL_i $ &	 $\LL_i L_{2k+1-i}L_{i+1} L_{2k-i}\dots  L_{k+2}L_k  L_{k+1} L_{k+1}L_k  L_{k+2} \dots  L_{2k+1-i} \LL_i  $ \\
\hline
\multirow{2}{*}{ $\LL_i \LL_j$ }	 &	 $\LL_i L_{2k+1-i}L_{i+1} L_{2k-i} \dots  \LL_j $	 										 \\
				                         &	 $ \LL_i L_{2k+1-i}L_{i+1} L_{2k-i}\dots   L_{k+2}L_k  L_{k+1} L_{k+1}L_k  L_{k+2} \dots L_{2k+1-j} \LL_j$ \\  
\end{tabular}
\caption{Interpolations for  $n=2k$, $1\leq i <j \leq n$.\label{interpolations2neven}}
\end{table}
\begin{table}
\begin{tabular}{c||l  }
{ $\LL_i \LL_i$ } 
                     &	$  \LL_i L_{2k+2-i}L_{i+1} L_{2k+1-i}\dots  L_k L_{k+2}  L_{k+1} L_{k+1}  L_{k+2} L_k \dots  L_{2k+2-i} \LL_i $ \\
\hline
\multirow{2}{*}{ $\LL_i \LL_j$ }	 
				                                           
				&    $\LL_i L_{2k+2-i}L_{i+1}L_{2k+1-i} \dots  \LL_j   \dots  \LL_j $   \\
& $\LL_i L_{2k+2-i}L_{i+1} L_{2k+1-i}\dots  L_k L_{k+2}  L_{k+1} L_{k+1}  L_{k+2} L_k \dots  L_{2k+2-j} \LL_j$ \\
\end{tabular}
\caption{Interpolations for  $n=2k+1$, $1\leq i <j \leq n$.\label{interpolations2nodd}}
\end{table}
Generation rules which generalize the ones in Figure \ref{generationrulesfig} can thus be obtained by labelling the graphs $\mathscr{D}_1^n, \dots, \mathscr{D}^n_{2n-1}$ with the interpolation words given in Tables \ref{interpolations2neven} and \ref{interpolations2nodd}, where the interpolation word is chosen from the multiple entries in the Tables as follows: the generation words which appear in $\mathscr{D}_k^{n}$, $1\leq k\leq 2n-1$ are such that all letters are sandwiched according to group $G_{[k/2]}$. This determines one interpolation word  in each entry of  Tables \ref{interpolations2neven} and \ref{interpolations2nodd}. These allow us to define generation operators $\gop{i}{j}$, $0 \leq i,j \leq 2n-1$, as in Definition \ref{generationopto0} and \ref{othergoperators}.
Using these generation operators to get the analogue of Proposition \ref{generatetype0}, one can prove a characterization of finite words in a cutting sequence and prove the following analogue of Theorem \ref{cuttingseqthm} and Proposition \ref{generationandinfcoherence} for regular $2n$-gons. 

\begin{thm}\label{closure2ngons} The closure of the set of cutting sequences for the regular $2n$-gon   coincides with the set of infinitely coherent words. Equivalently, the closure is described in terms of generation operators as 
\bes
 \bigcap_{k\in\mathbb{N}} \, \,   \bigcup_{\substack{ { 0\leq s_0 \leq 2n-1 } \\1\leq s_1, \dots, s_k   \leq 2n-1 }} \left\{ \,  \g{s_1}{s_0}{ \g{s_2}{s_1}{ \ldots \g{s_{k-1}}{s_{k-2}}{ \g{s_k}{s_{k-1}}{u}}}},  \,   \quad u \in \mathscr{A}_{s_k}  \,  \right\} . 
\ees
where $\mathscr{A}_{s}$ are the words admissible in $\mathscr{D}_s^{n}$, $0\leq s \leq 2n-1$.
 \end{thm}

\subsubsection{On the definition of coherence}\label{oncoherencesec}
Let us add a few comments to justify the definition $($C$1)'$ and sketch how to reproduce  the proof of Proposition \ref{infcoherentprop} for $O_{2n}$.
Consider the affine polygon $O_{2n}':=\Psi_{\gamma_n} \octcdot O_{2n}$. If we apply a shear which fixes horizontal lines and sends lines of slope $\pi/2n$ to vertical lines, the affine polygon can be represented as a staircase of rectangles, as shown in Figure \ref{staircase}, where opposite parallel sides are identified and the staircase contains exactly two copies of  $O_{2n}$. 
\begin{figure}[!h]
\centering
\includegraphics[width=.7\textwidth]{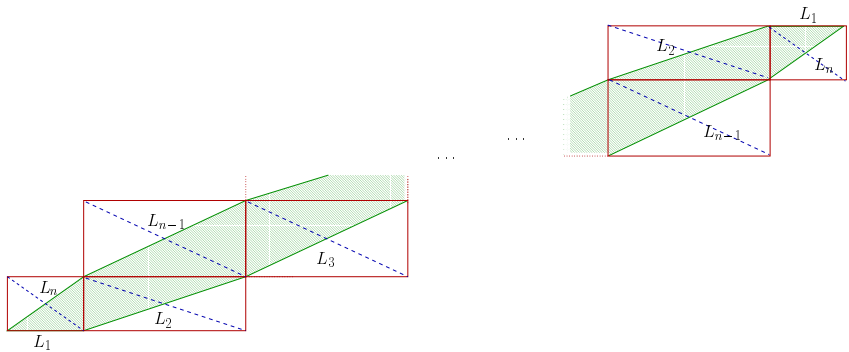}
\caption{The staircase representation of an affine $2n$-gon.\label{staircase}}
\end{figure}
While the bottom left-top right diagonals of each rectangle  are sides of  $O_{2n}'$, the opposite ones are sides of $O_{2n}$. We remark that the letters labelling the  two diagonals of each rectangle are same, one being a letter of the original cutting sequence, the other being a letter of the derived sequence.

The permutations  $\pi_n^{\rightarrow}$ and  $\pi_n^{\uparrow}$ 
contain the following information\footnote{These permutations are closely related to the ones used to describe glueing patterns in \cite{EO:asy} or \cite{SmW:glueing}.}. 
If the rectangle $R_j$ has a diagonal labelled by $L_j$, then the rectangle vertically adjacent to it has diagonals labelled $\pi_{\uparrow}(L_j)$, while the  diagonals of the rectangle horizontally next to it are labelled by $\pi_{\rightarrow}(L_j)$. 

Hence, denoting by $\theta'_j$ the angle of the side of $O_{2n}'$ labelled by $L_j$, reasoning as in the proof of Proposition \ref{infcoherentprop} in \S \ref{infcoherencesec}, one can see that if $\theta'> \theta'_j$ then the  letter $L_j$ is  $\pi_{\uparrow}(L_j)$-sandwiched while if $\theta'< \theta'_j$ then $L_j$ is  $\pi_{\rightarrow}(L_j)$-sandwiched.  As in \S \ref{infcoherencesec} (see Table \ref{coherencetable}), the group $G^n_l$ lists the sandwiched letters  corresponding to $\theta'_{n-l}< \theta' < \theta'_{n-l+1}$ (setting $\theta'_{n+1}:=\pi/2n$) and this explains the definition of groups given in $($C$1)'$.

\section*{Acknowledgments}
We would like to acknowledge the hospitality and the support given from  the Centro di Ricerca Matematica Ennio de Giorgi in Pisa during the intensive research period on ``Dynamical Systems and Number Theory'' and from the Mathematical Sciences Research Institute in Berkeley, CA  during the Program ``Teichm{\"u}ller Theory and Kleinian Group'', which gave us the opportunity of beginning  the present collaboration. 

\bibliographystyle{alpha}

\begin{thebibliography}{McM03}

\bibitem[AF91]{AF:mod}
Roy~L. Adler and Leopold Flatto.
\newblock Geodesic flows, interval maps, and symbolic dynamics.
\newblock {\em Bull. Amer. Math. Soc}, 25:229--334, 1991.

\bibitem[AH00]{AH:fra}
Pierre Arnoux and Pascal Hubert.
\newblock Fractions continues sur les surfaces de veech.
\newblock {\em Journal d'Analyse Math{\'e}mathique}, 81:35--64, 2000.

\bibitem[Arn02]{Ar:stur}
Pierre Arnoux.
\newblock {\em Sturmian sequences}, pages 143--198.
\newblock Lecture Notes in Math., 1794. Springer, Berlin, 2002.

\bibitem[Cal04]{C:vee}
Kariane Calta.
\newblock Veech surfaces and complete periodicity in genus two.
\newblock {\em J. Amer. Math. Soc.}, 17(4):871--908, 2004.

\bibitem[CFS80]{CFS:erg}
I.~P. Cornfeld, S.~V. Fomin, and Ya.~G. Sinai.
\newblock {\em Ergodic Theory}.
\newblock Springer-Verlag, 1980.

\bibitem[CH73]{CH:seq}
Ethan~M. Coven and G.~A. Hedlund.
\newblock Sequences with minimal block growth.
\newblock {\em Math. Systems Theory}, 7:138--153, 1973.

\bibitem[Chr75]{C:obs}
E.~B. Christoffel.
\newblock Observato arithmetica.
\newblock {\em Annali di Mathematia 2nd. series}, 6:148--152, 1875.

\bibitem[EG97]{CG:tei}
Clifford~J. Earle and Frederick~P. Gardiner.
\newblock {\em Teichm{\"u}ller disks and Veech's $\mathcal{F}$-structures.},
  pages 165--189.
\newblock Number 201 in Contemp. Math. Amer. Math. Soc., Providence, RI, 1997.

\bibitem[EO01]{EO:asy}
Alex Eskin and Andrei Okounkov.
\newblock Asymptotics of numbers of branched coverings of a torus and volumes
  of moduli spaces of holomorphic differentials.
\newblock {\em Invent. Math.}, 145(1):59--103, 2001.

\bibitem[HS01]{HS:inv}
Pascal Hubert and Thomas~A. Schmidt.
\newblock Invariants of translation surfaces.
\newblock {\em Ann. Inst. Fourier}, 51(2):461--495, 2001.

\bibitem[HS04]{HS:inf}
Pascal Hubert and Thomas~A. Schmidt.
\newblock Infinitely generated {V}eech groups.
\newblock {\em Duke Journal}, 123(1):49--69, 2004.

\bibitem[Leu67]{Leu:ube}
Armin Leutbecher.
\newblock {\"{U}}ber die {H}eckeschen gruppen ${G}(\lambda )$.
\newblock {\em Abh. Math. Sem. Univ. Hamburg}, 31:31199--205, 1967.
\newblock (German).

\bibitem[Leu74]{Leu:ube2}
Armin Leutbecher.
\newblock {\"{U}}ber die {H}eckeschen gruppen ${ G}(\lambda )$. {II}.
\newblock {\em Mathematische Annalen}, 211:63--86, 1974.
\newblock (German).

\bibitem[LM95]{LM:sym}
Douglas Lind and Brian Marcus.
\newblock {\em Symbolic Dynamics and Coding.}
\newblock Cambdridge University Press, 1995.

\bibitem[McM03]{Mc:tei}
Curtis~T. McMullen.
\newblock Teichm{\"{u}}ller geodesics of infinite complexity.
\newblock {\em Acta Mathematica}, 191(2):191--223, 2003.

\bibitem[MH40]{MH:sym}
Marston Morse and Gustav~A. Hedlund.
\newblock Symbolic dynamics {II}. {S}turmian trajectories.
\newblock {\em Amer. J. Math.}, 62:1--42, 1940.

\bibitem[Rau79]{Ra:ech}
G\'erard Rauzy.
\newblock \'{E}changes d'intervalles et trasformations induites.
\newblock {\em Acta Arithmetica}, XXXIV:315--328, 1979.

\bibitem[Rau96]{Ra:low}
G\'erard Rauzy.
\newblock {\em Dynamics of complex interacting systems (Santiago, 1994)}, pages
  147--177.
\newblock Low complexity and geometry. Kluwer Acad. Publ., Dordrecht, 1996.

\bibitem[Ser85a]{Se:geo}
Caroline Series.
\newblock The geometry of {M}arkoff numbers.
\newblock {\em The Mathematical Intelligencer}, 7(3):20--29, 1985.

\bibitem[Ser85b]{Se:mod}
Caroline Series.
\newblock The modular surface and continued fractions.
\newblock {\em Journal of the London Mathematical Society. (2)}, 31(1):69--80,
  1985.

\bibitem[Ser91]{Se:sym}
Caroline Series.
\newblock {\em Ergodic theory, symbolic dynamics, and hyperbolic spaces
  (Trieste, 1989)}, chapter Geometrical methods of symbolic coding., pages
  125--151.
\newblock Oxford Sci. Publ. Oxford Univ. Press, New York, 1991.

\bibitem[Smi77]{S:note}
H.~J.~S. Smith.
\newblock Note on continued fractions.
\newblock {\em Messenger of Mathematics 2nd. series}, 6:1--14, 1877.

\bibitem[SW]{SmW:glueing}
John Smillie and Barak Weiss.
\newblock Characterizations of lattice surfaces.
\newblock Preprint Arxiv: 0809.3729.

\bibitem[Vee82]{Ve:gau}
William~A. Veech.
\newblock Gauss measures for transformations on the space of interval exchange
  maps.
\newblock {\em Annals of Mathematics}, 115:201--242, 1982.

\bibitem[Vee89]{Ve:tei}
William~A. Veech.
\newblock Teichm{\"u}ller curves in moduli space, eisenstein series and an
  application to triangular billiards.
\newblock {\em Invent. Math.}, 97(3):553--583, 1989.

\bibitem[Vor96]{Vor:pla}
Ya.~B. Vorobets.
\newblock Plane structures and billiards in rational polygons: the veech
  alternative.
\newblock {\em Uspekhi Mat. Nauk}, 51(5):3--43, 1996.
\newblock (Translated in: \emph{Russian Math. Surveys}, 51:5:779--817, 1996).

\end{thebibliography}

{\small \rmfamily{  DEPARTMENT OF MATHEMATICS, CORNELL UNIVERSITY, ITHACA, NY, 14853, USA}}

{\it E-mail address: }\texttt{smillie@math.cornell.edu} 

\vspace{5mm}

{\small \rmfamily{SCHOOL OF MATHEMATICS, UNIVERSITY OF BRISTOL, BRISTOL, BS8 1TW, UK }}

{\it E-mail address: }\texttt{corinna.ulcigrai@bristol.ac.uk}

\end{document}